\documentclass[12pt,twoside]{article}
\usepackage{amssymb,amsmath,amscd,amstext,amsthm,amsfonts}
\usepackage[ansinew]{inputenc} 
\usepackage{graphicx}
\usepackage[mathscr]{eucal}
\usepackage[pdfnewwindow=true]{hyperref}

\newcommand{\R}{\mathbb{R}}

\newcommand{\C}{\mathbb{C}}

\newcommand{\N}{\mathbb{N}}
\newcommand{\Z}{\mathbb{Z}}

\newcommand{\Su}{\mathbb{S}}

\newcommand{\hfrak}{\mathfrak{h}}
\newcommand{\gfrak}{\mathfrak{g}}

\newcommand{\bmat}{\left( \begin{array}}
\newcommand{\emat}{\end{array}\right) }

\newcommand{\Ascr}{\mathscr{A}}

\newcommand{\Escr}{\mathscr{E}}

\newcommand{\Mscr}{\mathscr{M}}
\newcommand{\Kscr}{\mathscr{K}}

\newcommand{\Fscr}{\mathscr{F}}

\newcommand{\Pscr}{\mathscr{P}}

\newcommand{\Sscr}{\mathscr{S}}

\newcommand{\Grass}{\mathscr{G}}
\newcommand{\Sperm}{{\rm S}}
\newcommand{\SJ}{{\rm S}^{sp}}
\newcommand{\Flag}{\Fscr}
\newcommand{\Fliso}{\Fscr^{sp}}

\newcommand{\Diff}{{\rm Diff}}

\newcommand{\tr}{{\rm tr}\,}

\newcommand{\Un}{{\rm U}}

\newcommand{\Mat}{{\rm Mat}}

\newcommand{\spl}{{\rm sp}}
\newcommand{\Sp}{{\rm Sp}}
\newcommand{\SL}{{\rm SL}}
\newcommand{\SO}{{\rm SO}}
\newcommand{\slin}{{\rm sl}}
\newcommand{\GL}{{\rm GL}}
\newcommand{\ut}{{\rm ut}}
\newcommand{\UT}{{\rm UT}}
\newcommand{\UTpos}{{\rm UT}_+}
\newcommand{\gl}{{\rm gl}}

\newcommand{\FlStr}[1]{ {\rm S}_{\underline{#1}} }
\newcommand{\OOStr}[1]{ \widetilde{{\rm S}}_{\underline{#1}} }
\newcommand{\FlisoStr}[1]{ {\rm S}^{sp}_{\underline{#1}} }
\newcommand{\OspStr}[1]{ \widetilde{{\rm S}}^{sp}_{\underline{#1}} }

\newcommand{\supp}{{\rm supp}}
\newcommand{\spsupp}{\widehat{{\rm supp}}}

\newcommand{\OO}{{\rm O}}

\newcommand{\diag}{{\rm diag}}
\newcommand{\Diag}{{\rm Diag}}

\newcommand{\Pm}{{\rm D}}
\newcommand{\OD}{{\rm D}}

\newcommand{\Osp}{{\rm O}^{sp}}

\newcommand{\card}{\mbox{card}}
\newcommand{\depth}{\mbox{depth}}
\newcommand{\codepth}{ \mbox{co-depth} }

\newcommand{\pol}[3]{p^{{#1},{#2}}_{#3} }
\newcommand{\polsp}[3]{{}^{sp}\! p^{{#1},{#2}}_{#3} }

\newcommand{\matriz}[4]{ 
\left(
\begin{array}{cc}
{#1}&{#2}\\
{#3}&{#4}
\end{array}\right) }

\newcommand{\coluna}[2]{
\left(
\begin{array}{c}
{#1}\\
{#2}
\end{array}\right)
}

\newcommand{\Mod}[1]{\left\vert{#1}\right\vert}
\newcommand{\nrm}[1]{\left\|#1\right\|}

\newcommand{\figura}[2]{
  \begin{center}
 \includegraphics*[scale={#2}]{{#1}} 
 \end{center}
}

\newtheorem{defi}{Definition}

\newtheorem{lema}{Lemma}
\newtheorem{rmk}{Remark}
\newtheorem{prop}{Proposition}
\newtheorem{teor}{Theorem}
\newtheorem{coro}{Corollary}

\newcommand{\dem}{ \par\medbreak\noindent{\bf
Proof. }\enspace} 
 
\newcommand{\cqd}{\hfill
$\sqcup\!\!\!\!\sqcap\bigskip$}

\newcommand{\coment}[1]{}

{

\title{A Class of Morse Functions on Flag Manifolds}
\author{ P. Duarte}
\date{}

\begin{document}
\maketitle

\noindent
{\bf Abstract }
{\em Given a positive definite symmetric matrix in one of the groups $\SL(n,\R)$ or $\Sp(n,\R)$,
we analyse its actions on Flag Manifolds, 
proving these are diffeomorphisms which admit  as strict Lyapunov functions a special class of quadratic
functions, all of them  $\Z_2$-perfect Morse functions.
Stratifications on the  Flag Manifolds are provided which are invariant under 
both the diffeomorphisms  and  the gradient flows of the Lyapunov functions.
}

\bigskip

\section{Introduction}

Let $G\subseteq \SL(n,\R)$ be a subgroup, \, $(X,\mu)$ be a probability space,
$T:X\to X$ a measurable transformation which preserves measure $\mu$, and  $A:X\to G$ a measurable function.
These objects, $(X,\mu)$, $T$ and $A$, determine a {\em linear cocycle}, a name which is often given to the skew-product map
$F:X\times\R^n\to X\times \R^n$,  defined by $F(x,v)=(T x, A(x)\,v)$,
whose iterates are given by
$F^k(x,v)=(T^k x, A^k(x)\,v)$,
where $A^k(x)= A(T^{k-1}x)\, \cdots \, A(T x)\,A(x)$.
This  cocycle is said to be {\em integrable} when
$\int_X \log^+ \nrm{A(x)}\, d\mu(x) <\infty$, where
where \, $\log^+(x) = \max\{\log x, 0\}$.
A classical theorem of Oseledet states that if the cocycle is integrable
and we assume, for the sake of simplicity, that $T$ is ergodic w.r.t. $\mu$, 
then there are $n\in\N$,  $\lambda_1> \lambda_2> \cdots > \lambda_{m}$ real numbers, 
and $\R^n = E_1(x) \supset E_2(x) \supset \ldots \supset E_{m}(x) $ an $F$-invariant measurable filtration such that:
for $\mu$-almost every $x\in X$ \, (1)  $\lambda_i = \lim_{k\to\infty}\frac{1}{k}\,\log \nrm{ A^k(x)\, v}$, \,
for every $v\in E_i(x)-E_{i+1}(x)$, and \,
(2) $\sum_{i=1}^{m} \lambda_i\,(\dim E_i(x)-\dim E_{i+1}(x) ) = 0$.
The numbers $\lambda_1,\ldots,\lambda_m$ are called the Lyapunov exponents of the linear cocycle.
The computation of Lyapunov exponents is an active and important subject in the theory of
Linear Cocycles and Dynamical Systems. Although there are some formulas for the Lyapunov exponents,
in general they can not be explicitely evaluated.
One of the first formulas for the Lyapunov exponents appeared in the seminal work of
H. Furstenberg~\cite{F2}. There the author is concerned with proving "law of large number" theorems for  products of 
independent and identically distributed ramdom matrices, but his setting can be brought 
to the linear cocycles' framework by letting the probability space to be $(G^\N,\mu^ \N)$, 
where $\mu$ is some  given probability measure on $G$, $T:G^\N\to G^\N$,
$T(g_k)_{k\in\N}= (g_{k+1})_{k\in\N}$, be the shift transformation, and 
$A:G^\N\to G$ be the function $A(g_k)_{k\in\N}= g_0$.
In the language of the Linear Cocycle's Theory, Furstenberg gives a 
formula  for the largest Lyapunov exponent $\lambda_1$, and proves that $\lambda_1>0$ under very
general assumptions.
There are several extensions and generalizations of Furstenberg's theorems,
first by Y. Guivarch'h, A. Raugi and also Ta. Gol'dsheid, G. Margulis still in the context of
Random matrices, and later by A. Avila, C. Bonatti, and M. Viana for generic base maps.

We now make a short overview of the basic concepts behind Furstenberg's
formula and theorems. In~\cite{F2} the author considers a semisimple real Lie group  $G$ with finite center.
Two examples of such groups are (1) $G=\SL(n,\R)$, the special linear group
of square $n\times n$ matrices with determinant one,
and (2) $G=\Sp(n,\R)$, the symplectic linear group of square $(2n)\times (2n)$ symplectic matrices.
Under these assumptions the group $G$ has three subgoups:
$K$  maximal compact, $A$ abelian and $N$  nilpotent  such that $G=K\cdot A\cdot N$. 
Moreover, for each $g\in G$, there is a unique decomposition
$g=k\,a\,n$ with $k\in K$, $a\in A$ and $n\in N$.
The decomposition $G=K\cdot A\cdot N$ is called Iwasawa's decomposition, but it is not unique. 
The group $S=A\cdot N$ is solvable. On the examples above, 
the group factorization $G=K\cdot S$  comes from the QR-decomposition, which is
associated with the Gram-Schmidt orthogonalization algorithm.
Let us now identify each of these factor subgroups on the two previous examples.
On example (1),  $K=\SO(n,\R)$ is the special orthogonal group,
$A$  is the group of determinant one positive diagonal matrices,  $N$ is the group
of upper triangular matrices with  ones on the diagonal, while 
$S$ is the group of upper triangular matrices with positive diagonal and determinant one.
On example (2),  $K=\Osp(n,\R)$, the group of symplectic orthogonal matrices, which is isomorphic to $\Un(n,\C)$,
and whence called the unitary group,
$A$  is the group of symplectic positive diagonal matrices,  $N$ is the group
of symplectic matrices of the form
$\bmat{cc} u & v\\ 0 & u^{-T} \emat $ where $u$ and $v$ are square $n\times n$ matrices,
$u$ is upper-triangular with  ones on the diagonal, and $v$ is
such that $v^T\,u^{-T}$ is symmetric. We shall refer to these matrices as symplectic upper-triangular.
Finally, in this case $S$ is the group of symplectic upper-triangular matrices
with positive diagonal.
For a compact manifold $M$, $\Pscr(M)$ will denote the space of probability measures on $M$,
which is a convex compact set when 
equipped with the weak-$\ast$ topology.
A key concept is that of boundary of a Lie Group.
A  boundary  of $G$ is any compact manifold with transitive
action $G\times M\to M$ with the following property:
for every measure $\pi\in\Pscr(M)$ there is a sequence $g_k\in G$ such that
$g_k\,\pi$ converges weakly to a point mass $\delta_p$, with $p\in M$.
The following relation is defined among boundaries of $G$.
Given two boundaries $M$ and $M'$ of $G$, we write
$M\preceq M'$ \, iff \, $M$ is the surjective image of $M'$ by
some $G$-equivariant map $f:M'\to M$, where
a  map $f:M\to M'$ is said to be  $G$-equivariant  if
$f(g\,x)=g\,f(x)$, for every $x\in M$ and $g\in G$.
The relation $\preceq$ is a partial order on the set of 
 $G$-equivariant equivalence classes
of boundaries of $G$. The poset of  
 $G$-equivariant equivalence classes of boundaries of $G$ has a unique maximal element,
 whose representatives are referred as maximal boundaries of $M$.
The quotient  $G/K\simeq S$ is a non-compact symmetric space with a transitive action of $G$.
The "boundary" name is justified, see~\cite{F2} page 403, by the fact that 
every boundary of $G$ can be realized as part of the boundary of some compactification
of the symmetric space $G/K$.
A class of  boundaries for the group $G=\SL(n,\R)$ are  the so called flag manifolds that
we now describe. See section 5 of~\cite{F1}.
Given a sequence of numbers $\underline{k}=(k_1,k_2,\ldots,k_m)$ with 
$1\leq k_1< k_2 <\ldots < k_m\leq n$, any sequence $V_\ast = (V_1, \ldots , V_m)$ of linear spaces
$V_1\subset V_2\subset \ldots \subset V_m\subset \R^n$ such that
$\dim V_i=k_i$ for every $i=1,\ldots, m$ is called a {\em $\underline{k}$-flag} in $\R^n$,
and $\Flag_{n,\underline{k}} = \Flag_{\underline{k}}(\R^n)$ will denote the {\em $\underline{k}$-flag manifold}
consisting of all $\underline{k}$-flags  in $\R^n$. There is a natural action of $\SL(n,\R)$ on 
$\Flag_{n,\underline{k}}$ defined by
$\SL(n,\R)\times \Flag_{n,\underline{k}}\to \Flag_{n,\underline{k}}$, 
$A\,V_\ast=(A V_1,\ldots, A V_m)$ where $V_\ast=(V_1,\ldots, V_m)$.
Given $A\in\SL(n,\R)$ we denote by 
$\varphi_A:\Flag_{n,\underline{k}} \to \Flag_{n,\underline{k}}$ the diffeomorphism
$\varphi_A(V_\ast)= A V_\ast$.
Notice that $\varphi:\SL(n,\R) \to \Diff(\Flag_{n,\underline{k}})$, $A\mapsto \varphi_A$,  is
a representation of $\SL(n,\R)$ as a group of diffeomorphisms of the
flag manifold $\Flag_{n,\underline{k}}$. 
For each symmetric matrix $A\in\SL(n,\R)$  with simple spectrum 
the diffeomorphism $\varphi_A:\Flag_{n,\underline{k}} \to \Flag_{n,\underline{k}}$
 has a unique attractive fixed point whose basin of attraction has full measure.
This property alone easily implies that $\Flag_{n,\underline{k}}$ is a boundary of
$\SL(n,\R)$. Note that the flag manifolds include the Grassman ones\,
$\Grass_{n,k}=\Grass_k(\R^n)=\{\,V\subset\R^n\,:\,
V\,\text{ is a linear subspace, }\, \dim V=k\,\}$,
because $\Grass_{n,k}=\Flag_{n,(k)}$.
The maximal boundary of $\SL(n,\R)$ is the full flag manifold $\Flag_{n,\underline{k}}$
determined by the sequence $\overline{k}=(1,2,\ldots,n)$.
Next we define the class of isotropic flag manifolds which are boundaries for the group
$G=\Sp(n,\R)$.  We consider on $\R^{2n}$ the usual
symplectic structure defined by the matrix 
$$J=\matriz{O}{-I}{I}{O}\;.$$
Given a sequence $\underline{k}=(k_1,k_2,\ldots,k_m)$ with 
$1\leq k_1< k_2 <\ldots < k_m\leq n$,
a $\underline{k}$-flag $V_\ast$ in $\R^{2n}$ is called {\em isotropic}\, iff\,
the largest subspace $V_m\subset \R^{2n}$ is isotropic, i.e.,
$u^t J v=0$ for every $u,v\in V_m$. When $\underline{k}=(1,2,\ldots,n)$ the $\underline{k}$-flag $V_\ast$ 
is said to be a {\em Lagrangian flag}. We denote by $\Fliso_{n,\underline{k}}=\Fliso_{\underline{k}}(\R^{2n})$
the submanifold of all isotropic $\underline{k}$-flags in $\Flag_{\underline{k}}(\R^{2n})$.


\bigskip

To each matrix $B$ in $\SL(n,\R)$, resp. $\Sp(n,\R)$,
let us associate an {\em angle matrix} $R\in\OO(n)$, resp $R\in\Un(n)$,
and a vector $b=(b_1,\ldots, b_n)\in\R^n$ with $b_1\geq b_2\geq \ldots \geq b_n>0$
where $b$ consists of the ordered {\em singular values} of $B$ and $R$ is such that for some
other orthogonal matrix $S$ we have $B=R\,D_b\,S$. This is the singular value
decomposition of matrix $B$.
We look at the angle $R$ and the singular values $b$ as coordinates of $B$.
Now, let $A$ be another matrix in $\SL(n,\R)$, resp. $\Sp(n,\R)$, and let
$R'$ be the angle of $A\,B$, while $b'$ be the singular values of $A\,B$.
We would like to understand the transformation $B\mapsto A\,B$ in terms of these of coordinates,
i.e., to describe  $(R,b)\mapsto (R',b')$ in terms of matrix $A$. 
We have $R'=\varphi_A(R)=A\ast R$. Whence this first coordinate map relates to
the action of our group on the complete flag manifold, resp. Lagrangian flag manifold.
About the second coordinate map, note that
$\nrm{A\,B}_{{\rm hs}}=\nrm{D_{b'}}_{{\rm hs}}= Q_{A^t A,b}(R)^{1/2}$.
Therefore, the function $ Q_{A^t A,b}$ gives us valuable information on the
expansion of the singular values as we multiply by $A$.
The singular values' vector $b'$ can be computed as the singular values of the
nilpotent matrix $U\,D_b$, where $U$ is the matrix in the QR-decmposition 
$A\,R=K\,U$ of $A\,R$.
We belive that to analyze Lyapunov exponents it is worth to studing the geometry
of the tranformations $\varphi_A$ and of the functions
$Q_{A^t A,b}$ on the maximal boundary. In this article we only study the geometry of these
objects without any reference to possible implications on Lyapunov exponents.

\bigskip

We will denote by $G_n$ any of the groups $\SL(n,\R)$ and $\Sp(n,\R)$
and by  $\gfrak_n$ the corresponding Lie algebra.
The first is the group of real $n\times n$ matrices with determinant one, while the second
is the group of real $2n\times 2n$ symplectic matrices.
Every theorem below will hold  in two different contexts:
$G_n=\SL(n,\R)$ and $\gfrak_n=\slin(n,\R)$ on one hand, and  $G_n=\Sp(n,\R)$ and $\gfrak_n=\spl(n,\R)$
on the other.

Any sequence $V_\ast = (V_1, \ldots , V_k)$ of linear spaces
$V_1\subset V_2\subset \ldots \subset V_k\subset \R^n$ such that
$\dim V_i=i$ for every $i=1,\ldots, k$ is called a {\em $k$-flag} in $\R^n$,
and $\Flag_{n,k} = \Flag_k(\R^n)$ will denote the {\em $k$-flag manifold}
consisting of all $k$-flags  in $\R^n$. We consider on $\R^{2n}$ the usual
symplectic structure defined by the matrix 
$$J=\matriz{O}{-I}{I}{O}\;.$$
A $k$-flag $V_\ast$ in $\R^{2n}$ is called {\em isotropic}\, iff\,
the subspace $V_k\subset \R^{2n}$ is isotropic, i.e.,
$u^t J v=0$ for every $u,v\in V_k$. When $k=n$ the $k$-flag $V_\ast$ 
is said to be a {\em Lagrangian flag}. We denote by $\Fliso_{n,k}$
the submanifold of all isotropic $k$-flags in $\Flag_k(\R^{2n})$.
The subscript $sp$ will be omitted on symplectic contexts, i.e., 
we  adopt the convention that in each context where $G_n$ and $\Flag_{n,k}$ appear,
$\Flag_{n,k}$ stands for the isotropic $k$-flag manifold $\Fliso_{n,k}$
whenever $G_n=\Sp(n,\R)$. There is a natural action of $G_n$ on $\Flag_{n,k}$
defined by
\[ G_n\times \Flag_{n,k}\to \Flag_{n,k},\quad
(\, A, \,V_\ast=(V_1,\ldots, V_k)\,)\mapsto A\,V_\ast=(A V_1,\ldots, A V_k) \;.\]
Given $A\in G_n$ we denote by 
$\varphi_A:\Flag_{n,k} \to \Flag_{n,k}$ the diffeomorphism
$\varphi_A(V_\ast)= A V_\ast$.
Notice that $\varphi:G_n \to \Diff(\Flag_{n,k})$, $A\mapsto \varphi_A$,  is
a representation of $G_n$ as a group of diffeomorphisms of the
flag manifold $\Flag_{n,k}$. Given $A\in\gfrak_n$, the family
$\{\,\varphi_{e^{t A}}\,\}_{t\in\R}$ is a flow on $\Flag_{n,k}$ that we shall
simply denote by $\varphi_A^t$. Note that   $\varphi_A^1 = \varphi_{e^A}$.

We introduce now a family of quadratic functions on the flag manifolds $\Flag_{n,k}$.
Given two linear subspaces $E\subset F$ denote by $F\ominus E$ the orthogonal complement
$F\ominus E=E^\perp\cap F$.
Given a matrix $A\in \gfrak_n$ and a vector $b\in\R^k_+$ we define the
function $Q_{A,b}:\Flag_{n,k}\to\R$ as follows:
For each $V_\ast\in \Flag_{n,k}$ take a vector $v_1 \in V_1$ with
$\nrm{v_1}=b_1$, and  a vector $v_i \in V_i\ominus V_{i-1}$ with
$\nrm{v_i}=b_i$ ($2\leq i \leq k$) and then set
$$ Q_{A,b} (V_\ast)=\frac{1}{k} \,\sum_{i=1}^k  \langle A\,v_i, v_i\rangle \;. $$
Because $V_1$ and  $V_i\ominus V_{i-1}$ ($2\leq i \leq k$) are one-dimensional spaces
this definition is clearly independent of the choice of  vectors $v_i$.

We denote by $\Diag_n^+$ the subgroup of positive diagonal matrices in $G_n$.
Any subgroup of $G_n$ which is conjugated to $\Diag_n^+$ by some orthogonal matrix in $G_n$
will be referred as a {\em commutative positive symmetric subgroup}.
Notice that matrices in such subgroups are always  symmetric and positive definite.
The Lie algebra of a commutative positive symmetric subgroup will be referred as
a {\em commutative symmetric subalgebra}.

\bigskip

We shall call {\em stratification on a manifold } $\Fscr$ to any collection $\Sscr$ of closed 
\footnote{compact without boundary} connected submanifolds, that we refer as {\em strata}  of $\Fscr$,   such that 
\begin{enumerate}
\item $\Fscr\in\Sscr$ is the unique stratum of dimension  $N=\dim \Fscr$;
\item There are several strata of dimension zero;
\item Each stratum $S\in\Sscr$ of dimension $d<N$ is contained in
some stratum $S'\in\Sscr$ of dimension $d+1$;
\item Given $S, T\in\Sscr$ such that 
$S\subset T$ and $S\neq T$ then  $\dim S < \dim T$;
\item Given $S, T\in\Sscr$,\, either $S\cap T=\emptyset$
 or else $S\cap T\in\Sscr$.
\end{enumerate}
We denote by $\Sscr^i$ the union of all strata $S\in\Sscr$ with $\dim S=i$.
It follows from 2. and 3. that each $\Sscr^i$ is non-empty.
By 3. we have $\Sscr^0 \subset \Sscr^1 \subset \ldots \subset \Sscr^N=\Fscr$,
and by 4. and 5. each $\Sscr^i-\Sscr^{i-1}$ is a (disconnected) manifold of dimension $i$.
We say that a stratification $\Sscr$ on $\Fscr$ is {\em invariant under a diffeomorphism }
$\varphi:\Fscr\to\Fscr$\, iff\, $\varphi(S)=S$ for every stratum $S\in\Sscr$.

\bigskip

We prove the following results about the objects $\Flag_{n,k}$, $\varphi_A^t$ and $Q_{A,b}$.
We say that two symmetric matrices $A,B\in\gfrak_n$ {\em share the same ordered
eigen-directions} \, iff\, there is a common basis of eigenvectors $v_1,\ldots, v_n$
such that the eigenvalues of $A$, $A\,v_i=\lambda_i(A)\,v_i$, 
are ordered in the same way as the eigenvalues of $B$,  $B\,v_i=\lambda_i(B)\,v_i$, which means that
$\lambda_i(A)>\lambda_j(A)$ $\Leftrightarrow$ $\lambda_i(B)>\lambda_j(B)$,
for $i,j=1,\ldots, n$.
For instance, $A$ and $e^A$ share the same ordered
eigen-directions.

\bigskip

\noindent
{\bf Theorem A} {\em 
Given  $b_1>b_2>\ldots >b_k>0$, \,
and two symmetric matrices $A,H\in\gfrak_n$  sharing the same ordered
eigen-directions then $Q_{A,b}:\Flag_{n,k}\to\R$ is a strict Lyapunov function for 
$\varphi_H^t:\Flag_{n,k}\to\Flag_{n,k}$.
}

\bigskip

\noindent
{\bf Theorem B} {\em Given a commutative symmetric subalgebra $\hfrak\subset \gfrak_n$,
there is a stratification on $\Flag_{n,k}$ which is invariant
under both $\varphi_H^t:\Flag_{n,k}\to\Flag_{n,k}$ and the gradient
flow of the function $Q_{H,b}:\Flag_{n,k}\to\R$, for each $H\in \hfrak$ and each $b\in\R^k_+$.
}

\bigskip

\noindent
{\bf Theorem C} {\em 
Given a symmetric
matrix $A\in \gfrak_n$ with simple spectrum  and $b_1>b_2>\ldots >b_k>0$,\,
$Q_{A,b}:\Flag_{n,k}\to\R$ is a $\Z_2$-perfect Morse function.
}

\bigskip

For reader's reference we recall the key definitions above.
A matrix $A\in \gfrak_n$ is said to have {\em simple spectrum} \, iff\,
the eigenvalues of $A$ are all simple.
A function $Q:\Fscr\to\R$ is said to be a {\em  Lyapunov function}
for a flow $\varphi^t:\Fscr\to \Fscr$\, iff\, 
$Q(\varphi^t(x))\geq Q(x)$, for all $x\in \Fscr$ and $t\geq 0$.
Function $Q$ is called a {\em  strict Lyapunov function} when furthermore it satisfies
$Q(\varphi^t(x)) = Q(x)$ \, $\Leftrightarrow$\, $\varphi^t(x)=x$, for every $x\in \Fscr$ and $t>0$.
A smooth function $Q:\Fscr\to\R$ is called a Morse function \,iff\, all its critical points
are non-degenerate. Function $Q$ is said to be $\Z_2$-perfect \, iff\, the Poincar\'{e}
polynomial of $\Fscr$ with coefficients in the field $\Z_2$,\,
$\Pscr_t(\Fscr,\Z_2)=\sum_{i=0}^{\dim \Fscr} \dim H_i(\Fscr;\Z_2)\,t^i$,\, coincides with the
Morse polynomial of $Q$. Finally, the Morse polynomial of $Q$ is defined by \,
$\Mscr_t(Q)= \sum_{i=0}^{\dim \Fscr} c_i(Q)\,t^i$, \,  where $c_i(Q)$
is the number of critical points with index $i$,
the index of a critical point being the number of negative eigenvalues of its
Hessian matrix.

\bigskip

\section{Notation}

All the notation below is defined in the text's body, but we gather it here for
an easy reference.
\subsection{Matrices}
\begin{tabular}{cl}
 $\Mat_{n\times k}(\R)$ &  the space of real $n\times k$  matrices\\
 $(X)_i$ & the column $i$ of matrix $X$\\
 $[X]_i$ &  the submatrix formed by the first $i$ columns of $X$\\
 $[X]_I$ & the submatrix formed by the columns of $X$ with index $i\in I$\\
 $\langle X\rangle_I$ & the space spanned by the columns of $X$ with index $i\in I$\\
 $\delta_{i,j}$ & the Kronecker symbol,\,
 $\delta_{i,j}=\left\{\begin{array}{ccc} 1 & \text{ if } & i=j \\ 0 & \text{ if } & i\neq j \end{array}\right.$\\
 $I_{n,k}$ & the $n\times k$ matrix with entries $\delta_{i,j}$\\
 $X^\triangleleft$ & the matrix
 $(x^\triangleleft_{i,j})_{i,j}$, where 
$x^\triangleleft_{i,j}=\left\{
\begin{array}{ccl}
x_{i,i} & \text{ if } & i=j \\
x_{i,j} + x_{j,i} & \text{ if } & i < j \\
0 & \text{ if } &  j< i
\end{array}\right. $ for $X=(x_{i,j})_{i,j}$.
\end{tabular}

\subsection{Inner Products}
\begin{tabular}{cl}
$\langle X,Y\rangle_{{\rm hs}}$  &  the  Hilbert-Schmidt inner product of matrices $X$ and $Y$\\
$\nrm{X}_{{\rm hs}}$  &   the Hilbert-Schmidt norm of matrix $X$\\
\end{tabular}

\subsection{Groups}
\begin{tabular}{cl}
$\OO(n)$ &   the real orthogonal group\\
$\Un(n)$ &   the unitary group\\
$\UT(n)$ &   the upper triangular group\\
$\UT_+(n)$ &   the positive diagonal upper triangular group\\
\end{tabular}

\subsection{Homogeneous Spaces}
\begin{tabular}{cl}
$\underline{k}$ &  the signature $(k_1,k_2,\ldots, k_m)$, with $1\leq k_1< k_2< \ldots < k_m\leq n$, of a flag\\
$\Flag_{n,\underline{k}}$ & the $\underline{k}$-flag manifold\\
$\Fliso_{n,\underline{k}}$ & the isotropic  $\underline{k}$-flag manifold\\
$\OO_{n,k}$ & the $k$-orthogonal frame manifold\\
$\Osp_{n,k}$ & the $k$-unitary frame manifold\\
\end{tabular}

\bigskip

\bigskip

\section{Flag Manifolds}

A sequence of numbers $\underline{k}=(k_1,k_2,\ldots,k_m)$ with 
$1\leq k_1< k_2 <\ldots < k_{m}\leq  n$, is called a {\em signature } of length $m$
in the set $\{1,2,\ldots, n\}$.
We denote by $\Lambda(n)$ the set of all signatures in $\{1,2,\ldots, n\}$, which is partially ordered as
follows: We say that $\underline{k}\geq \underline{k}'$ \, iff\, $\underline{k}'$ is a subsequence
of $\underline{k}$. The sequence $(1,2,\ldots, n)$ is the unique maximal element of $\Lambda(n)$,
while every sequence of length one $(i)$, with $1\leq i\leq n$, is a minimal element of $\Lambda(n)$.
Given a signature $\underline{k}\in\Lambda(n)$ of length $m$,
any sequence of linear spaces $V_\ast = (V_1, \ldots , V_{m})$ such that
$V_1\subset V_2\subset \ldots \subset V_{m}\subset \R^n$ and
$\dim V_i=k_i$ for every $i=1,2,\ldots, m$ is called a {\em $\underline{k}$-flag in $\R^n$}.
We denote by $\Flag_{n,\underline{k}} = \Flag_{\underline{k}}(\R^n)$   the space of all
 $\underline{k}$-flags  in $\R^n$. 
 
\begin{prop}
The space $\Flag_{n,\underline{k}}$ is a compact connected manifold 
 of dimension $\frac{n(n-1)}{2}- \sum_{i=1}^{m+1} \frac{n_i(n_i-1)}{2}$,
 where $n_1=k_1$, $n_i=k_i-k_{i-1}$ for $i=2,\ldots, m$ and $n_{m+1}=n-k_m$.
\end{prop}

The order on $\Lambda(n)$ induces a hierarchy on flag manifolds.
Given $\underline{k}, \underline{k}'\in\Lambda(n)$, if $\underline{k}\geq \underline{k}'$
and $k'_i=k_{s_i}$, for $i=1,\ldots, m'$, then 
there is a natural projection map
$\pi:\Flag_{n,\underline{k}}\to \Flag_{n,\underline{k}'}$
defined by $\pi(V_1,\ldots, V_m)=(V_{s_1},\ldots, V_{s_m})$.
The flag $\pi(V_\ast)$ is called the {\em restriction } of $V_\ast$.
This projection map is a smooth submersion.
Note that if the signature $\underline{k}'$ satisfies
$k'_{m}<n$  and $\underline{k}$ is obtained appending $n$ at the end of $\underline{k}'$, 
i.e., $\underline{k} = (k_1',\ldots, k'_m,n)$, then the projection
 $\pi:\Flag_{n,\underline{k}}\to \Flag_{n,\underline{k}'}$
is actually a  diffeomorphism.
Whence, we identify these two flag manifolds and restrict our attention to flag manifolds with signatures in the
set $\Lambda'(n)=\{\,\underline{k}\in\Lambda(n)\,:\, k_m<n\,\}$.
The flag manifold associated with the maximum signature $\underline{k}=(1,\ldots, n-1)$
is called the {\em complete flag manifold}, which as explained in the introduction,
is the maximal boundary of the group $\SL(n,\R)$.
Flag manifolds include an important subclass of varieties called
the Grassmann manifolds. For each $k=1,\ldots, n-1$,
$\Grass_{n,k}:=\Flag_{n,(k)}$ is the Grassmann manifold of all $k$-dimensional
subspaces $V\subset\R^n$.
Note that for $k=1$ the Grassmann manifold $\Grass_{n,1}$ is the real projecticve space.
 The action of the group $\SL(n,\R)$ on the flag manifold
$\Flag_{n,\underline{k}}$ is
$\SL(n,\R)\times \Flag_{n,\underline{k}}\to \Flag_{n,\underline{k}}$, 
$g\,V_\ast=(g V_1,\ldots, g V_m)$ with $V_\ast=(V_1,\ldots, V_m)$.
Two important, but obvious, facts about this action are:

\begin{prop}
The action of $\SL(n,\R)$ on $\Flag_{n,\underline{k}}$ is transitive, and 
for $\underline{k}\geq \underline{k}'$, the
 projection  $\pi:\Flag_{n,\underline{k}}\to \Flag_{n,\underline{k}'}$
 is $\SL(n,\R)$-equivariant, i.e.,
 $\pi(g\,V_\ast)=g\,\pi(V_\ast)$, for every $g\in\SL(n,\R)$ and $V_\ast\in\Flag_{n,\underline{k}}$.
\end{prop}

\bigskip

\section{Isotropic Flag Manifolds}

Given a signature $\underline{k}\in\Lambda(n)$, a $\underline{k}$-flag
$V_\ast\in\Flag_{2n,\underline{k}}=\Flag_{\underline{k}}(\R^{2n})$ is 
called {\em isotropic} \, iff\, the largest subspace $V_{k_m}$ is isotropic,
which means that the restriction of the linear symplectic structure to $V_{k_m}$ is
the zero $2$-form. We denote by $\Fliso_{n,\underline{k}}=\Fliso_{\underline{k}}(\R^{2n})$
the space of all isotropic $\underline{k}$-flags on $\R^{2n}$.

\begin{prop}
The space $\Fliso_{n,\underline{k}}$ is a compact connected manifold 
 of dimension $n^{2}- \sum_{i=1}^{m+1} \frac{n_i(n_i-1)}{2}$,
 where $n_1=k_1$, $n_i=k_i-k_{i-1}$ for $i=2,\ldots, m$ and $n_{m+1}=n-k_m$.
\end{prop}

Note that if $\underline{k}\geq \underline{k}'$  then the projection
$\pi:\Flag_{2n,\underline{k}}\to \Flag_{2n,\underline{k}'}$ maps
$\Fliso_{n,\underline{k}}$ onto $\Fliso_{n,\underline{k}'}$ and, therefore, restricts to a map we still denote by
$\pi:\Fliso_{n,\underline{k}}\to \Fliso_{n,\underline{k}'}$ .
Clearly, the action of the subgroup $\Sp(n,\R)\subset\SL(2n,\R)$  on $\Flag_{2n,\underline{k}}$
leaves invariant each isotropic flag manifold
$\Fliso_{n,\underline{k}}$, and whence induces by restriction an action
$\Sp(n,\R)\times \Fliso_{n,\underline{k}}\to \Fliso_{n,\underline{k}}$.

\begin{prop}
The action of $\Sp(n,\R)$ on $\Fliso_{n,\underline{k}}$ is transitive, and 
for $\underline{k}\geq \underline{k}'$, the
 projection  $\pi:\Fliso_{n,\underline{k}}\to \Fliso_{n,\underline{k}'}$
 is $\Sp(n,\R)$-equivariant.
\end{prop}

For $\underline{k}=(1,2,\ldots, n)\in\Lambda(n)$, the $\underline{k}$-flags are
called {\em Lagrangian flags}, and the corresponding manifold $\Fliso_{n,\underline{k}}$,
which as we have see is the maximal boundary of the group $\Sp(n,\R)$, is called the
{\em  Lagrangian flag manifold}.

\bigskip

\section{Orthogonal Frame Manifolds} \label{orth:frames}

Given  $1\leq k \leq n$, we denote by $\Mat_{n\times k}(\R)$ the space of all
real $n\times k$ matrices.
A matrix $X\in \Mat_{n\times k}(\R)$ is called {\em orthogonal} \, iff\,
$X^t\,X=I$, i.e., the $k$ columns of $X$ form an orthonormal
family of vectors in $\R^n$. Such matrices will be referred as 
$k$-ortho-frames.
We denote by $\OO_{n,k}$ the  manifold of $k$-ortho-frames
\[ \OO_{n,k}=\{\, X\in\Mat_{n\times k}(\R)\,:\, X\; \text{ is orthogonal}\;\}\;.\]
This is a compact connected manifold of dimension $k\,(2n-k-1)/2$.
When $k=n$, $\OO_{n,n}$ is identified with the {\em orthogonal group} $\OO(n)$.
The group $\OO(n)$ acts by left multiplication on the space $\OO_{n,k}$ and
 $\OO(k)$  acts by right multiplication on $\OO_{n,k}$.
Each space $\OO_{n,k}$ has a distinguished element $I_{n,k}\in\OO_{n,k}$,
which is the matrix formed by the first $k$ columns of the identity matrix
$I_n\in\Mat_{n\times n}(\R)$.
Given a matrix $X\in\OO_{n,k}$ let $(X)_i$ be its $i^{\text{th}}$-column.
Given a set $I\subset\{1,2,\ldots,n\}$ we define
$\langle X\rangle_I$ to be the linear span of the columns
$(X)_i$ with $i\in I$. Finally, given $\underline{k}\in\Lambda(n)$,
we define the projection $p:\OO_{n,k_m}\to\Flag_{n,\underline{k}}$  by \,
$p(X)=\left(\,\langle X\rangle_{\{1,\ldots,k_1\}}, \,\langle X\rangle_{\{1,\ldots,k_2\}}, \,
\ldots\, ,\,\langle X\rangle_{\{1,\ldots,k_m\}}  \,\right)$.
We shall call  $V_{\underline{k}}=p(I_{n,k})$ the {\em canonical $\underline{k}$-flag}.
Consider the subgroup $\OD_{\underline{k}}$ of  $\OO(k_m)$ consisting of block diagonal matrices
of the form $R=\diag(R_1,\ldots, R_m)$ with  $R_1\in\OO(k_1)$ and $R_i\in\OO(k_i-k_{i-1})$
for $i=2,\ldots, m$.
Let $\OO_{n,k}/\OD_{\underline{k}}$  be the quotient orbifold, i.e., the space of orbits
of the right action of $\OD_{\underline{k}}$ on $\OO_{n,k}$.
Note the projection $p:\OO_{n,k_m}\to\Flag_{n,\underline{k}}$ defined above 
is invariant under the action of $\OD_{\underline{k}}$.
Whence it induces a quotient map $\overline{p}:\OO_{n,k_m}/\OD_{\underline{k}}\to\Flag_{n,\underline{k}}$,
and it is not difficult to see that in fact

\begin{prop}
$\overline{p}:\OO_{n,k_m}/\OD_{\underline{k}}\to\Flag_{n,\underline{k}}$ \, is a diffeomorphism.
\end{prop}

Given $1\leq k' \leq k$ and a matrix $E\in\Mat_{n\times k'}(\R)$ we denote by $[E]_k\in\Mat_{n\times k}(\R)$
the matrix formed by the first $k$ columns of $E$. With this notation we define a
projection $\pi:\OO_{n,k}\to\OO_{n,k'}$ by $\pi(X)=[X]_{k'}$.  Later, after having defined a natural action
of $\SL(n,\R)$ on $\OO_{n,k}$, we shall show this projection is $\SL(n,\R)$-equivariant.
For now we make two remarks about $\pi:\OO_{n,k}\to\OO_{n,k'}$. First, 
given $\underline{k}\geq \underline{k}'$ in $\Lambda(n)$, the following diagram commutes
$$\begin{CD}
\OO_{n,k_m}    @>\pi>> \OO_{n,k'_m}  \\
 @VpVV       @VpVV      \\
\Flag_{n,\underline{k}}   @>\pi>> \Flag_{n,\underline{k}'}   \\
\end{CD}\;.$$
Let  $\phi:\OD_{\underline{k}}\to \OD_{\underline{k}'}$ be the group homomorphism
which maps each  $k_m\times k_m$ matrix in $\OD_{\underline{k}}$ to the $k'_m\times k'_m$ submatrix
with indices in the set $\{1,\ldots,k'_m\}$, which clearly belongs to $\OD_{\underline{k}'}$.
The second remark is that $\pi:\OO_{n,k}\to\OO_{n,k'}$
is $\phi$-equivariant, in the sense that
$\pi(X\,g)=\pi(X)\,\phi(g)$ for every $X\in\OO_{n,k}$ and $g\in\OD_{\underline{k}}$.

\bigskip

\section{Unitary Frame Manifolds} \label{un:frames} 

Given  $1\leq k \leq n$, a matrix $X\in\Mat_{2n\times k}(\R)$ is called {\em symplectic} \, iff\,
$X^t\,J\,X=0$, i.e., the $k$ columns of $X$ span an isotropic subspace of $\R^{2n}$. 
A matrix  $X\in\Mat_{2n\times k}(\R)$ is called {\em unitary} \, iff\, $X$ is
orthogonal and symplectic. Such matrices will be referred as 
$k$-unitary-frames.
We denote by $\Osp_{n,k}$ the  manifold of $k$-unitary-frames
\[ \Osp_{n,k}=\{\, X\in\Mat_{2n\times k}(\R)\,:\, X\; \text{ is unitary}\;\}\;.\]
This a compact connected manifold of dimension $k\,(2n-k)$.

Let $\Un(n)$ denote the {\em unitary group} which is formed of all matrices
$A\in\Mat_{2n\times 2n}(\R)$ such that $A^t\,J\,A=J$ and $A^t\,A=I$.
When $k=n$, we write $\Osp_n$ in stead of $\Osp_{n,n}$.
Any matrix $X=\coluna{X_1}{X_2}\in\Osp_n$ determines a unique matrix
$\matriz{X_1}{-X_2}{X_2}{X_1}$  which belongs to $\Un(n)$.
We can therefore identify
$\Osp_n$ with $\Un(n)$.
The group $\Un(n)$ acts by left multiplication on the manifold $\Osp_{n,k}$ and
 $\OO(k)$  acts by right multiplication on $\Osp_{n,k}$.
Note that if $X$ is symplectic, resp. unitary,
 and $R\in\OO(k)$ then $X\,R$ is also symplectic, resp. unitary.
Whence the right multiplication action of the group $\OO(k)$ leaves $\Osp_{n,k}$ invariant.
Given $\underline{k}\in\Lambda(n)$, by restriction the projection $p:\OO_{2n,k_m}\to\Flag_{2n,\underline{k}}$ 
induces a projection $p:\Osp_{n,k_m}\to\Fliso_{n,\underline{k}}$, which is $\OD_{\underline{k}}$-invariant. 
Therefore, it induces
a quotient map $\overline{p}:\Osp_{n,k_m}/\OD_{\underline{k}} \to\Fliso_{n,\underline{k}}$,
which is again a diffeomorphism.
We have a commutative diagram
$$\begin{CD}
\Osp_{n,k_m}    @>\pi>> \Osp_{n,k'_m}  \\
 @VpVV       @VpVV      \\
\Fliso_{n,\underline{k}}   @>\pi>> \Fliso_{n,\underline{k}'}   \\
\end{CD} $$
where the top map is obtained restricting the projection\,
$\pi:\OO_{2n,k_m} \to \OO_{2n,k'_m}$.
We shall later define a natural action
of $\Sp(n,\R)$ on $\Osp_{n,k}$ for which this projection is $\Sp(n,\R)$-equivariant.
For now we can say that the map $\pi:\Osp_{n,k_m}\to\Osp_{n,k'_m}$ in the above diagram
is $\phi$-equivariant, where $\phi:\OD_{\underline{k}}\to \OD_{\underline{k}'}$ is the same group homomorphism 
introduced in the previous section.

\bigskip

\bigskip

\section{Action on Frame Manifolds}


Let $\UT(k,\R)$ be the group of upper triangular real matrices
with non-zero diagonal, and $\UTpos(k,\R)$ denote the subgroup of upper triangular real matrices
with positive diagonal. 
The  action of $\SL(n,\R)$ on  $\OO_{n,k}$ uses the QR-decomposition.

\begin{teor} [QR-decomposition] Given $A\in \Mat_{n\times k}(\R)$ there are unique matrices
$K\in\OO_{n,k}$ and $U\in \UTpos(k,\R)$  such that $A=K\,U$.
\end{teor}

The QR-decomposition is used to define the following projection map
\begin{equation} \label{OF:projection}
\Pi:\Mat_{n\times k}(\R)\to \OO_{n,k}\,,\; \text{by } \; \Pi(A)=K\,,
\end{equation}
and with it the multiplication $\ast:\SL(n,\R)\times \OO_{n,k}\to \OO_{n,k} $
\begin{equation} \label{OF:action}
A\ast X=\Pi(A\,X)\,.
\end{equation}

\bigskip

\begin{prop} \label{ast:action}
The operation $\ast$ defined above is a left action of $\SL(n,\R)$ on $\OO_{n,k}$.
\end{prop}

\dem
Given $A\in\SL(n,\R)$ and $X\in \OO_{n,k}$, 
consider the QR-decompositions 
$A\,X=K_1\,U_1$ of $A\,X$, and $B\,K_1=K_2\,U_2$ of $B\,K_1$.
Notice that $U=U_2\,U_1\in\UTpos(k,\R)$  so that
$B\,A\,X=B\,K_1\,U_1=K_2\,U_2\,U_1 = K_2\,U$ is the QR-decomposition
of $B\,A\,X$.
Therefore
\[ B\ast(A\ast X) = B\ast K_1 = \Pi(B\,K_1) = K_2 = \Pi(B\,A\,X) = (B\,A)\ast X\;.
\]
\cqd

\bigskip

Given $1\leq k\leq n$, denote by $I_{n,k}$ the $n\times k$ matrix whose columns are the first $k$ columns
of the $n\times n$ identity matrix. Notice that for any matrix $X\in\Mat_{m\times n}(\R)$, the 
product $X\,I_{n,k}$ is the submatrix of $X$ formed of its first $k$ columns.

\begin{prop} \label{SL:equivariance} The following projections are $\SL(n,\R)$-equivariant:
\begin{enumerate}
\item  $\pi:\OO_{n,k}\to\OO_{n,k'}$ for any $k\geq k'\geq 1$, and
\item  $p:\OO_{n,k_m}\to\Flag_{n,\underline{k}}$ for any
$\underline{k}\in\Lambda(n)$.
\end{enumerate}
\end{prop}

\dem
Let us prove first that  $\pi:\OO_{n,k}\to\OO_{n,k'}$ is $\SL(n,\R)$-equivariant.
Take $A\in \SL(n,\R)$, $X\in \OO_{n,k}$ and let $A\,X=K\,U$ be its QR-decomposition with
$K\in \OO_{n,k}$ and $U\in\UTpos(k',\R)$. Then
$$ A\,X\, I_{k,k'} =K\,U\, I_{k,k'} = K\,I_{k,k'}\,(I_{k,k'})^t\, U\, I_{k,k'}
=(K\,I_{k,k'})\,((I_{k,k'})^t\, U\, I_{k,k'}) = K'\, U'\;, $$
where $K'=K\,I_{k,k'}$ is orthogonal and
$U'=(I_{k,k'})^t\, U\, I_{k,k'}\in\UTpos(k',\R)$.
Therefore, this is the QR-decomposition of $A\,X\, I_{k,k'}$.
Finally, since $\pi(X)=X\, I_{k,k'}$, we get
\begin{align*}
A\ast \pi(X) &= A\ast (X\,I_{k,k'}) = K' = K\,I_{k,k'} \\
&= (A\ast X)\,I_{k,k'} =  \pi(A\ast X) \;.
\end{align*}

Now, for the equivariance of $p:\OO_{n,k}\to\Flag_{n,k}$, 
first remark that given $K\in\OO_{n,k}$ and $U\in\UTpos(k)$ we have
$p(K\,U)=p(K)$, and given $A\in\SL(n,\R)$ and $X\in\OO_{n,k_m}$ we
have $A\,p(X)=p(A\,X)$. Considering then the QR-decomposition $A\,X=K\,U$,
we get $p(A\,X)=p(K\,U)=p(K)=p(A\ast X)$, and whence
$p(A\ast X)=A\,p(X)$.
\cqd

\bigskip

Let us now turn to unitary frame manifolds.
Consider the set $\Sp_{n, k}$ of symplectic
$2n\times k$ matrices. 
Symplectic matrices, defined in section~\ref{un:frames},
are those whose columns span isotropic subspaces of $\R^{2n}$.
Clearly, $\Sp_{n, k}$ is not a linear space, but it has
two important obvious properties whose proofs are left to the reader.

\begin{prop} \label{Sp:invariance}
The space $\Sp_{n, k}$ is invariant under:
\begin{enumerate}
\item the left multiplication action of $\Sp(n,\R)$;
\item the right multiplication action of $\GL(k,\R)$.
\end{enumerate}
\end{prop}

Consider the projection $\Pi:\Mat_{n\times k}(\R)\to \OO_{n,k}$
defined in~(\ref{OF:projection}).

\begin{prop} We have
\begin{enumerate}
\item[(a)] $\Pi( \Sp_{n,k} )\subseteq \Osp_{n,k}$,\, and
\item[(b)] $\Osp_{n,k}$ is invariant under the left action of $\Sp(n,\R)$.
\end{enumerate}
\end{prop}

\dem
Let us prove (a).
Given $A\in \Sp_{n,k}\subset \Mat_{2n\times k}(\R)$, let $A=K\,U$ be its QR-decomposition.
Then $K$ is orthogonal and $U\in\UTpos(k,\R)$, and, 
since $\Osp_{n,k}=\OO_{2n,k}\cap \Sp_{n,k}$, we are left to prove
that $K$ is symplectic. 
This follows by item 2 of 
proposition~\ref{Sp:invariance}  because $K=A\,U^{-1}$.
Finally, (b) follows combining (a) with
item 1 of proposition~\ref{Sp:invariance}.
\cqd

Therefore, by restriction we obtain a left action of $\Sp(n,\R)$ on $\Osp_{n,k}$.
The projections $\pi:\OO_{2n,k}\to\OO_{2n,k'}$, for $k\geq k'\geq 1$, and
$p:\OO_{2n,k_m}\to\Flag_{2n,\underline{k}}$, for $\underline{k}\in\Lambda(n)$,
respectively induce by restriction maps
$\pi:\Osp_{n,k}\to\Osp_{n,k'}$  and
$p:\Osp_{n,k_m}\to\Fliso_{n,\underline{k}}$. It is now easy to check that
both these projections are $\Sp(n,\R)$-equivariant.


\bigskip

\section{Tangent Actions}

Let $\ut(k,\R)$ be the Lie algebra of upper triangular real matrices.
Given a matrix $X=(x_{i,j})_{i,j}\in \Mat_{k\times k}(\R)$,  denote by
$X^\triangleleft=(x^\triangleleft_{i,j})_{i,j}$ the matrix in $\ut(k,\R)$
defined by  
$$x^\triangleleft_{i,j}=\left\{
\begin{array}{ccl}
x_{i,i} & \text{ if } & i=j \\
x_{i,j} + x_{j,i} & \text{ if } & i < j \\
0 & \text{ if } &  j< i
\end{array}\right. $$

Then $X\mapsto X^\triangleleft$ is a linear projection operator.
An important property of this operator is that for all
$A\in \Mat_{k\times k}(\R)$,
\begin{equation}\label{up:op:property}
(A^\triangleleft)^t + A^\triangleleft = A^t + A\;.
\end{equation}

\bigskip

We notice that given $X\in\OO_{n,k}$, 
\[ T_X \OO_{n,k}=\{\, V\in \Mat_{n\times k}(\R)\,:\, X^t V+V^t X = 0 \,\}\;,\]
and $T_U \UT_+(k,\R)= \ut(k,\R)$, for every $U\in\UTpos(k,\R)$, since
$\ut(k,\R)$ is a linear space and $\UT_+(k,\R)$ an open subset of $\ut(k,\R)$.

\begin{teor}[Tangent QR-decomposition]\label{tangent:QR:decomp}
Given $X,V\in\Mat_{n\times k}(\R)$,\; 
let $X=K\,U$ be the QR-decomposition of $X$, where
$K\in\OO_{n,k}$ and  $U\in\UTpos(k,\R)$. Then there are unique matrices
$V_0\in T_K \OO_{n,k}$ and $V_1\in \ut(n,k)$ such that\;
$ V=V_0\,U + K\,V_1 $. 
\begin{eqnarray} \label{V0:eq}
V_0 &=& V\,U^{-1} - K\,(K^t V U^{-1})^\triangleleft \\
V_1 &=&  (K^t V U^{-1})^\triangleleft\,U \label{V1:eq}
\end{eqnarray}
\end{teor}

\dem
Consider the QR-decomposition $X=K\,U$.
Then $X$ is the unique point in the transversal intersection
of the two orbits
$\OO_{n,k}\cdot U$ and $K\cdot\UTpos(k,\R)$.
An easy argument shows that given $V\in \Mat_{n\times k}(\R)$ there are unique
vectors $V_0\in T_K\OO_{n,k}$ and $V_1\in T_U\UTpos(k,\R)=\ut(k,\R)$
such that \, $V= V_0\,U+K\,V_1$.
It follows that $V_0$ satisfies
$K^t V_0+V_0^t K =0$, and $V_1$ is upper triangular.

\begin{figure}[h]
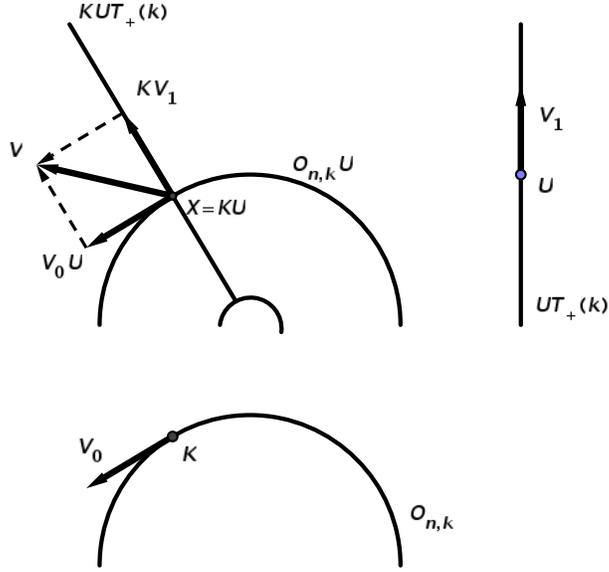

\figura{decomp}{.8} 
\caption{The Tangent QR-decomposition}
\end{figure}

We have
\begin{align*}
K^t V\,U^{-1} + (K^t V\,U^{-1})^t  &=
K^t V\,U^{-1}  +  U^{-t} V^t  K \\
&= K^t (V_0 U + K V_1) \,U^{-1}  +  U^{-t}(V_0 U + K V_1)^t  K\\
&= K^t V_0  + V_1 U^{-1}  +  V_0^t K  +  U^{-t} V_1^t\\
&= \underbrace{ K^t V_0 + V_0^t K }_{=0} + V_1 U^{-1}+ (V_1 U^{-1})^t\\
&=  V_1 U^{-1}+ (V_1 U^{-1})^t
\end{align*}
and since $V_1 U^{-1}$ is upper triangular it follows
from the identity~(\ref{up:op:property})  that 
$V_1 U^{-1} = (K^t V U^{-1})^\triangleleft$, which proves~(\ref{V1:eq}).
Replacing the value~(\ref{V1:eq}) for $V_1$ in the equality
$V= V_0\,U+K\,V_1$ we get~(\ref{V0:eq}).
\cqd

This theorem gives us formulas for the derivatives of the canonical
projections $\Pi:\Mat_{n\times k}(\R)\to\OO_{n,k}$
and $\Pi':\Mat_{n\times k}(\R)\to\UTpos(k,\R)$ associated
with the QR-decomposition.
Given matrices $X, V\in\Mat_{n\times k}(\R)$,\;
assume $X=K\,U$ is the QR-decomposition of $X$. Then we have
\begin{enumerate}
\item $D\Pi_X(V) = V\,U^{-1} - K\,(K^t V U^{-1})^\triangleleft$
\item $D\Pi'_X(V) = (K^t V U^{-1})^\triangleleft\,U$
\end{enumerate}

When $X\in\OO_{n,k}$ these formulas simplify even further
\begin{enumerate}
\item[1'.] $D\Pi_X(V) = V - X\,(X^t V)^\triangleleft$
\item[2'.] $D\Pi'_X(V) = (X^t V)^\triangleleft$
\end{enumerate}

\bigskip

Given a matrix $A\in\slin(n,\R)$ and $\underline{k}\in\Lambda(n)$, 
we define a flow on the $\underline{k}$-flag manifold,
$\varphi_A^t:\Flag_{n,\underline{k}}\to \Flag_{n,\underline{k}}$ by
$\varphi^t_A(V)=e^{t\,A}\,V$. Analogously, given $k\geq 1$
we define a flow on the $k$-ortho frame manifold $\OO_{n,k}$,
$\varphi_A^t:\OO_{n,k}\to \OO_{n,k}$ by
$\varphi^t_A(X)=e^{t\,A}\ast X$. Proposition~\ref{SL:equivariance} shows these two
flows are semiconjugate by the submersion
$p:\OO_{n,k_m}\to\Flag_{n,\underline{k}}$.
We denote by $F_A$ the vector field on $\OO_{n,k}$ associated to the flow $\varphi^t_A$,
which is defined by
$F_A(X)=\frac{d}{dt}\left[ e^{t\,A}\ast X\right]_{t=0}$ for $X\in\OO_{n,k}$.
We have the following explicit formula:

\begin{prop} \label{FA:explicit} Given  $A\in\slin(n,\R)$ and $X\in\OO_{n,k}$,\;
$F_A(X) = A\,X-X\,(X^t\,A\,X)^\triangleleft$.
\end{prop}

\dem From remark 1' above we get
\begin{align*}
F_A(X) &= \frac{d}{dt}\left[ e^{t\,A}\ast X\right]_{t=0}
= \frac{d}{dt}\Pi\left( e^{t\,A}\, X\right)_{t=0} \\
&= \left[ D\Pi_{e^{t\,A}\, X}(A\,e^{t\,A}\, X)\right]_{t=0}
= D\Pi_{X}(A\, X)\\
&= A\,X-X\,(X^t\,A\,X)^\triangleleft\;.
\end{align*}
\cqd

Given $A\in\spl(n,\R)$, the flows $\varphi^t_A$ on $\OO_{2n,k}$
and $\Flag_{2n,\underline{k}}$ respectively induce by res\-triction flows
on the invariant submanifolds $\Osp_{n,k}$  and $\Fliso_{n,\underline{k}}$,
which we still denote by $\varphi_A^t$.
We also denote by $F_A$ the vector field on $\Osp_{n,k}$ associated with the flow 
$\varphi^t_A:\Osp_{n,k}\to \Osp_{n,k}$.

\begin{rmk}
Proposition~\ref{FA:explicit} holds with the same expression for
the induced flow $\varphi^t_A$ on the unitary frame manifold $\Osp_{n,k}$.
\end{rmk}

\bigskip

Assume now that $A\in\slin(n,\R)$ is a symmetric matrix.
\begin{prop} \label{FA:OO:eigen}  Given   $X\in\OO_{n,k}$,\;
$F_A(X) = 0$\; $\Leftrightarrow$\; every column of $X$ is an eigenvector of $A$.
\end{prop}

\dem
Assume every column of $X$ is an eigenvector of $A$.
Then if $D\in\Mat_{k\times k}(\R)$ is the diagonal matrix with the
corresponding eigenvalues \, $A\,X=X\,D$, and whence
$F_A(X)= A\,X-X\,(X^t\,A\,X)^\triangleleft = X\,D-X\,(X^t\,X\,D)^\triangleleft
=X\,D-X\,D^\triangleleft=0$, because $D^\triangleleft=D$.
Conversely, if $N=(X^t\,A\,X)^\triangleleft$ and $0=F_A(X) = A\,X-X\,N$,
it follows that $\langle X\rangle_{\{1,\ldots,i\}}$ is an invariant subspace under $A$,
for $i=1,\ldots, k$. Since $A$ is symmetric and the columns of $X$ are pairwise orthogonal, 
every column of $X$ must be an eigenvector of $A$.
\cqd

\begin{prop} \label{FA:Flag:fp}  Given  $\underline{k}\in\Lambda(n)$ and $V\in\Flag_{n,\underline{k}}$,\;
$V=\varphi^t_A(V)$ for some $t>0$\; $\Leftrightarrow$\; there is some matrix $X\in\OO_{n,k_m}$
such that $p(X)=V$ and $F_A(X)=0$.
\end{prop}

\dem
Assume $V=p(X)$ for some $X\in\OO_{n,k_m}$ with $F_A(X)=0$.
Then $X=\varphi^t_A(X)$ for all $t$, and since $p:\OO_{n,k_m}\to\Flag_{n,\underline{k}}$ 
semiconjugates the flows $\varphi^t_A$ on $\OO_{n,k_m}$ and $\Flag_{n,\underline{k}}$
it follows that $V=\varphi^t_A(V)$ for all $t$.
Conversely, assume that $V=\varphi^t_A(V)$ for some $t>0$.
Then each subspace $V_i$ in the flag $V$ is invariant under $A$.
Because $A$ is symmetric, the same is true about 
the subspaces $W_1=V1$ and $W_i=V_i\ominus V_{i-1}$, for $i=2,\ldots, m$.
Whence each $W_i$ is a direct sum of eigenspaces of $A$, and we can find
an orthonormal basis for $W_i$  formed by eigenvectors of $A$.
Putting these basis together, as columns of a matrix $X$, we have that
$V=p(X)$ and $F_A(X)=0$.
\cqd

\begin{rmk}
Propositions~\ref{FA:OO:eigen} and~\ref{FA:Flag:fp} also hold 
for isotropic flags and unitary frames provided the symmetric
matrix $A$ is chosen in the symplectic Lie algebra $\spl(n,\R)$.
\end{rmk}

\bigskip

\section{Riemannian Metrics}

Consider the Hilbert-Schmidt inner product on $\Mat_{n\times k}(\R)$
\[ \langle E,F\rangle_{{\rm hs}} = \frac{1}{k}\,\tr(E^t\,F) \]
with its associated norm
\[ \nrm{E}_{{\rm hs}}= \sqrt{ \frac{1}{k}\,\tr(E^t\,E) }  \;. \]
We take the Riemannian structures induced by this inner product on the
ortho frame and unitary frame manifolds $\OO_{n,k}$ and $\Osp_{n,k}$.
We shall refer to them as the {\em Hilbert-Schmidt} metrics on these manifolds.

\begin{lema} \label{hs:nrm}
Given $X\in \OO_{n,k}$, 
$E\in \Mat_{n\times k}(\R)$, $R\in\SO(n)$ and $S\in\OO(k)$,
\begin{enumerate}
\item $\nrm{X}_{{\rm hs}}=1$, 
\item $\nrm{R\,E\,S}_{{\rm hs}} = \nrm{E}_{{\rm hs}}$.
\end{enumerate}
\end{lema}

From this lemma, we see that both the left action of $\SO(n)$ and the right action of $\OO(k)$
on $\OO_{n,k}$ leave the Hilbert-Schmidt metric invariant, which means that both these actions are
isometric. Since $\Osp_{n,k}$ is a submanifold of $\OO_{2n,k}$ invariant under both 
the left action of $\Un(n)$ and the right action of $\OO(k)$, these actions too are isometric
for the Hilbert-Schmidt metric on $\Osp_{n,k}$.
We consider on $\Flag_{n,\underline{k}}$ and $\Fliso_{n,\underline{k}}$
the unique Riemannian metrics which respectively turn the projections
$p:\OO_{n,k_m}\to\Flag_{n,\underline{k}}$ and $p:\Osp_{n,k_m}\to\Fliso_{n,\underline{k}}$
to Riemannian submersions, i.e., such that the tangent maps
$Dp_X:T_X\OO_{n,k_m}\to T_{p(X)}\Flag_{n,\underline{k}}$
and $Dp_X:T_X\Osp_{n,k_m}\to T_{p(X)}\Fliso_{n,\underline{k}}$
are orthogonal projections. A linear map $A:E\to F$ between Euclidean spaces
$E$ and $F$ is said to be an orthogonal projection\, iff \,$\langle A\,X,A\,Y\rangle= \langle X,Y\rangle$
for every $X,Y\in \mbox{Ker}(A)^\perp$.

\bigskip

The following propositions describe the tangent-normal
space decompositions over the orthogonal and unitary frame manifolds.

\begin{prop}\label{TOrth:decomp}
Given $X\in\OO_{n,k}$,
\begin{enumerate}
\item[(a)] $T_X \OO_{n,k}=\{\, S\in\Mat_{n\times k}(\R)\,:\, S^t\,X+X^t\,S=0\,\} $,
\item[(b)] $T_X^\perp \OO_{n,k}=\{\, S\in\Mat_{n\times k}(\R)\,:\, S^t\,X-X^t\,S=0\;
\text{ and }\; (I-X\,X^t)\,S=0\;\} $.
\end{enumerate}
The corresponding projections\,
$\Pi^T:\Mat_{n\times k}(\R)\to T_X \OO_{n,k}$\,
and\, $\Pi^\perp:\Mat_{n\times k}(\R)\to T_X^\perp \OO_{n,k}$\,
are given by :
\begin{enumerate}
\item[(a')] $\Pi^T(B)= \frac{1}{2}\, X\,(X^t\,B-B^t\, X)+ (I-X\,X^t)\,B$,
\item[(b')] $\Pi^\perp(B)= \frac{1}{2}\, X\,(X^t\,B+B^t\, X)$.
\end{enumerate}
\end{prop}

\dem
Given $X\in\OO_{n,k}$ define
\begin{align*}
E^+_{n,k}(X) &= \{\, S\in\Mat_{n\times k}(\R)\,:\, S^t\,X+X^t\,S=0\,\} \;,\\
E^-_{n,k}(X)&= \{\, S\in\Mat_{n\times k}(\R)\,:\, S^t\,X-X^t\,S=0\,\} \;.
\end{align*}
These two spaces are transversal, i.e.,
$\Mat_{n\times k}(\R) = E^+_{n,k}(X) + E^-_{n,k}(X)$, 
and their intersection is
$ E^0_{n,k}(X) = E^+_{n,k}(X)\cap E^-_{n,k}(X) =
\{\, S\in\Mat_{n\times k}(\R)\,:\, X\,X^t\,S=0\,\} $.
Because $X\,X^t$ is the orthogonal projection onto the linear subspace
of $\R^n$ spanned by the columns of $X$, it follows that the orthogonal
complement of $ E^0_{n,k}(X)$ is
$ E^0_{n,k}(X)^\perp =
\{\, S\in\Mat_{n\times k}(\R)\,:\, (I-X\,X^t)\,S=0\,\} $. This implies that
\begin{equation}\label{Epm:orth:compl}
 E^\pm_{n,k}(X)^\perp = E^\mp_{n,k}(X)\cap E^0_{n,k}(X)^\perp \;. 
\end{equation}
Now, it is obvious that (a)\, $T_X \OO_{n,k}=E^+_{n,k}(X)$, and
then~(\ref{Epm:orth:compl}) implies (b).
It is also clear that $B=\Pi^T(B)+\Pi^\perp(B)$.
Finally, one can easily check that $R=\Pi^T(B)$ defined in (a') satisfies
$R^t\,X=-X^t\,R$, and that $S=\Pi^\perp(B)$
defined in (b') satisfies $S^t\,X=X^t\,S$, $(I-X\,X^t)\,S=0$.
\cqd

\bigskip

In the symplectic case we have

\begin{prop}\label{TUnitary:decomp}
Given $X\in\Osp_{n,k}$,
\begin{enumerate}
\item[(a)] $T_X \Osp_{n,k}=\{\, S\in\Mat_{2n\times k}(\R)\,:\, S^t X+X^t S=0\;
\text{ and }\; S^t J X+X^t J S=0 \,\} $,
\item[(b)] $T_X^\perp \Osp_{n,k}$ is the sum of two subspaces:
\begin{enumerate}
	\item[i.] $\{\, S\in\Mat_{2n\times k}(\R)\,:\, S^t X-X^t S=0\;
\text{ and }\; (I-X X^t) S=0\;\} $, 
	\item[ii.] $\{\, S\in\Mat_{2n\times k}(\R)\,:\, S^t J X-X^t J S=0\;
\text{ and }\; (I+ J X X^t J) S=0\;\} $.
\end{enumerate}
\end{enumerate}
The corresponding projections\,
$\Pi^T:\Mat_{2n\times k}(\R)\to T_X \Osp_{n,k}$\,
and\, $\Pi^\perp:\Mat_{2n\times k}(\R)\to T_X^\perp \Osp_{n,k}$\,
are given by :
\begin{enumerate}
\item[(a')] $\Pi^T(B)= \frac{1}{2}  X ( X^t B - B^t X ) +
	\frac{1}{2}  J X ( B^t J X - X^t J B)  + 
	(I - X X^t + J X X^t J)\,B$,
\item[(b')] $\Pi^\perp(B)= \frac{1}{2} X (X^t B + B^t X) - 
	\frac{1}{2} J X ( X^t J B + B^t J X )$.
\end{enumerate}
\end{prop}

\dem
Differentiating the two defining equations of $\Osp_{n,k}$
we obtain that $T_X \Osp_{n,k}=E^+_{n,k}(X)\cap E^-_{n,k}(J X)$,
which proves (a). Therefore, by~(\ref{Epm:orth:compl}),
we get that $T_X  {\Osp_{n,k}}^\perp$ is the sum of the subspaces
$i.$ and $ii.$ above
\begin{align*}
T_X {\Osp_{n,k}}^\perp  &= (\, E^+_{n,k}(X)\cap E^-_{n,k}(J X)\,)^\perp 
= E^+_{n,k}(X)^\perp +  E^-_{n,k}(J X)^\perp \\
&=   \underbrace{ E^-_{n,k}(X)\cap E^0_{n,k}(X)^\perp }_{ {i.}} \; +\;
\underbrace{  E^+_{n,k}(J X)\cap E^0_{n,k}(J X)^\perp }_{ {ii.}} \;.
\end{align*}

One can easily check that $B=\Pi^T(B)+\Pi^\perp(B)$.
To finish the proof we need, as in the previous lemma, to show
that $R=\Pi^T(B)$ defined in (a') satisfies
$R^t X+X^t R=0$ and $R^t J X+X^t J R=0$,  and also to show that writing
$S=\Pi^\perp(B)$ defined in (b') as $S=S_1+S_2$ with
$S_1=\frac{1}{2} X (X^t B + B^t X)$, $S_2=-\frac{1}{2} J X ( X^t J B + B^t J X )$,
then  $S_1$  belongs to space $i.$ while
$S_2$  belongs to  space $ii.$. Because this argument does not shade any light on
formulas (a') and (b'),  we describe now a more geometric approach on how to derive 
these projection formulas. First we need a couple of definitions. 
We denote by $P_E$ the orthogonal projection onto a subspace $E$ of some euclidean space.
Let us  say that two subspaces $E$ and $F$ are {\em perpendicular}\, iff\, the inner product of any vector
$u\in (E\cap F)^\perp\cap E$ with any vector $v\in (E\cap F)^\perp\cap F$ is always zero.
When this occurs, the following general formulas hold
$$ P_{E\cap F} = P_E\circ P_F \quad \text{ and  }\quad
P_{(E\cap F)^\perp} = P_{F^\perp} + P_{E^\perp}\circ P_F\;. $$
Finally, formulas (a') and (b') can be driven from these abstract formulas
applied to the subspaces
$E=E^-_{n,k}(J X)$ and $F=E^+_{n,k}(X)$, showing first they are perpendicular.
\cqd

\bigskip

\section{Quadratic Functions}

Let $\Delta_k=\{\,b=(b_1,\ldots, b_k)\in\R^k\,:\,
b_1\geq b_2\geq \ldots \geq b_k >0 \,\}$, and for each $b\in\Delta_k$
denote by $D_b$ the diagonal matrix  
$$ D_b = \left(
\begin{array}{cccc}
b_1 & 0 & \cdots & 0 \\
0 & b_2 & \cdots & 0 \\
\vdots & \vdots & \ddots & \vdots \\
0 & 0 & \cdots & b_k
\end{array}
\right) \;. $$
\begin{defi}
Given a symmetric matrix $A\in\slin(n,\R)$ and $b\in\Delta^k$, we define the quadratic functions 
$Q_{A,b}:\OO_{n,k}\to\R$ and $Q_{A}:\OO_{n,k}\to\R$ respectively by \,
$$Q_{A,b}(X)= \langle A\,X\,D_b,X\,D_b\rangle_{{\rm hs}}\quad\text{ and }\quad
Q_{A}(X)= \langle A\,X,X\rangle_{{\rm hs}}\;.$$
\end{defi}
For each $\underline{k}\in\Lambda(n)$ \, set
\begin{align*}
\Delta_{\underline{k}}&=\{\, b\in\Delta_{k_m}\,:\,
b_1=\cdots = b_{k_1}\,<\, b_{k_1+1}=\cdots = b_{k_2}\,<\,\ldots \,< \,b_{k_{m-1}+1} = \cdots = b_{k_m}\,\}\\
\overline{\Delta}_{\underline{k}}&=\{\, b\in\Delta_{k_m}\,:\,
b_1=\cdots = b_{k_1}\,\leq\, b_{k_1+1}=\cdots = b_{k_2}\,\leq\,\ldots \,\leq \,b_{k_{m-1}+1} = \cdots = b_{k_m}\,\} 
\end{align*}
Note that $Q_{A}=Q_{A,b}$ with $b=(1,\ldots, 1)$. Also, given $k\geq 1$, the family $\{\Delta_{\underline{k}}\}_{\underline{k}}$ \, 
indexed over all $\underline{k}\in\Lambda(n)$
with last entry $k_m=k$\, is a partition of $\Delta_k$.

\bigskip

\begin{prop}
Given a symmetric matrix $A\in\slin(n,\R)$, $\underline{k}\in\Lambda(n)$  and $b\in\overline{\Delta}_{\underline{k}}$,
the function $Q_{A,b}:\OO_{n,k_m}\to\R$ is invariant under the right action of
$\OD_{\underline{k}}$.
\end{prop}

\dem
Use lemma~\ref{hs:nrm} and the fact that for $b\in\overline{\Delta}_{\underline{k}}$  
the diagonal matrix $D_b$ commutes with every matrix in $\OD_{\underline{k}}$.
\cqd

This proposition shows that given a symmetric matrix $A\in\slin(n,\R)$, $\underline{k}\in\Lambda(n)$  and $b\in\overline{\Delta}_{\underline{k}}$, 
$Q_{A,b}:\OO_{n,k_m}\to\R$ induces a quotient function $Q_{A,b}:\Flag_{n,\underline{k}}\to\R$.
For symplectic matrices $A\in\spl(n,\R)$ we denote by
$Q_{A,b}:\Osp_{n,k_m}\to\R$ and $Q_{A,b}:\Fliso_{n,\underline{k}}\to\R$ the restriction
functions. Note that $\Osp_{n,k_m}$ and $\Fliso_{n,\underline{k}}$ are submanifolds respectively of
$\OO_{2n,k_m}$ and $\Flag_{2n,\underline{k}}$.

\bigskip

\begin{prop}\label{QAb:formula}
Given $A\in  \slin(n,\R)$, resp. $A\in\spl(n,\R)$, and $b=(b_1,\ldots, b_k)\in\Delta_k$,
if we write  $\omega_i = i\,(b_i^2-b_{i+1}^2)/k >0$ \, for $i=1,\ldots, k-1$
\, and $\omega_k=b_k^2$,\, then the function
$Q_{A,b}:\OO_{n,k}\to\R$, resp. $Q_{A,b}:\Osp_{n,k}\to\R$, is given by
\begin{equation}\label{Qab:decomp}
Q_{A,b}(X)=  \sum_{i=1}^k \omega_i \,(Q_A\circ \pi_i)(X) \;,
\end{equation}
where each $\pi_i:\OO_{n.k}\to\OO_{n,i}$, resp. $\pi_i:\Osp_{n.k}\to\Osp_{n,i}$,
denotes the canonical projection $\pi_i(X)=[X]_i$.
\end{prop}

\dem
With the convention that  $b_{k+1}= 0$,
we have
\begin{align*}
Q_{A,b}(X) &= \langle A\,X\,D_b, X\,D_b \rangle_{{\rm hs}} 
= \frac{1}{k}\sum_{i=1}^k   (A\,X\,D_b)_i \cdot  (X\,D_b)_i  \\
&= \frac{1}{k}\sum_{i=1}^k b_i^2  (A\,X)_i \cdot  (X)_i  \\
&= \frac{1}{k}\sum_{i=1}^k (b_i^2-b_{i+1}^2)\,\left( 
 (A\,X)_1\cdot (X)_1  +\cdots + (A\,X)_i\cdot (X)_i \,\right)\\
&= \sum_{i=1}^k \frac{i}{k} (b_i^2-b_{i+1}^2) \,\langle [A\,X]_i, [X]_i \rangle_{{\rm hs}}
= \sum_{i=1}^k \frac{i}{k} (b_i^2-b_{i+1}^2) \,\langle A\,[X]_i, [X]_i \rangle_{{\rm hs}} \\
&= \sum_{i=1}^k \omega_i \,Q_A(\,\pi_i(X)\,)\;.
\end{align*}
\cqd

If $b\in\Delta_{\underline{k}}$ with $\underline{k}=(k_1,\ldots, k_m)$
then $\omega_{k_i}>0$ for $i=1,\ldots, m$, and $\omega_j=0$ for all $j\neq k_i$.
This shows that $Q_{A,b}:\Flag_{n,\underline{k}}\to\R$ is a sum of pullbacks of functions
defined over the Grassmannian manifolds $\Grass_{n,k_i}=\Flag_{n,(k_i)}\preceq \Flag_{n,\underline{k}}$.

\begin{lema}\label{QA:grads}
The gradient of $Q_{A}:\OO_{n,k}\to\R$ at $X$ is given by
$$ \bigtriangledown (Q_{A})(X) = 2\,(I-X\,X^t)\,A\, X\;. $$
\end{lema}

\dem
Consider the natural quadratic extension of $Q_A$ to the linear space of all matrices:\,
$Q_A(X)=\tr (X^t\,A \,X)$.
The gradient of this extension at $X$ is  $2\,A X$.
Then by proposition~\ref{TOrth:decomp},
\begin{align*}
\bigtriangledown (Q_{A})(X) &= \Pi^T(\, 2\,A X\, ) \\
&= X\,\underbrace{ (X^t\,A \,X-(A \,X)^t\, X)}_{=0} + 2\,(I-X\,X^t)\,A \,X \\
&=  2\,(I-X\,X^t)\,A \, X\;.
\end{align*}
In the symplectic case we compute the following expression for the gradient
$$ \bigtriangledown (Q_{A})(X) = 2\, (J X) (J X)^t A  X  + 2\,(I-X\,X^t+J X X^t J)\,A \, X $$
which simplified gives exactly the same expression \, $2\,(I-X X^t) A X$.
\cqd

\begin{rmk}
Let $P^X$ be the orthogonal projection onto $\langle X\rangle_{\{1,\ldots, k\}}$. 
Then we  can write \,
$ \bigtriangledown (Q_{A})(X) = 2\,(I-P^X)\,A \, X $.
\end{rmk}

\bigskip

\bigskip

Next we derive an explicit formula for the gradient of
$Q_{A,b}$. Given $X\in\OO_{n,k}$,
we denote by $P^X_i=P_i$ ($1\leq i\leq k+1$)
the following orthogonal projection matrices:
For $k\leq i$, $P_i$ is the projection onto the linear span $\langle X \rangle_{\{i\}}$, 
while for $i=k+1$,
$P_{k+1}$ is the orthogonal projection onto $(\langle X\rangle_{\{1,\ldots, k\}})^\perp$.
Notice that\,
$P_1+\ldots + P_{k+1}=I$.
We write as above $\omega_i = i\,(b_i^2-b_{i+1}^2)/k >0$.
Notice that for $i\leq j$ we have
$\frac{k}{i}\,\omega_i+\ldots +\frac{k}{j}\,\omega_j = b_i^2- b_{j+1}^2$.
Let $E_{i}(b)\in\Mat_{k\times k}(\R)$ denote the matrix
\begin{equation}\label{Ei:b}
 E_{i}(b) = \left(
\begin{array}{ccccccc}
b_1^2-b_{i+1}^2 & 0 & \cdots & 0 &  0 & \cdots & 0 \\
0 & b_2^2-b_{i+1}^2 & \cdots & 0 &  0 & \cdots & 0\\
\vdots & \vdots & \ddots & \vdots &  \vdots & \ddots & \vdots\\
0 & 0 & \cdots & b_i^2-b_{i+1}^2 &  0 & \cdots & 0\\
0 & 0 & \cdots & 0 &  0 & \cdots & 0\\
\vdots & \vdots & \ddots & 0 &  \vdots & \ddots & \vdots\\
0 & 0 & \cdots & 0 &  0 & \cdots & 0
\end{array}
\right)\;, 
\end{equation}
and let $E_{i} \in\Mat_{k\times k}(\R)$ be the matrix representing the orthogonal projection
onto the linear span $\langle e_1,\ldots, e_{i}\rangle$.

\bigskip

\begin{prop}\label{grad:QAb:decomp}
The gradient of $Q_{A,b}:\OO_{n,k}\to\R$ at $X$ is given by
\[ \bigtriangledown (Q_{A,b})(X) = 2\,\sum_{i=1}^{k+1} P^X_i\, A \,X \, E_{i-1}(b)\;.
\]
\end{prop}

\dem 
Notice that $I_{k, i}\,I_{i\, k}=E_i$.
The projection $\pi_i$ can be written as $\pi_i(X)=X\, I_{k, i}$. Its derivative
$(D\pi_i)_X\, V= V \, I_{k, i}$ has adjoint 
$(D\pi_i)_X^\ast\, W= \frac{k}{i}\, W \, I_{i, k }$.
Therefore, using lemma~\ref{QA:grads}, we get that
\begin{align*}
 \bigtriangledown (Q_A\circ \pi_i)(X)  &=
( D\pi_i)_X^\ast (\bigtriangledown  Q_A)( \pi_i(X) ) \\
 &= \frac{k}{i}\, (\bigtriangledown  Q_A)( X I_{k, i} )\, I_{i\times k} \\
&= \frac{2k}{i}\,\left( I-(P_1+\cdots +P_i) \right)\,\,A \, X\,E_i\\
&= \frac{2k}{i}\,\left( P_{i+1}+\cdots + P_{k+1}  \right)\,A \, X\,E_i \;.
\end{align*}
Thus
\begin{align*}
\bigtriangledown (Q_{A,b})(X) &=
\sum_{i=1}^k \omega_i \,\bigtriangledown (Q_A\circ \pi_i)(X)\\
&= 2\,
\sum_{i=1}^k \frac{k}{i}\, \omega_i\,  \left( \sum_{j=i+1}^{k+1}  P_{j} \right) \,A \, X\,E_i\\
&= 2\,
\sum_{j=2}^{k+1}    P_{j}\,A \, X\,\left( \sum_{i=1}^{j-1} \frac{k}{i}\, \omega_i\,E_i \right)\\
&= 2\,
\sum_{j=2}^{k+1}    P_{j}\,A \, X\, E_{j-1}(b)\;.
\end{align*}
\cqd

\bigskip

\begin{prop}\label{crit:point:char}
Given $A\in  \gfrak^{\rm sym}_{n}$,  $b_1> b_2>\ldots > b_k> 0$, and $X\in\OO_{n,k}$,\,
$X$ is a critical point of $Q_{A,b}$ \, iff\, 
every column of $X$ is an eigenvector of $A$.
\end{prop}

\dem
It is clear that if each column of $X$ is an eigenvector of $A$
then $X$ is a critical point of $Q_A\circ \pi_i$, for every $i=1,\ldots, k$.
Therefore, $X$ is a critical point of $Q_{A,b}$.
Assume now that $\bigtriangledown (Q_{A,b})(X) = 0$.
Fix $j=1,\ldots, k$.
Since
\[
0=\frac{1}{2}\,\bigtriangledown (Q_{A,b})(X)\,e_j  = 
 \sum_{i=j+1}^{k+1} (b_j^2-b_i^2)\,P_i\, A \,X \, e_j =
\sum_{i=j+1}^{k+1} (\overbrace{b_j^2-b_i^2}^{>0})\,P_i\, A \,(X)_j\;,
\]
it follows that $P_i\, A \,(X)_j=0$ for every $i=j+1,\ldots, k+1$,
which implies that $A \,(X)_j\in \langle [X]_j \rangle$.
By induction we derive (see proof of lemma~\ref{fixed:point:char})
that $(X)_j$ is an eigenvector of  $A$.
\cqd

\bigskip

\section{Lyapunov functions}

We shall denote by $G^{\rm sym}_{n,+}$ the space of symmetric positive definite matrices
in $G_n$ and by $\gfrak^{\rm sym}_{n}$ the space of symmetric matrices
in $\gfrak_n$. We write  $\SL^{\rm sym}_{n,+}$,  $\Sp^{\rm sym}_{n,+}$, 
$\slin^{\rm sym}_{n}$ or  $\spl^{\rm sym}_{n}$ to
emphasize that $G_n=\SL(n,\R)$, $G_n=\Sp(n,\R)$, 
$\gfrak_n=\slin(n,\R)$ or $\gfrak_n=\spl(n,\R)$.

To establish Theorem A it is enough proving that the lift 
$Q_{A,b}:\OO_{n,k}\to\R$ is a Lyapunov function for lifted flow $\varphi_H^t:\OO_{n,k}\to\OO_{n,k}$
when $A,H\in\gfrak_n$ are symmetric matrices sharing the same ordered
eigen-directions.

\bigskip

\begin{prop}\label{AX:R}
Given $B\in\SL(n,\R)$ and $X\in\OO_{n,k}$, let 
$B\,X=R\,D\,R'$ be the singular value decomposition of $B\,X$, 
with $R\in\OO_{n,k}$, $R'\in\OO_{k}$
orthogonal, and $D$ a diagonal positive  $k\times k$ matrix.
Then there is an orthogonal matrix $S\in\OO_{k}$ such that\;  $B\ast X = R\,S$.
\end{prop}

\dem
Let $B\,X=K\,U$ be the QR-decomposition of $B\,X$.
We have $B\ast X = K=R\,S$ with $S=D\,R'\,U^{-1}\in\SL(k,\R)$.
Since $K$ and $R$ are both orthogonal, it follows
\[ I=K^t\,K = (R\,S)^t\,(R\,S)= S^t\,(R^t\,R)\,S = S^t\,S\;,\]
which proves that $S\in\OO_{k}$.
\cqd

\begin{rmk}\label{nrm:AX:R}
With the notation of the previous proposition it follows that
\[ Q_A(B\ast X) = \langle A\,R, R\rangle_{{\rm hs}} \;.\]
\end{rmk}

\bigskip

\begin{lema} \label{trEF:E,F>=0}
Given matrices $E, F\in\Mat_{n\times n}(\R)$,
if $E$ and $F$ are symmetric positive semi-definite, i.e., $E,F\geq 0$, then
$\tr(E\,F)\geq 0$, with
$\tr(E\,F)=0$\, iff\, $E\,F=0$.
\end{lema}

\dem
We can write $E=R_1^t\, D_1\,R_1$ and  $F=R_2^t\, D_2\,R_2$ 
with $D_i$ diagonal positive semi-definite and $R_i$ orthogonal,
for $i=1,2$.
Then
\begin{align*}
\tr(E\,F) &= \tr ( \, (R_1^t\, D_1\,R_1)\,(R_2^t\, D_2\,R_2)\,) \\
&= \tr ( \, \sqrt{D_1}\,R_1\,R_2^t\, D_2\,R_2\,R_1^t\, \sqrt{D_1}\,) \\
&= \tr ( \, \sqrt{D_1}\,(R_2\,R_1^t)^t\, D_2\,(R_2\,R_1^t)\, \sqrt{D_1}\,) \\
&= n\,\nrm{ \sqrt{D_2}\,\,(R_2\,R_1^t)\, \sqrt{D_1} }^2_{{\rm hs}} \geq 0\;.
\end{align*}
Since the matrices $E\,F$ and $\sqrt{D_2}\,\,(R_2\,R_1^t)\, \sqrt{D_1}$ are
conjugated,
$\tr(E\,F)=0$ implies $\sqrt{D_2}\,\,(R_2\,R_1^t)\, \sqrt{D_1}=0$,
which in turn forces $E\,F=0$. 
\cqd

\bigskip

\begin{lema}\label{fixed:point:char}
Given $B\in G^{\rm sym}_{n,+}$  and $X\in\OO_{n,k}$,
$B\ast X=X$\, iff\, every column of $X$ is an eigenvector of $B$.
\end{lema}

\dem
Assume every column of $X$ is an eigenvector of $B$.
Then there is some positive diagonal matrix $D$ such that $B\,X=X\,D$,
which implies that $B\ast X=\Pi(B\,X)=\Pi(X\,D)=X$.
Conversely, if $B\ast X=\Pi(B\,X)=X$ there is some
 $U\in\UTpos(k,\R)$  such that $B\,X=X\,U$.
Denote by $X_i$ the $i^{\text{th}}$-column of $X$ and let
$V_i=\langle X_1,\ldots, X_i\rangle$ be the linear span of the first $i$
columns of $X$.
Then $\{0\}=V_0\subset V_1\subset V_2\subset \ldots \subset V_k$ is a flag invariant by $B$.
Since $X_i\in V_i\cap V_{i-1}^\perp$, for each $i=1,\ldots, k$, it follows that
$X_i$ must be an eigenvector of $B$.
\cqd

\begin{rmk}\label{QA:XS:old}
 Let $B\in G^{\rm sym}_{n,+}$, $A\in\gfrak^{\rm sym}_n$, $X\in\OO_{n,k}$ and $S\in\OO_{k}$.
\begin{enumerate}
\item If $B\ast X=X$ then  in general $B\ast (X\,S)\neq X\,S$, but
\item if $Q_A(B\ast X)=Q_A(X)$ then $Q_A(B\ast (X\,S) )=Q_A(X\,S)$.
\end{enumerate}
\end{rmk}

\bigskip


Let $P\in\gl(n,\R)$ be an orthogonal projection matrix, i.e.,
$P=P^t =P^2$.
We define a bilinear form
$\xi_P:\gfrak_n^{\rm sym}\times \gfrak_n^{\rm sym}\to \R$ by
$\xi_P(A,H)=\tr( P\,A\,(I-P)\,H)$.

\begin{prop}\label{xiP(A;H)}
Given symmetric matrices $A,H\in\gfrak_n$  sharing the same ordered
eigen-directions,\, $\xi_P(A,H)\geq 0$,
and $\xi_P(A,H)=0$\, $\Leftrightarrow$\, $A\,P=P\,A$ and $H\,P=P\,H$.
\end{prop}

\dem
We can with no loss of generality assume that $A$ and $H$ are diagonal
matrices.  Denote by $D_x$ the diagonal matrix with diagonal
entries $x=(x_1,\ldots, x_n)\in\R^n$.
Then the bilinear form $\xi_P$ can be seen as the 
bilinear form $\xi_P:\R^n\times \R^n\to\R$ defined by
$\xi_P(x,y)= \tr( P\,D_x\,(I-P)\,D_y)$.
Consider the cone $\Gamma=\{\,x\in\R^n\,: \, x_1\geq \ldots \geq x_n\,\}$.
We want to see that given $a,h\in\Gamma$, $\xi_P(a,h)\geq 0$,
and $\xi_P(a,h)=0$\, $\Leftrightarrow$\, $D_a\,P=P\,D_a$ and $D_h\,P=P\,D_h$.
We define the vectors
$v_i=(1,\ldots,1,0,\ldots, 0)\in\Gamma$  with the first $i$ coordinates equal to $1$
for $1\leq i \leq n$, which form a basis of $\R^n$.
Notice that because $D_{v_n}=I$ is the identity matrix we get that
$\xi_P(v_n,x)=0=\xi_P(x,v_n)$ for every $x\in\R^n$.
If $i\leq j\leq n-1$
\begin{align*}
\xi_P(v_i,v_j) &= \sum_{k=1}^i p_{k,k} -\sum_{k=1}^i \sum_{r=1}^j (p_{k,r})^2 =
\sum_{k=1}^i \left( p_{k,k}-\sum_{r=1}^j (p_{k,r})^2\right) \\
&=\sum_{k=1}^i  \sum_{r=j+1}^n (p_{k,r})^2\geq 0 \;,
\end{align*}
and otherwise, if $j < i \leq n-1$,
\begin{align*}
\xi_P(v_i,v_j) &= \sum_{k=1}^j p_{k,k} -\sum_{r=1}^i \sum_{k=1}^j (p_{r,k})^2 =
\sum_{k=1}^j \left( p_{k,k}-\sum_{r=1}^j (p_{r,k})^2\right) \\
&=\sum_{k=1}^j  \sum_{r=i+1}^n (p_{r,k})^2 \geq 0\;.
\end{align*}
Now, given $a,h\in \Gamma$ we can write
$$ a= a_n \,v_n + \sum_{i=1}^{n-1} (a_i-a_{i+1})\,v_i \quad
\text{ and } \quad
h= h_n \,v_n + \sum_{i=1}^{n-1} (h_i-h_{i+1})\,v_i\;. $$
Whence
$$\xi_P(a,h)=\sum_{i,j=1}^{n-1} \underbrace{(a_i-a_{i+1})}_{\geq 0}\,
\underbrace{(h_j-h_{j+1})}_{\geq 0}\,\underbrace{\xi_P(v_i,v_j)}_{\geq 0} \geq 0\;.$$
Assume now that $\xi_P(a,h)=0$, and $a_i>a_j$ $\Leftrightarrow$ $h_i>h_j$,
for every $1\leq i<j\leq n$.
Then for every $i=1,\ldots, n$ such that $a_i>a_{i+1}$ we have
$0=\xi_P(v_i,v_i)= \sum_{k=1}^i  \sum_{r=i+1}^n (p_{r,k})^2$,
which implies that $p_{r,k}=p_{k,r}=0$ for all $k\leq i$ and $r>i$.
This shows that $P$ is a block-diagonal matrix whose blocks
correspond in $A=D_a$ and $H=D_h$ to multiples of identity matrices.
Therefore $P$ commutes with $A=D_a$ and $H=D_h$.
\cqd

\begin{lema}\label{tr:B:A:triangleleft}
Given $A,B\in\gl(k,\R)$, if $A^t=A$ and $B^t=B$\, then\,
$\tr( B\,A^\triangleleft ) = \tr (B\,A)$.
\end{lema}

\dem
Let $S=A^\triangleleft-A$.
The entriy $s_{i,j}$ of $S$ is\,
$s_{i,j}= a_{i,j}+a_{j,i} - a_{i,j} = a_{i,j}$ when $i <j$,
$s_{i,j}= a_{i,i} - a_{i,i} = 0$ when $i = j$, and
$s_{i,j}= 0 - a_{i,j} = -a_{i,j}$ if $i > j$.
Therefore $S$ is skew-symmetric and
$$ \tr(B\,S) = \tr (S\,B) = \tr( -S^t\,B^t) = -\tr( (B\,S)^t) 
= -\tr(B\,S)\;, $$
which implies $\tr(B\,(A^\triangleleft-A))=\tr(B\,S)=0$
and proves the lemma.
\cqd

\begin{lema}\label{sing:char}
Given $H\in \gfrak^{\rm sym}_{n}$  and $X\in\OO_{n,k}$,
$F_H(X)=0$\, iff\, every column of $X$ is an eigenvector of $H$.
\end{lema}

\dem
Assume every column of $X$ is an eigenvector of $H$.
Then there is some positive diagonal matrix $D$ such that $H\,X=X\,D$,
which implies that $$F_H(X)=
H\,X-X\,(X^t\,H\,X)^\triangleleft =
X\,D- X\,(X^t\,X\,D)^\triangleleft =
X\,D- X\,D^\triangleleft = 0\;.$$
Conversely, assume $0=F_H(X)=H\,X-X\,U$ with
 $U = (X^t\,H\,X)^\triangleleft \in\UTpos(k,\R)$, so  that $H\,X=X\,U$.
Denote by $X_i$ the $i^{\text{th}}$-column of $X$ and let
$V_i=\langle X_1,\ldots, X_i\rangle$ be the linear span of the first $i$
columns of $X$.
Then $\{0\}=V_0\subset V_1\subset V_2\subset \ldots \subset V_k$ is a flag invariant by $H$.
Since $X_i\in V_i\cap V_{i-1}^\perp$, for each $i=1,\ldots, k$, it follows that
$X_i$ must be an eigenvector of $H$.
\cqd

\begin{rmk}\label{QA:XS}
Given $A,H\in\gfrak^{\rm sym}_n$, $X\in\OO_{n,k}$ and $S\in\OO_{k}$.
\begin{enumerate}
\item If $F_H(X)=0$ then  in general $F_H(X\,S)\neq 0$, but
\item if $Q_A(\varphi_H^t(X))=Q_A(X)$ then $Q_A(\varphi_H^t(X\,S) )=Q_A(X\,S)$.
\end{enumerate}
\end{rmk}

\bigskip

\begin{prop}\label{QA:Lyap:varphi_H^t}
Given $A,H \in \gfrak^{\rm sym}_{n}$ sharing the same ordered
eigen-directions, we have for every $X\in\OO_{n,k}$ and $t>0$,
\[ Q_A(\varphi_H^t(X)) \geq Q_A(X) \;,\]
with equality \, iff\, for some $S\in\OO_{k}$,\, $F_H(X\,S)= 0$.
\end{prop}

\dem
Given $X\in\OO_{n,k}$, using lemma~\ref{tr:B:A:triangleleft} and proposition~\ref{xiP(A;H)}
we get
\begin{align*}
D (Q_A)_X\,F_H(X) &= 2\,\langle A\,X, F_H(X) \rangle_{{\rm hs}}\\
&= 2\,k^{-1}\,\tr ( X^t\,A\, F_H(X)  )\\
&= 2\,k^{-1}\,\tr ( X^t\,A\,(H\,X-X\,(X^t\,H\,X)^\triangleleft)   )\\
&= 2\,k^{-1}\,\tr ( X^t\,A\,H\,X) -
2\,k^{-1}\,\tr ( (X^t\,A\,X)\,(X^t\,H\,X)^\triangleleft)   )\\
&= 2\,k^{-1}\,\tr ( X^t\,A\,H\,X) - 2\,k^{-1}\,\tr (  X^t\,A\,X \, X^t\,H\,X   )\\
&= 2\,k^{-1}\,\tr (  X^t\,A\,(I-X \, X^t)\,H\,X   )\\
&= 2\,k^{-1}\,\tr (  X\,X^t\,A\,(I-X \, X^t)\,H   )\\
&= 2\,k^{-1}\,\tr (  P\,A\,(I-P)\,H   )=  2\,k^{-1}\,\xi_P(A,H)\geq 0 \;,
\end{align*}
where $P=X\,X^t$ is the orthogonal projection onto the subspace $\langle X\rangle$ spanned by
the columns of $X$.
Whence $Q_A$ is a Lyapunov function for the flow $\varphi_H^t$ of $F_H$.
To see it is a strict Lyapunov function assume 
$Q_A(\varphi_H^t(X)) = Q_A(X)$ for some $t>0$.
Then $D (Q_A)_X\,F_H(X)=0$ which by proposition~\ref{xiP(A;H)}
implies that $H\,P= P\,H$.
This commutativity shows that the subspace $\langle X\rangle$ is $H$-invariant.
Therefore, there is an orthogonal matrix $S\in\OO_k$ such that every
column of $X\,S$ is an eigenvector of $H$, and by lemma~\ref{sing:char} it follows that $F_H(X\,S)=0$.
\cqd

\begin{rmk} \label{equality:invariance} The argument above also proves that
$Q_A(\varphi_H^t(X))= Q_A(X)$ \, iff\, $\langle X\rangle$
is an $H$-invariant subspace.
\end{rmk}

\bigskip

\begin{coro} Given $A,H \in \gfrak^{\rm sym}_{n}$ sharing the same ordered
eigen-directions, the function 
$Q_A:\OO_{n,k}\to\R$ is a Lyapunov function for 
$\varphi_A:\OO_{n,k}\to\OO_{n,k}$, but it is only a
strict Lyapunov function when $k=1$.
\end{coro}

\bigskip

\begin{teor}
Given $A,H \in \gfrak^{\rm sym}_{n}$ sharing the same ordered
eigen-directions and  $b_1>b_2>\ldots >b_k>0$, 
the function $Q_{A,b}:\OO_{n,k}\to\R$ is a strict Lyapunov function for
the map $\varphi_A:\OO_{n,k}\to\OO_{n,k}$.
\end{teor}

\dem
By propositions~\ref{projs:equivariance},~\ref{QA:Lyap:varphi_H^t} and~\ref{QAb:formula}, 
$Q_{A,b}$ is a sum of $k$ Lyapunov functions
for the flow $\varphi_H^t$,
and, therefore, it is also a Lyapunov function of $\varphi_H^t$.
Assume now that $Q_{A,b}(\,\varphi_H^t(X)\,)=Q_{A,b}(X)$ for some $t>0$.
Then 
$Q_A(\,\pi_i\circ\varphi_H^t(X) \,)=Q_A(\,\pi_i(X)\,)$,
 for every $i=1,\ldots, k$.
By remark~\ref{equality:invariance} this implies that
$\langle [X]_i\rangle$ is $H$-invariant for all $i=1,\ldots, k$.
By induction we can easily prove that each column of $X$ is an eigenvector of $H$.
Therefore, by lemma~\ref{sing:char}, $F_H(X)=0$.
\cqd

\bigskip

\section{The Stratifications}
\label{the:stratifications}

Fix a commutative positive symmetric subgroup $H_n\subset\SL(n,\R)$,
and let $E_1, E_2,\ldots, E_n$ be its eigen-directions.
Each $E_i$ is a one dimensional subspace such that $A\,E_i=E_i$
for all $A\in H_n$, and we have $\R^n=E_1\oplus E_2\oplus \ldots \oplus E_n$.
We also assume that a commutative positive symmetric subgroup $H_n\subset\Sp(n,\R)$ is fixed.
In this case $H_n$ has $2n$  eigen-directions denoted by $E_1, E_2,\ldots, E_{2n}$, and 
we assume the directions $E_i$ and $E_{i+n}$ to be conjugate
in the sense that $J E_{i} = E_{i+n}$, for every $i=1,\ldots, n$. 
Given a set $I\subset \{1,2,\ldots, n\}$ we shall write
$E_I := \bigoplus_{i\in I} E_i$. 
In all examples given we shall assume that
$H_n=\Diag_n^+$.

In this section we define the stratifications associated to $H_n$ whose existence is asserted in Theorem B. 
\begin{defi}
We call {\em $(n,k)$-tree} to any sequence $\underline{P}=(P_1,P_2,\ldots, P_k)$
of subsets $P_i\subset \{1,2,\ldots, n\}$ such that
for every $1\leq i<j\leq k$ either $P_i\cap P_j=\emptyset$ or $P_i\subseteq P_j$.
In other words, $(n,k)$-trees are just trees in a set of $k$ nodes consisting of  subsets of $\{1,2,\ldots, n\}$.
We say that $k$ is the length of the $(n,k)$-tree $\underline{P}$ and write
$\Mod{\underline{P}}:=k$. The union $\cup_{j=1}^k P_j$ is called
the {\em support of} $\underline{P}$ and denoted by $\supp(\underline{P})$.
\end{defi}
We shall borrow the following terminology from Graph Theory.
We say that $P_j$ is {\em a parent of} $P_i$ \, iff\,
$P_i\subseteq P_j$. 
We say that $P_i$ is {\em the mother of  } $P_j$ \, iff\,
$j<i$, $P_j\subseteq P_i$ and there is no $j< s <i$ with $P_j\subseteq P_s\subseteq P_i$.
We say that $P_i$ is {\em a leaf } \, iff\,
there is no $1\leq s <i$ with $P_s\subseteq P_i$. In other words, $P_i$ is a node 
with no daughters.
We say that $P_i$ is {\em a root } \, iff\,
there is no $i<s\leq k$ such that $P_i\subseteq P_s$, which means that $P_i$ is a node
with no mother.
Given an increasing map $\pi:\{1,\ldots, r\}\to\{1,\ldots, k\}$
and a $(n,k)$-tree  $\underline{Q}$ we define a new $(n,r)$-tree setting 
$\pi^\ast \underline{Q}=(Q_{\pi(1)},\ldots, Q_{\pi(r)})$.
We shall refer to it a subtree  of  $\underline{Q}$.
A connected component subtree  of $\underline{Q}$ means a
$(n,r)$-tree $\underline{P}=\pi^\ast\underline{Q}$
where  $\pi:\{1,\ldots, r\}\to\{1,\ldots, k\}$ is an increasing map such that
$Q_{\pi(r)}$ is a root of $\underline{Q}$ and
$\{\pi(1),\ldots, \pi(r)\}=\{\,j\in\{1,\ldots, k\}\,:\, Q_j\subseteq Q_{\pi(r)}\, \}$.
The number of connected components of $\underline{P}$ 
is the number of roots of $\underline{P}$.
Given a tree $\underline{P}$ we define the two functions
 $\card_{\underline{P}}, \depth_{\underline{P}}:\{1,2,\ldots, k\}\to\N$ by
\begin{align*}
\card_{\underline{P}}(i) &=\# P_i\quad \text{ and } \\
\depth_{\underline{P}}(i) &=\# \{\, s\in\{1,\ldots, i-1\}\,:\, P_s\subseteq P_i\,\} \;.
\end{align*}
The second one gives the {\em depth of a node } in the tree $\underline{P}$,
which can be recursively defined as follows: leaves have depth zero. For all other nodes,
the depth is  the number of its daughter nodes plus the sum of all its daughter nodes' depths.
Alternatively,
$\depth_{\underline{P}}(i)+1$ is the length of the subtree of 
$\underline{P}$ consisting of all nodes $P_j$ such that $P_j \subseteq P_i$
($P_i$ included).

\bigskip

Given a $(n,k)$-tree $\underline{P}$ we define
\begin{align*}
\FlStr{P} &= \{\, V\in\Flag_{n,k}\,:\,
V_i\subseteq E_{P_1\cup\cdots \cup P_i},\;\forall 1\leq i\leq k\,\}\,, \\
\OOStr{P} &= \{\, X\in\OO_{n,k}\,:\,
(X)_i\in E_{P_i},\;\forall 1\leq i\leq k\,\}\,.
\end{align*}
It follows easily that  $\OOStr{P}=p^{-1}(\FlStr{P})$ for each $(n,k)$-tree
$\underline{P}$,
where $p$ is the canonical projection from $\OO_{n,k}$ onto $\Flag_{n,k}$ defined in section~\ref{covering:spaces}.
These sets will be referred as $(n,k)$-strata.

\begin{prop} \label{consistentcy}
Given a $(n,k)$-tree $\underline{P}$,
the following are equivalent:
\begin{enumerate}
\item the $(n,k)$-stratum  $\FlStr{P}$ is non-empty\;, 
\item $\depth_{\underline{P}}(i)+1\leq \card_{\underline{P}}(i)$ for every $i=1,\ldots, k$.
\end{enumerate}
\end{prop} 

\dem
Assume first $\FlStr{P}$ is non-empty and take $X \in \OOStr{P}$.
We have $\card_{\underline{P}}(i) = \dim E_{P_i}$.
Denote by $C_i(\underline{P},X)$ linear span of the columns of $(X)_j$ of $X$
with $j<i$  and $P_j\subset P_i$.
Then $\depth_{\underline{P}}(i)=\dim C_i(\underline{P},X)$,
and because $E_{P_i}$ must contain $C_i(\underline{P},X)\oplus \langle (X)_i\rangle$
item 2. follows.
The converse is proved recursively constructing (column by column) an orthogonal matrix 
$X \in \OOStr{P}$.
\cqd

This motivates the definition

\begin{defi}
We say that a $(n,k)$-tree $\underline{P}$ is {\em consistent}\, iff\,
$ \depth_{\underline{P}}(i)+1 \leq \card_{\underline{P}}(i)$ for every $i=1,\ldots, k$.
We say that a node $P_i$ {\em is full in} $\underline{P}$\, iff\,
$ \depth_{\underline{P}}(i)+1 = \card_{\underline{P}}(i)$.
\end{defi}

\bigskip

Given a $(n,k)$-tree $\underline{P}=(P_1,\ldots, P_k)$ and $1\leq i\leq k$ we define 
$\pi_i(\underline{P}):= (P_1,\ldots, P_i)$, which clearly is a  $(n,i)$-tree.
Consider the projection $\pi_i:\Flag_{n,k}\to\Flag_{n,i}$ defined in the
previous section. One can easily check that the $(n,k)$-strata are preserved by these
projections in the following sense:

\begin{prop} For any $(n,k)$-tree $\underline{P}$ and $1\leq i\leq k$ ,\,
$\pi_i(\FlStr{P})={\rm S}_{\pi_i(\underline{P})}$.
\end{prop}

\bigskip 

We are now going to prove that strata are always closed manifolds.

\begin{prop}\label{stratum:dim}
Given a consistent $(n,k)$-tree $\underline{P}$, the $(n,k)$-stratum $\FlStr{P}$ is a compact 
connected manifold (without boundary) of dimension
\begin{equation} \label{dim:P}
\dim (\FlStr{P}) = 
 -k -\kappa(\underline{P}) + \sum_{i=1}^k \# P_i   \;,
\end{equation}
where 
$\kappa(\underline{P}):=\#\{\,(i,j)\,:\, 1\leq i<j\leq k,\; P_i\subseteq P_j\,\}$.
The stratum $\OOStr{P}$ is also a compact 
manifold without boundary of the same dimension, which in general may be disconnected.
\end{prop}

\dem
Consider the linear space 
\[ \widehat{{\rm S}}_{n,k}(\underline{P})=\{\, B\in\Mat_{n\times k}(\R)\,:\,
(B)_i\in E_{P_i},\; \forall\, i=1,\ldots, k\;\}\;,\]
which has dimension $\dim  \widehat{{\rm S}}_{n,k}(\underline{P})= \sum_{i=1}^k \#(P_i)$.
Consider the set
\[ \Sigma =\{\,(i,j)\,:\, 1\leq i<j\leq k,\; P_i\subseteq P_j\,\}\;,\]
and define $\Phi: \widehat{{\rm S}}_{n,k}(\underline{P})\to \R^k\times\R^\Sigma$ by
\[ \Phi(X) = \left(\, \{\,\nrm{ (X)_i}^2 -1\,\}_{1\leq i\leq k} \, ,\;
\{\, (X)_i\cdot (X)_j\, \}_{(i,j)\in\Sigma} \; \right)
\]
Notice that $\OOStr{P} = \widehat{{\rm S}}_{n,k}(\underline{P})\cap \OO_{n,k}
=\Phi^{-1}(0,0)$.
To finish the proof it is enough showing that every point $X\in \OOStr{P}$
is a regular point of $\Phi$, in the sense that the derivative
$D\Phi_X : \widehat{{\rm S}}_{n,k}(\underline{P})\to \R^k\times\R^\Sigma$ is surjective.
Then $\OOStr{P}$ being a regular level set of $\Phi$ will be a closed manifold of the claimed
dimension.
The derivative of $\Phi$ is explicitly given by
\[ D\Phi_X(Y) = \left(\, \{\,2\,(X)_i\cdot(Y)_i\, \}_{1\leq i\leq k} \, ,\;
\{\, (X)_i\cdot (Y)_j + (Y)_i\cdot (X)_j\, \}_{(i,j)\in\Sigma} \; \right)\;.
\]
Given $(\underline{d},\underline{c})= \left(\, \{\,d_i\, \}_{1\leq i\leq k} \, ,\;
\{\, c_{i,j}\, \}_{(i,j)\in\Sigma} \; \right)\in \R^k\times\R^\Sigma$,
we shall prove by induction that there is a matrix $Y\in\Mat_{n\times k}(\R)$
such that $D\Phi_X(Y)= (\underline{d},\underline{c})$.
Fix $X\in \OOStr{P}=\Phi^{-1}(0,0)$ and 
consider the following property in the index $i$:
{\em there are vectors $Y_1,\ldots, Y_i$ in $\R^n$ such that 
\begin{enumerate}
\item $2\,(X)_j\cdot Y_j = d_j$, \, for all $j=1,\ldots, i$,
\item $(X)_s\cdot Y_j + Y_i\cdot (X)_j = d_{s,j}$, \, for all
$(s,j)\in\Sigma$ with $j\leq i$,
\item $Y_i\in E_{P_i}$,\, for all $i=1,\ldots, k$.
\end{enumerate}
}
For $i=1$ this is true: condition (2) is empty, and 
condition(1) is obvious because $(X)_1\neq 0$.
Assume this property holds for $i-1$ with vectors
$Y_1,\ldots, Y_{i-1}$. To find $Y_i$ such that (1)-(3) above hold,
we have to solve the following system of equations in the unknown $Y_i\in E_{P_i}$.
\[ \left\{
\begin{array}{ccl}
2\,(X)_i\cdot Y_i &=& d_i\\
(X)_j\cdot  Y_i &=& c'_{j,i} \quad \text{ for every }\; (j,i)\in\Sigma
\end{array}
\right. \quad,\]
where $c'_{j,i} = d_{j,i} - (X)_i\cdot Y_j$ is already determined.
This system is determined because the columns of $X$ are orthonormal.
In particular, the columns $(X)_i$ and $(X)_j$ with $(j,i)\in\Sigma$
form an orthonormal system in $E_{P_i}$.
Therefore property above holds for $i$, which completes the induction,
and proves $D\Phi_X$ is surjective.

The strata $\OOStr{P}$ are not connected in general.
For instance if $\underline{P}=(\{1\})$ then $\OOStr{P}$ consists of two points,
and if $\underline{P}=(\{1\},\{2,3\})$ then $\OOStr{P}$  consists of two circles.
The proof that every strata $\FlStr{P}$ is connected goes by induction on $k$.
For $k=1$, every stratum is diffeomorphic to some projective space and
is therefore connected.
In general the projection $\pi_{k-1}:\Flag_{n,k}\to\Flag_{n,k-1}$ preserves strata.  
For each consistent $(n,k)$-tree  $\underline{P}$
let $\underline{P}'$ be the $(n,k-1)$-tree
$\underline{P'}:=(P_1,P_2,\ldots, P_{k-1})$.
Then $\pi_{k-1}(\FlStr{P}) =\FlStr{P'}$ and the projection
$\pi_{k-1}: \FlStr{P}\to\FlStr{P'}$ is a fibration with connected fibers.
Knowing (by induction hypothesis) that $\FlStr{P'}$ is connected it follows that
$\FlStr{P}$ is also connected. Notice that each fiber of $\pi_{k-1}: \FlStr{P}\to\FlStr{P'}$
is the submanifold of all $k$-dimensional vector subspaces  $V\in \Grass_{m,k}$ 
($m=\#(\cup_{i=1}^k P_i)$) that
contain a given $k-1$ dimensional subspace, which is clearly a connected space.
\cqd

\begin{coro}
Let  $\underline{P}$ be a consistent $(n,k)$-tree  and $2\leq i\leq k$.
Then the node $P_i$ is full in $\underline{P}$ \, iff\,
$\dim {\rm S}_{\pi_i(\underline{P})} = \dim {\rm S}_{\pi_{i-1}(\underline{P})}$.
\end{coro}

\bigskip

Given two $(n,k)$-trees  $\underline{P}=(P_1,P_2,\ldots, P_k)$ 
and  $\underline{Q}=(Q_1,Q_2,\ldots, Q_k)$,
we say that $\underline{Q}$ is {\em contained in } $\underline{P}$,
and write $\underline{Q}\leq \underline{P}$, \, iff\,
$Q_i\subseteq P_i$ for every $i=1,\ldots, k$.
The relation $\leq $ is a partial order on the set of all $(n,k)$-trees. 
The following monotony is obvious.

\begin{prop}
Given $(n,k)$-trees $\underline{P}, \underline{Q}$,\, if \,
 $\underline{P}\leq \underline{Q}$\; then \;
$\FlStr{P}\subseteq \FlStr{Q}$. In particular,
$\dim \FlStr{P}\leq \dim \FlStr{Q}$.  
\end{prop}

\bigskip

Consider two trees: a $(n,r)$-tree  $\underline{P}=(P_1,\ldots, P_r)$ and
a $(n,s)$-tree   $\underline{Q}=(Q_1,\ldots, Q_s)$. We say they are disjoint
\, iff\,  $\supp(\underline{P}) \cap \supp(\underline{Q})=\emptyset$.
In this case we define their concatenation to be $(n,r+s)$-tree
$\underline{P}\diamond \underline{Q}=(P_1,\ldots, P_r,Q_1,\ldots, Q_s)$.

\begin{prop}\label{tree:product}
Given two disjoint trees $\underline{P}$ and  $\underline{Q}$,\,
$\widetilde{{\rm S}}_{\underline{P}\diamond \underline{Q}} \simeq \OOStr{P}\times\OOStr{Q}$.
\end{prop}

\dem
The diffeomorphism takes a matrix $X\in \widetilde{{\rm S}}_{\underline{P}\diamond \underline{Q}}$
to the pair of matrices $([X]_r,[X]_s)\in \OOStr{P}\times\OOStr{Q}$ 
where $[X]_r$ is formed by the first $r$ columns of $X$
and $[X]_s$ is formed by the last $s$ columns of $X$.
\cqd

\begin{coro}
For any tree $\underline{P}$, \,
$\OOStr{P}$ is diffeomorphic to the cartesian product of its connected component strata.
In particular, $\dim \FlStr{P}$ is the sum of dimensions of the connected component strata  of $\underline{P}$.
\end{coro}

\bigskip

\begin{defi}
We say that a $(n,k)$-tree $\underline{P}$ is {\em irreducible}\, iff\,
for each full node in $\underline{P}$  is a root of $\underline{P}$.
\end{defi}

Given a $(n,k)$-tree $\underline{P}$, we denote by $\underline{P}^\ast=(P_1^\ast,\ldots, P_k^\ast)$
the tree obtained from $\underline{P}$ as follows: if $\underline{P}$ is irreducible then we set
$\underline{P}^\ast:=\underline{P}$. Otherwise take the first 
$i=1,\ldots, k$ such that
$ \depth_{\underline{P}}(i)+1 = \card_{\underline{P}}(i)$
and $P_i$ is not a root of $\underline{P}$. Define $\underline{P}'=(P_1',\ldots, P_k')$
where $P_j'=P_j$ when $j\leq i$, and
$P_j'=P_j-P_i$ for all other $j$. It is not difficult to check that $\underline{P}'$
is still a $(n,k)$-tree with $\underline{P}'\leq \underline{P}$. If it is irreducible we set $\underline{P}^\ast:=\underline{P}'$. 
Otherwise we repeatedly apply this procedure until we reach an irreducible tree, 
which we set to be $\underline{P}^\ast$.
We denote by $\Kscr_{n,k}$ the set of all consistent irreducible $(n,k)$-trees.
Then $\underline{P} \mapsto \underline{P}^\ast$ is a projection operator 
mapping consistent $(n,k)$-trees onto the set $\Kscr_{n,k}$.

\begin{prop}\label{reduction:strata}
For any consistent $(n,k)$-tree $\underline{P}$,\;
$\FlStr{P}= {\rm S}_{\underline{P}^\ast}$.
\end{prop}

\dem
The inclusion ${\rm S}_{\underline{P}^\ast}\subseteq \FlStr{P}$ follows because
$\underline{P}^\ast\leq \underline{P}$. For the reverse inclusion it is enough
proving that $\FlStr{P}\subseteq {\rm S}_{\underline{P}'}$ for each one step
operation $\underline{P} \mapsto \underline{P}'$ in the reduction procedure.
Take $X\in \OOStr{P}$.
We shall use the notation introduced in the proof of proposition~\ref{consistentcy}.
Let $i=1,\ldots, k$ be the first index such that
$ \depth_{\underline{P}}(i)+1 = \card_{\underline{P}}(i)$.
Then $E_{P_i} = C_i(\underline{P},X)\oplus \langle (X)_i\rangle$ and for every  $j>i$,
because $X$ is orthogonal, we have $(X)_j\in E_{P_j-P_i}=E_{P_j'}$.
Therefore\, $X\in \widetilde{{\rm S}}_{\underline{P}'}$.
\cqd

Therefore, restricting to irreducible trees  we still get all $(n,k)$-strata .

\bigskip

\begin{lema}\label{supp:length}
Given a consistent tree $\underline{P}$,
\begin{enumerate}
\item $\# \supp(\underline{P})\geq \Mod{\underline{P}}$\,,
\item $\# \supp(\underline{P})\geq \Mod{\underline{P}}+1$\, when 
$\underline{P}$ has no full nodes in $\underline{P}$.
\end{enumerate}
\end{lema}

\dem There is an easy proof by induction in the length of $\underline{P}$.
\cqd

\begin{prop} \label{irred:dim}
If $\underline{P}, \underline{Q}\in \Kscr_{n,k}$,\,
$\underline{P}\leq \underline{Q}$\, and \, $\dim \FlStr{P}= \dim \FlStr{Q}$\,
then \; $\underline{P} = \underline{Q}$.
\end{prop}

\dem
Because strata are connected manifolds,
it is clear that
$\underline{P}\leq \underline{Q}$\; and \; $\dim \FlStr{P}= \dim \FlStr{Q}$
imply that  $\FlStr{P}=\FlStr{Q}$.
The proof goes by induction in $k$.
It is trivial for $k=1$. Assume it holds for $k-1$,
and take $\underline{P}, \underline{Q}\in \Kscr_{n,k}$,\, such that
$\underline{P}\leq \underline{Q}$\, and \, $\dim \FlStr{P}= \dim \FlStr{Q}$.
Then $\pi_{k-1}(\underline{P})\leq \pi_{k-1}(\underline{Q})$
and $\pi_{k-1}(\FlStr{P}) = \pi_{k-1}(\FlStr{P})$. Therefore
$$\dim {\rm S}_{\pi_{k-1}(\underline{P}) } = \dim \pi_{k-1}(\FlStr{P}) = 
\dim \pi_{k-1}(\FlStr{P}) = \dim {\rm S}_{\pi_{k-1}(\underline{Q}) }\;,$$
and by induction hypothesis $\pi_{k-1}(\underline{P}) = \pi_{k-1}(\underline{Q})$.
It remains to prove that $P_k=Q_k$.
Let $\underline{T}^{(i)}$ ($1\leq i \leq \ell$) be the (distinct) connected component subtrees
of $\underline{P}$ such that $\supp(\underline{T}^{(i)})\subseteq Q_k-P_k$.
Because their supports are disjoint, using lemma~\ref{supp:length}, we have
\begin{align*}
0 &= \dim \FlStr{Q} -\dim \FlStr{P} =\#(Q_k-P_k)-\sum_{j=1}^\ell \Mod{\underline{T}^{(j)}} \\
&\geq \sum_{j=1}^\ell \#\supp(\underline{T}^{(j)}) -\Mod{\underline{T}^{(j)}}\geq \ell \geq 0\;.
\end{align*}
Whence $\ell=0$, i.e.,
there is no connected component subtree of $\underline{P}$ with support contained in $Q_k-P_k$.
Therefore, $\#(Q_k-P_k)=0$ and $P_k=Q_k$.
\cqd

There is a unique maximum tree in $\Kscr_{n,k}$, which corresponds to take
$P_i=\{1,\ldots, n\}$ for every $i=1,\ldots, k$.
We have $\FlStr{P}=\Flag_{n,k}$ and $\OOStr{P}=\OO_{n,k}$ for the maximum stratum.
 The dimension of $\FlStr{P}$ is of course  $k\,(2n-k-1)/2=\dim \Flag_{n,k}$.

There are many minimal  trees in $(\Kscr_{n,k},\leq)$.
They are the trees where each $P_i$ is a singular set,\, $P_i=\{n_i\}$,
i.e., $\underline{P}=(\{n_1\},\{n_2\},\ldots, \{n_k\}) $.
The set $\{n_1,n_2,\ldots, n_k\}$ must be formed by $k$-distinct elements
from $\{1,2,\ldots, n\}$. 
The minimal trees correspond to zero dimensional  strata:
$\FlStr{P}$ is a one point set, while $\OOStr{P}$ is a set with $2^k$ points.

\bigskip

\begin{prop}
Given trees $\underline{P}, \underline{Q}\in \Kscr_{n,k}$,
then either  $\FlStr{P}\cap\FlStr{Q}=\emptyset$  or else
the infimum $P\wedge Q$ exists in $(\Kscr_{n,k},\leq )$ and\,
$\FlStr{P}\cap\FlStr{Q}={\rm S}_{ P\wedge Q}$.
\end{prop}

\dem
Given two $(n,k)$-trees $\underline{P}$ and $\underline{Q}$ define 
 $\underline{P}\cap \underline{Q}:=(P_1\cap Q_1, \ldots, P_k\cap Q_k)$.
This may not be a tree because $P_i\cap Q_i=\emptyset$ for some $i=1,\ldots, k$,
but if all these intersections are non-empty then it is straightforward  checking that
 $\underline{P}\cap \underline{Q}$ is a tree.
We claim that $\FlStr{P}\cap \FlStr{Q}={\rm S}_{\underline{P}\cap \underline{Q}}$ and
$\OOStr{P}\cap \OOStr{Q}=\widetilde{{\rm S}}_{\underline{P}\cap \underline{Q}}$, which 
follow easily from the fact \, $E_{P_i}\cap E_{Q_i} = E_{P_i\cap Q_i}$.
Now, if $\FlStr{P}\cap \FlStr{Q}$ is non-empty then the $(n,k)$-tree $\underline{P}\cap \underline{Q}$
is consistent, and the tree $\underline{P}\wedge \underline{Q} := (\underline{P}\cap \underline{Q})^\ast$
is both consistent and irreducible. This tree is the infimum of $\underline{P}$ and $\underline{Q}$
in the partially ordered set $(\Kscr_{n,k},\leq)$. By proposition~\ref{reduction:strata},
${\rm S}_{\underline{P}\wedge \underline{Q}} = {\rm S}_{\underline{P}\cap \underline{Q}}
=\FlStr{P}\cap \FlStr{Q}$.
\cqd

\bigskip

\begin{prop}
Given $\underline{P} \in \Kscr_{n,k}$
with $\dim \FlStr{P}<\dim \Flag_{n,k}$ there is $\underline{P'}\in\Kscr_{n,k}$
such that $\underline{P} \leq \underline{P}'$ and $\dim \FlStr{P'}=1+\dim \FlStr{P}$\,.
\end{prop}

\dem
The proof goes by induction in $k$.
It is obvious for $k=1$ since in this case the dimension of a stratum
corresponds to the cardinal of the unique tree node.
Assume the statement of this proposition holds for $k-1$,
and let us prove it holds for $k$ too.
Along the proof we use the following definition:
we say a tree $\underline{P}=(P_1,\ldots, P_k)$ is a full chain when
$P_1= \ldots = P_k$.
We shall consider four exhausting cases:

\noindent
{\em Case 1. } $\supp(\underline{P}) \subsetneqq \{1,2,\ldots, n\}$.
In this case we take $\alpha\in \{1,2,\ldots, n\}-\supp(\underline{P})$ and define
$\underline{P}'=(P_1',\ldots, P_k')$ setting 
$$ P_j'=\left\{
\begin{array}{lcr}
P_j & \text{ for } & 1\leq j< k\\
P_i\cup\{\alpha\} & \text{ for } & j = k
\end{array}\right.\;. $$
It is then clear that $\dim \FlStr{P'}=1+\dim \FlStr{P}$.

\noindent
{\em Case 2. } $\supp(\underline{P}) = \{1,2,\ldots, n\}$ and all the 
tree connected components of $\underline{P}$ are full chains.
Notice that in this case $\underline{P}$ must have more than one connected components,
because otherwise  $\underline{P}$ would be the maximum $(n,k)$-tree 
whose dimension equals that of $\Flag_{n,k}$.
Assume $P_i$ and $P_k$ ($i<k$) are two distinct roots of  $\underline{P}$.
Notice that $P_i\cap P_k=\emptyset$.
Choose any element $\alpha\in P_k$ and define
$\underline{P}'=(P_1',\ldots, P_k')$ setting 
$$ P_j'=\left\{
\begin{array}{lcl}
P_j & \text{ if } & 1\leq j< i\\
P_i\cup\{\alpha\} & \text{ if } & j=i\\
P_j & \text{ if } & i<j\leq k \; \text{ and }\; P_j\neq P_k\\
P_i\cup P_k & \text{ if } & i<j\leq k \; \text{ and }\; P_j=P_k
\end{array}\right.\;. $$
Let $n_i = \# P_i$ and $m=\#\{\, i<j\leq k\,:\,  P_j=P_k\,\}$.
Then
\begin{align*}
\dim {\rm S}_{\underline{P}'} -\dim \FlStr{P} &=
\sum_{j=1}^k \#( P_j')-\#( P_j) -\kappa(\underline{P}') + \kappa(\underline{P})\\
&= 1+m\,n_i -\#\{\,(j,j')\,:\, j\leq i<j',\;  P_j=P_i\; \text{ and }\; P_{j'}=P_k\,\}\\
&= 1+m\,n_i -m\,n_i = 1\;.
\end{align*}

\bigskip

The next two cases make use of the induction hypothesis.

\noindent
{\em Case 3. } $\supp(\underline{P}) = \{1,2,\ldots, n\}$ 
and the tree  $\underline{P}$ is connected.
In this case $P_k=\{1,2,\ldots, n\}$.
Take $\underline{Q}=\pi_{k-1}(\underline{P})$.
By induction hypothesis there is a $(n,k-1)$-tree $\underline{Q}'$ 
such that $\underline{Q}\leq \underline{Q}'$ and  $\dim \FlStr{Q'}=1+\dim \FlStr{Q}$.
Set $\underline{P}'=(Q_1',\ldots, Q_{k-1}', P_k)$. Then,
because $P_k'=P_k=\{1,2,\ldots, n\}$,
$$ \dim {\rm S}_{\underline{P}'} -\dim \FlStr{P}  =
\dim {\rm S}_{\underline{Q}'} -\dim \FlStr{Q} = 1\;. $$

\noindent
{\em Case 4. } $\supp(\underline{P}) = \{1,2,\ldots, n\}$,
 the tree  $\underline{P}$ is disconnected, and at least one of its
connected components is a not a full chain.
We apply the induction hypothesis to the connected component of the tree  $\underline{P}$ 
which is not a full chain. Keeping all other tree connected components of  $\underline{P}$ 
unchanged we obtain a new tree  $\underline{P}'$.
Using proposition~\ref{tree:product} we get that 
$\dim \FlStr{P'}=1+\dim \FlStr{P}$.   
\cqd

\bigskip

We define $\Sscr_{n,k}:= \{\, \FlStr{P}\,:\, \underline{P}\in\Kscr_{n,k}\,\}$ and
$\widetilde{\Sscr}_{n,k}:= \{\, \OOStr{P}\,:\, \underline{P}\in\Kscr_{n,k}\,\}$.
Collecting information above we see that

\begin{prop}
$\Sscr_{n,k}$ is a stratification on $\Fscr_{n,k}$, while 
$\widetilde{\Sscr}_{n,k}$ is a stratification on $\OO_{n,k}$
with possibly disconnected strata.
\end{prop}

\bigskip

In the rest of this section we discuss the symplectic case.
Some notation is needed: Given a subset $P\subset \{1,2,\ldots, 2n\}$, 
let $\overline{P}$ be the conjugate set
$\overline{P}:=\{\,\overline{i}=i+n\,({\rm mod}\,2n)\,:\, i\in P\,\} $,
$P^\#$ be the reduced set
$P^\#:=\{\, i\,({\rm mod}\, n)\,:\, i\in P\,\}$,
and $\widehat{P}$ be the conjugation saturated set
$\widehat{P}:=P\cup \overline{P}$. 
\begin{defi}
We call {\em symplectic $(n,k)$-tree} to any sequence $\underline{P}=(P_1,P_2,\ldots, P_k)$
of subsets $P_i\subset \{1,2,\ldots, 2n\}$ such that
for every $1\leq i<j\leq k$ 
\begin{enumerate}
\item either $P_i\cap P_j=\emptyset$ or $P_i\subseteq P_j$,
\item either $\overline{P_i}\cap P_j=\emptyset$ or $\overline{P_i}\subseteq P_j$.
\end{enumerate}
\end{defi}
The {\em full support} of a symplectic  $(n,k)$-tree $\underline{P}$ is the set
$\spsupp(\underline{P}):=\widehat{\cup_{j=1}^k P_j}$.
Given a symplectic $(n,k)$-tree $\underline{P}$   define\,
$\FlisoStr{P} := \FlStr{P}\cap \Fliso_{n,k}$
and\,
$\OspStr{P} :=  \OOStr{P}\cap \Osp_{n,k}$.
We shall refer to these as {\em symplectic } $(n,k)$-strata.
Consider the conjugate-depth function $\codepth_{\underline{P}}:\{1,2,\ldots, k\}\to\N$ defined by
$$\codepth_{\underline{P}}(i) = \# \{\,
j<i\,:\, \overline{P_j}\subseteq P_i\, \} \;.$$

\begin{prop} \label{sp:consistentcy}
Given a symplectic $(n,k)$-tree $\underline{P}$,
the following are equivalent:
\begin{enumerate}
\item the symplectic $(n,k)$-stratum  $\FlisoStr{P}$ is non-empty\;, 
\item $\depth_{\underline{P}}(i)+\codepth_{\underline{P}}(i) + 1\leq 
\card_{\underline{P}}(i)$ for every $i=1,\ldots, k$.
\end{enumerate}
\end{prop} 

\dem The proof is analogous to that of proposition~\ref{consistentcy}.
Denote by $C_i(\underline{P},X)$ the linear span of the columns of $(X)_j$ of $X$
with $j<i$  and $P_j\subset P_i$, and by
$\overline{C}_i(\underline{P},X)$ the linear span of the columns of $(X)_j$ of $X$
with $j<i$  and $\overline{P_j}\subset P_i$.
Assume $\FlisoStr{P}$ is non-empty and take $X \in \OspStr{P}$.
The key point is that $E_{P_i}$ must contain the subspace
$$C_i(\underline{P},X)\oplus \overline{C}_i(\underline{P},X)\oplus\langle (X)_i\rangle\;,$$
whose dimension is  $\depth_{\underline{P}}(i)+\codepth_{\underline{P}}(i)+1$.
\cqd

This motivates the definition

\begin{defi}
We say that a symplectic $(n,k)$-tree $\underline{P}$ is {\em symplectic consistent}\, iff\,
$ \depth_{\underline{P}}(i)+ \codepth_{\underline{P}}(i) +1 \leq \card_{\underline{P}}(i)$ for every $i=1,\ldots, k$.
When $ \depth_{\underline{P}}(i)+ \codepth_{\underline{P}}(i)+ 1 = \card_{\underline{P}}(i)$ we shall say that
$P_i$ is a {\em full isotropic node} of $\underline{P}$.
\end{defi}

We are now going to prove that symplectic strata are always closed manifolds.

\begin{prop}\label{dim:sympl:stratum}
Given a symplectic consistent $(n,k)$-tree  $\underline{P}$, the symplectic  $(n,k)$-stratum $\FlisoStr{P}$ is a compact 
connected manifold (without boundary) of dimension
\begin{equation} \label{sympl:dim:P}
\dim (\underline{P})  
= \, -k
-\kappa(\underline{P})- \mu(\underline{P}) + \sum_{i=1}^k \card_{\underline{P}}(i)  \;,
\quad \text{ where }
\end{equation}
\begin{align*}
\kappa (\underline{P})&:=\#\{\,(i,j)\,:\, 1\leq i<j\leq k,\; P_i\subseteq P_j  \,\}\,,\; \text{ and }\\
\mu (\underline{P})&:=\#\{\,(i,j)\,:\, 1\leq i<j\leq k,\; \overline{P_i}\subseteq P_j  \,\}\;.
\end{align*}
The stratum $\OspStr{P}$ is also a compact 
manifold without boundary of the same dimension, which in general may be disconnected.
\end{prop}

\dem
Consider the linear space $\widehat{{\rm S}}_{2n,k}(\underline{P})$ defined in the proof of
proposition~\ref{stratum:dim}, and the sets $\Sigma$, $\Gamma$ where
\begin{align*}
\Sigma &=\{\,(i,j)\,:\, 1\leq i<j\leq k,\; P_i\subseteq P_j \,\}\;,\\
\Gamma &=\{\,(i,j)\,:\, 1\leq i<j\leq k,\; \overline{P_i}\subseteq P_j  \,\}\;.
\end{align*}
Then define $\Phi: \widehat{{\rm S}}_{2n,k}(\underline{P})\to \R^k\times\R^{\Sigma}\times\R^{\Gamma}$ by
\[ \Phi(X) = \left(\, \{\,\nrm{ (X)_i}^2 -1\,\}_{1\leq i\leq k} \, ,\;
\{\, (X)_i\cdot (X)_j\, \}_{(i,j)\in\Sigma}\,,\;
\{\, J(X)_i\cdot (X)_j\, \}_{(i,j)\in\Gamma}  \; \right)
\]
Notice that $\OspStr{P} = \widehat{{\rm S}}_{2n,k}(\underline{P})\cap \Osp_{n,k}
=\Phi^{-1}(0,0)$.
To finish the proof it is enough showing that every point $X\in \OspStr{P}$
is a regular point of $\Phi$, in the sense that the derivative
$D\Phi_X : \widehat{{\rm S}}_{2n,k}(\underline{P})\to \R^k\times\R^\Sigma\times\R^{\Gamma}$ is surjective.
Then $\OspStr{P}$ being a regular level set of $\Phi$ will be a closed manifold of the claimed
dimension.
The derivative of $\Phi$ is explicitly given by
\begin{align*}
 D\Phi_X(Y) & = \left(\, \{\,2\,(X)_i\cdot(Y)_i\, \}_{1\leq i\leq k} \, , \right.\\
& \qquad \{\, (X)_i\cdot (Y)_j + (Y)_i\cdot (X)_j\, \}_{(i,j)\in\Sigma} \,, \\
& \qquad \left. \{\, J(X)_i\cdot (Y)_j + J(Y)_i\cdot (X)_j\, \}_{(i,j)\in\Gamma}
\; \right)\;.
\end{align*}
Given $(\underline{d},\underline{c},\underline{e})= \left(\, \{\,d_i\, \}_{1\leq i\leq k} \, ,\;
\{\, c_{i,j}\, \}_{(i,j)\in\Sigma} \,,\; \{\, e_{i,j}\, \}_{(i,j)\in\Gamma} \; \right)\in 
\R^k\times\R^\Sigma\times \R^{\Gamma}$,
we shall prove by induction that there is a matrix $Y\in\Mat_{n\times k}(\R)$
such that $D\Phi_X(Y)= (\underline{d},\underline{c},\underline{e})$.
Fix $X\in \OspStr{P}=\Phi^{-1}(0,0)$ and 
consider the following property in the index $i$:
{\em there are vectors $Y_1,\ldots, Y_i$ in $\R^n$ such that 
\begin{enumerate}
\item $2\,(X)_j\cdot Y_j = d_j$, \, for all $j=1,\ldots, i$,
\item $(X)_s\cdot Y_j + Y_i\cdot (X)_j = d_{s,j}$, \, for all
$(s,j)\in\Sigma$ with $j\leq i$,
\item $J (X)_s\cdot Y_j + J Y_i\cdot (X)_j = e_{s,j}$, \, for all
$(s,j)\in\Gamma$ with $j\leq i$,
\item $Y_i\in E_{P_i}$,\, for all $i=1,\ldots, k$.
\end{enumerate}
}
For $i=1$ this is true: condition (2) and (3) are empty, and 
condition(1) is obvious because $(X)_1\neq 0$.
Assume this property holds for $i-1$ with vectors
$Y_1,\ldots, Y_{i-1}$. To find $Y_i$ such that (1) and (2) above hold,
we have to solve the following system of equations in the unknown $Y_i\in E_{P_i}$.
\[ \left\{
\begin{array}{ccl}
2\,(X)_i\cdot Y_i &=& d_i\\
(X)_j\cdot  Y_i &=& c'_{j,i} \quad \text{ for every }\; (j,i)\in\Sigma\\
J(X)_j\cdot  Y_i &=& e'_{j,i} \quad \text{ for every }\; (j,i)\in\Gamma\\
\end{array}
\right. \quad,\]
where $c'_{j,i} = d_{j,i} - (X)_i\cdot Y_j$ and
 $e'_{j,i} = e_{j,i} + J(X)_i\cdot Y_j$ are already determined.
This system has a solution because the columns of $X$ and $J X$ are orthogonal.
In particular, the columns 
$(X)_i$, $(X)_j$ with $(j,i)\in\Sigma$ and $J(X)_j$ with $(j,i)\in\Gamma$
form an orthonormal system in $E_{P_i}$.
Therefore property above holds for $i$, which completes the induction
and proves that $D \Phi_X$ is surjective.

The connectedness of strata $\FlisoStr{P}$ 
is proved as in proposition~\ref{stratum:dim}.
\cqd

\begin{coro}
Let  $\underline{P}$ be a symplectic consistent $(n,k)$-tree  and $2\leq i\leq k$.
Then the node $P_i$ is full isotropic in $\underline{P}$ \, iff\,
$\dim {\rm S}^{sp}_{\pi_i(\underline{P})} = \dim {\rm S}^{sp}_{\pi_{i-1}(\underline{P})}$.
\end{coro}

\begin{defi}
We say that a symplectic consistent  $(n,k)$-tree $\underline{P}$ is {\em symplectically irreducible}\, iff\,
each full isotropic node $P_i$ is a root of $\underline{P}$. 
\end{defi}

\bigskip

The partial order among symplectic $(n,k)$-trees is the same
 $\underline{Q}\leq \underline{P}$  \, iff\, $Q_i\subseteq P_i$ for every $i=1,\ldots, k$.
Then the following monotony is obvious.

\begin{prop}
Given symplectic $(n,k)$-trees $\underline{P}, \underline{Q}$,\, if \,
 $\underline{P}\leq \underline{Q}$\; then \;
$\FlisoStr{P}\subseteq \FlisoStr{Q}$. In particular,
$\dim \FlisoStr{P}\leq \dim \FlisoStr{Q}$.  
\end{prop}

\bigskip

Concatenation of symplectic trees $\underline{P}\diamond \underline{Q}$
is defined when their full supports are disjoint, i.e.,
 $\spsupp(\underline{P}) \cap \spsupp(\underline{Q})=\emptyset$.

\begin{prop}\label{sp:tree:product}
Given two symplectic trees $\underline{P}$ and  $\underline{Q}$ with disjoint full supports,\,
$\widetilde{{\rm S}}^{sp}_{\underline{P}\diamond \underline{Q}} \simeq \OspStr{P}\times\OspStr{Q}$.
\end{prop}

\begin{coro}
For any symplectic tree $\underline{P}$ whose roots satisfy
$\widehat{P_i}\cap \widehat{P_j}=\emptyset$ for every pair of distinct roots $P_i,P_j$,\,
$\OspStr{P}$ is diffeomorphic to the cartesian product of its connected component strata.
In particular, $\dim \FlisoStr{P}$ is the sum of dimensions of the connected component strata  of $\underline{P}$.
\end{coro}

\bigskip

Given a symplectic $(n,k)$-tree $\underline{P}$, we denote by $\underline{P}^\ast=(P_1^\ast,\ldots, P_k^\ast)$
the tree obtained from $\underline{P}$ as follows: if $\underline{P}$ is symplectically irreducible then we set
$\underline{P}^\ast:=\underline{P}$. Otherwise take the first 
$i=1,\ldots, k$ such that
$ \depth_{\underline{P}}(i)+\codepth_{\underline{P}}(i)+ 1 = \card_{\underline{P}}(i)$
and $P_i$ is not a root of $\underline{P}$.
Define $\underline{P}'=(P_1',\ldots, P_k')$ setting
$\underline{P}'=(P_1',\ldots, P_k')$ setting 
$$ P_j'=\left\{
\begin{array}{ll}
P_j & \text{ if } \; 1\leq j\leq i\; \\ 
P_j-\widehat{P_i}  & \text{ otherwise }
\end{array}\right.\;. $$
It is not difficult to check that $\underline{P}'$
is still a symplectic $(n,k)$-tree with $\underline{P}'\leq \underline{P}$. If it is symplectically irreducible we set $\underline{P}^\ast:=\underline{P}'$. 
Otherwise we repeatedly apply this procedure until we reach a symplectically irreducible tree, 
which we set to be $\underline{P}^\ast$.
We denote by $\Kscr^{sp}_{n,k}$ the set of all symplectic consistent and symplectically irreducible $(n,k)$-trees.
Then $\underline{P} \mapsto \underline{P}^\ast$ is a projection operator 
mapping symplectic consistent $(n,k)$-trees onto the set $\Kscr^{sp}_{n,k}$.

\begin{prop}\label{reduction:sympl:strata}
For any symplectic consistent $(n,k)$-tree $\underline{P}$,\;
$\FlisoStr{P}= {\rm S}^{sp}_{\underline{P}^\ast}$.
\end{prop}

\dem
The inclusion ${\rm S}^{sp}_{\underline{P}^\ast}\subseteq \FlisoStr{P}$ follows because
$\underline{P}^\ast\leq \underline{P}$. For the reverse inclusion it is enough
proving that $\FlisoStr{P}\subseteq {\rm S}_{\underline{P}'}$ for each one step
operation $\underline{P} \mapsto \underline{P}'$ in the symplectic reduction procedure.
Take $X\in \OspStr{P}$.
We shall use the notation introduced in the proof of proposition~\ref{sp:consistentcy}.
Let $i=1,\ldots, k$ be the first index such that
$ \depth_{\underline{P}}(i)+\codepth_{\underline{P}}(i)+1 = \card_{\underline{P}}(i)$.
Then $C_i(\underline{P},X)\oplus \overline{C}_i(\underline{P},X)\oplus\langle (X)_i\rangle = E_{P_i}$ 
and for every  $j$ such that 
$j> i$, 
because $X$ is unitary, we have $(X)_j\in E_{P_j-\widehat{P_i}}=E_{P_j'}$.
Therefore\, $X\in \widetilde{{\rm S}}^{sp}_{\underline{P}'}$.
\cqd

This means that restricting to symplectic irreducible trees  we still get all symplectic $(n,k)$-strata .

\bigskip

\begin{prop} \label{sympl:irred:dim}
If $\underline{P}, \underline{Q}\in \Kscr^{sp}_{n,k}$,\,
$\underline{P}\leq \underline{Q}$\, and \, $\dim \FlStr{P}= \dim \FlStr{Q}$\,
then \; $\underline{P} = \underline{Q}$.
\end{prop}

Because this section is already long we omit the proof  of this  and of the following
propositions, which are simple adaptations of the general case arguments.

\bigskip

There is a unique maximum symplectic tree in $\Kscr^{sp}_{n,k}$, which corresponds to take
$P_i=\{1,\ldots, 2n\}$ for every $i=1,\ldots, k$.
We have $\FlisoStr{P}=\Fliso_{n,k}$ and $\OspStr{P}=\Osp_{n,k}$ for the maximum stratum.
 The dimension of $\FlisoStr{P}$ is $k\,(2n-k)$.

There are many minimal  trees in $(\Kscr^{sp}_{n,k},\leq)$.
They are the symplectic trees where each $P_i$ is a singular set,\, $P_i=\{n_i\}$,
i.e., $\underline{P}=(\{n_1\},\{n_2\},\ldots, \{n_k\}) $.
The set $\{n_1, \overline{n_1}, n_2, \overline{n_2}, \ldots, n_k,  \overline{n_k}\}$  
must be formed by $2k$-distinct elements
from $\{1,2,\ldots, 2n\}$. 
The minimal symplectic trees correspond to zero dimensional  strata:
$\FlisoStr{P}$ is a one point set, while $\OspStr{P}$ is a set with $2^k$ points.

\bigskip

\begin{prop}
Given trees $\underline{P}, \underline{Q}\in \Kscr^{sp}_{n,k}$,
then either  $\FlisoStr{P}\cap\FlisoStr{Q}=\emptyset$  or else
the infimum $P\wedge Q$ exists in $(\Kscr^{sp}_{n,k},\leq )$ and\,
$\FlisoStr{P}\cap\FlisoStr{Q}={\rm S}^{sp}_{ P\wedge Q}$.
\end{prop}

\bigskip

\begin{prop}
Given $\underline{P} \in \Kscr^{sp}_{n,k}$
with $\dim \FlisoStr{P}<\dim \Fliso_{n,k}$ there is $\underline{P}'\in\Kscr^{sp}_{n,k}$
such that $\underline{P} \leq \underline{P}'$ and $\dim {\rm S}^{sp}_{\underline{P}'}=1+\dim \FlisoStr{P}$\,.
\end{prop}

\bigskip

We define $\widetilde{\Sscr}_{n,k}:= \{\, \FlisoStr{P}\,:\, \underline{P}\in\Kscr^{sp}_{n,k}\,\}$ and
$\widetilde{\Sscr}_{n,k}:= \{\, \OspStr{P}\,:\, \underline{P}\in\Kscr^{sp}_{n,k}\,\}$.
Collecting information above we see that

\begin{prop}
$\Sscr^{sp}_{n,k}$ is a stratification on $\Fliso_{n,k}$, while 
$\widetilde{\Sscr}^{sp}_{n,k}$ is a stratification on $\Osp_{n,k}$
with possibly disconnected strata.
\end{prop}

\bigskip
\bigskip

\section{Invariance of the Stratifications }

The aim of this section is to establish the invariance of the previous section' stratifications
under both the diffeomorphisms $\varphi_A$ and the gradient
flows of the functions  $Q_{A,b}$.

\begin{prop}
 The stratification $\Sscr_{n,k}$ is  invariant under $\varphi_A$ for any
given $A\in H_n$.
\end{prop}

\dem
Given $P\in\Kscr_{n,k}$ it is enough seeing that
$\OOStr{P}=\varphi_A(\OOStr{P})$ in order to prove that $\FlStr{P}$
is invariant under $\varphi_A$. Take $X\in \OOStr{P}$.
We need to prove that $(A\ast X)_i\in E_{P_i}$, for every $i=1,\ldots, k$,
which is done by induction. 
We shall write $x\doteq y$ to express colinearity of the vectors $x$ and $y$.
For $i=1$ we have
$(A\ast X)_1\doteq  A\,(X)_1 \in E_{P_1}$ because $A\,E_{P_1}=E_{P_1}$.
Assume now that $(A\ast X)_j \in E_{P_j}$ for every $j<i$, and let
$I=\{\, j<i\,:\,P_j\subseteq P_i\,\}$.
Using this assumption we get
\begin{align*}
(A\ast X)_i &\doteq (A\,X)_i -
\sum_{j<i} \frac{(A\,X)_i\cdot (A\ast X)_j}{\nrm{ (A\ast X)_j }^2}\,(A\ast X)_j\\
&=  A\,(X)_i -
\sum_{j\in I} \frac{A\,(X)_i\cdot (A\ast X)_j}{\nrm{ (A\ast X)_j }^2}\,
\underbrace{ (A\ast X)_j }_{\in E_{P_j}\subseteq E_{P_i} }\;,
\end{align*}
which proves that $(A\ast X)_i\in E_{P_i}$ because $A\,E_{P_i}=E_{P_i}$.
The invariance in the symplectic case follows from the general invariance, because
the submanifold $\Fliso_{n,k}$ of $\Flag_{2n,k}$ is invariant under the diffeomorphism 
$\varphi_A:\Flag_{2n,k}\to\Flag_{2n,k}$ 
 and by definition of symplectic strata
$\FlisoStr{P}=\FlStr{P}\cap \Fliso_{n,k}$.
\cqd

\bigskip

We consider the manifolds  $\OO_{n,k}$ and $\Osp_{n,k}$ with Riemannian metrics
induced from the Hilbert-Schmidt inner product on the respective matrix spaces,
$\Mat_{n\times k}(\R)$ and $\Mat_{2n\times k}(\R)$.
Then, using the covering mappings $p:\OO_{n,k}\to\Flag_{n,k}$
and $p:\Osp_{n,k}\to\Fliso_{n,k}$, we project these Riemannian structures
onto the flag manifolds $\Flag_{n,k}$ and $\Fliso_{n,k}$.

\begin{prop}
 The stratification $\Sscr_{n,k}$ is  invariant under  
the gradient flow of the function $Q_{A,b}:\Flag_{n,k}\to\R$
for any $A\in H_n$ and $b\in\R^k$.
\end{prop}

\dem
We have only to prove that that the gradient of $Q_{A,b}$ is tangent to $\Flag_{n,k}$,
and we shall do it on the covering space $\OO_{n,k}$.
It is enough to prove this proposition for the function $Q_{A}:\OO_{n,k}\to\R$,
because of~(\ref{Qab:decomp}) and the fact that the derivatives of the the projections $\pi_i$
are self-adjoint operators (see the proof of proposition~\ref{grad:QAb:decomp}).
The matrix $P^X=I-X\,X^t\in\Mat_{n\times n}(\R)$ represents the orthogonal
projection onto the orthogonal complement $\langle X\rangle^\perp$ of the linear space
spanned by the columns of $X$.
Consider a $(n,k)$-stratum $\underline{P}$  and let
 $\widehat{{\rm S}}_{n,k}(\underline{P})$  denote the linear space
defined in proposition~\ref{stratum:dim}.
Then  $\OOStr{P}= \widehat{{\rm S}}_{n,k}(\underline{P})\cap \OO_{n,k}$
and $T_X \OOStr{P} = \widehat{{\rm S}}_{n,k}(\underline{P}) \cap T_X \OO_{n,k}$.
By lemma~\ref{QA:grads}, where we have computed the gradient of $Q_A$, we only need to prove that
\[X\in \OOStr{P}  \quad \Rightarrow \quad P^X\,A^2\,X\in \widehat{{\rm S}}_{n,k}(\underline{P}) \;.\]
Take $X\in \OOStr{P}$  so that $(X)_i\in E_{P_i}$ for every $i=1,\ldots, k$.
Because the linear spaces $E_i$ are eigen-directions of the symmetric matrix $A$
it follows that $A^2 X\in \widehat{{\rm S}}_{n,k}(\underline{P})$.
Whence it is enough to show that
$P^X\,E_{P_i}\subset E_{P_i}$, for every $i=1,\ldots, k$.
Let us say that a basis $(w_1,\ldots, w_n)$ of $\R^n$ is adapted to $\underline{P}$
when each space $E_{P_i}$ is spanned by the vectors $w_j$ ($1\leq j\leq n$)  with $w_j\in E_{P_i}$.
By induction in $k$ we can easily prove that the columns of $X\in \OOStr{P}$ can always be
extended to form a basis  $(w_1,\ldots, w_n)$  
adapted to $\underline{P}$.
Therefore for every $j=1,\ldots, n$, either $P\,w_j=0$ (when $w_j$ is a column of $X$),
or else $P\,w_j=w_j$, which implies that $P^X\,E_{P_i}\subset E_{P_i}$.

The invariance in the symplectic case follows from the general invariance, because
the submanifold $\Fliso_{n,k}$ of $\Flag_{2n,k}$ is invariant under the gradient flow of 
$Q_{A,b}:\Flag_{2n,k}\to\R$ 
 and by definition of symplectic strata
$\FlisoStr{P}=\FlStr{P}\cap \Fliso_{n,k}$.

\cqd

\bigskip
\bigskip

\section{The Stratifications Skeleton Graphs }

In this section we describe the structure of zero and one dimensional strata.

\bigskip

Given $1\leq k\leq n$ we denote by  $\Sperm_{n,k}$ the set of all
permutations $\pi=(\pi_1,\ldots, \pi_k)$ with length $k$ of the set $\{1,2,\ldots, n\}$.
$\Sperm_{n,k}$ is a set with $n!/(n-k)!$ elements.
When $k=n$ we write $\Sperm_n$ instead of $\Sperm_{n,n}$.
This set has a group structure: it is the so called $n^{\text{th}}$ symmetric group.
Given a permutation $\pi\in \Sperm_{n,k}$ we denote by $\underline{\pi}$
the $(n,k)$-stratum $\underline{\pi}:=(\{\pi_1\},\{\pi_2\},\ldots, \{\pi_k\} )$,
we denote by $V_\pi$ the flag in $\Flag_{n,k}$ defined by
$$V_\pi =\left(\, E_{\pi_1},\, E_{\pi_1}\oplus E_{\pi_2},\,\ldots,\,
E_{\pi_1}\oplus \cdots \oplus E_{\pi_k}\,\right)\;, $$
and we denote by $X_\pi$ any matrix in $\OO_{n,k}$ whose $i^{\text{th}}$-column
belongs to $E_{\pi_i}$ ($1\leq i\leq k$).
There are $2^k$ such matrices, all them in the coset $X_\pi\,\Pm_k$. We have
$p(X_\pi)=V_\pi \in \FlStr{\pi}$.

\begin{prop}
A  stratum $\underline{P}\in \Kscr_{n,k}$ has dimension zero \, iff\,
there is a permutation $\pi\in\Sperm_{n,k}$ such that
$\underline{P}=\underline{\pi}$.
\end{prop}

This shows that the correspondence $\pi\mapsto \underline{\pi}$ is one-to-one from
$\Sperm_{n,k}$ onto the set of zero dimensional strata in $\Kscr_{n,k}$.

For the symplectic case we make the following definition:
A permutation $\pi\in \Sperm_{2n,k}$ is called {\em isotropic}\, iff\, the set
$\{\,\pi_i\,({\rm mod}\,n)\,:\, i\in\{1,\ldots, k\}\,\}$ has $k$ distinct elements in $\Z_n$.
We denote by  $\SJ_{n,k}$ the set of all
$(2n)(2n-2) \cdots (2n-2k+2)$ isotropic permutations of length $k$.
We say that  $\pi\in \Sperm_{2n}$ is an {\em isotropic permutation}\, iff\,
$\pi_{\overline{i}}=\overline{\pi_i}$ for every every $i=1,2,\ldots, 2n$,
where $\overline{i}=i+n({\rm mod}\,2n)$. We denote by $\SJ_n$
the subgroup of all isotropic permutations in $\Sperm_{2n}$. This group has
order $2^n n!$. When $k=n$ the set $\SJ_{n,n}$ can be identified with the group
$\SJ_n$ because each permutation of length $n$ in $\SJ_{n,n}$ can be
uniquely extended to an isotropic permutation in $\SJ_{n}$.
Given an isotropic permutation $\pi\in \SJ_{n,k}$, the flag $V_\pi$ is isotropic,
i.e., $V_\pi\in \Fliso_{n,k}$, the matrix $X_\pi\in \OO_{2n,k}$ is unitary,
i.e., $X_\pi\in \Osp_{n,k}$, and the same relation 
$p(X_\pi)=V_\pi \in \FlisoStr{\pi}$ holds.

\begin{prop}
A  symplectic stratum $\underline{P}\in \Kscr^{sp}_{n,k}$ has dimension zero \, iff\,
there is an isotropic  permutation $\pi\in\SJ_{n,k}$ such that
$\underline{P}=\underline{\pi}$.
\end{prop}

Therefore  $\pi\mapsto \underline{\pi}$ is a one-to-one correspondence from
$\SJ_{n,k}$ onto the set of zero dimensional strata in $\Kscr^{sp}_{n,k}$.

\bigskip

We now turn our attention to one-dimensional strata.

\begin{prop} \label{one:dim:stratum:charac}
A  stratum $\underline{P}\in \Kscr_{n,k}$ has dimension one \, iff\,
either 
\begin{enumerate}
\item $\supp(\underline{P})$ has $k+1$ elements, all but one node has
exactly one element, while the exceptional node has two elements, 
or else 
\item  $\supp(\underline{P})$ has $k$ elements, all but two nodes have
exactly one element, while the two exceptional nodes share the same two elements.
\end{enumerate}
In any case $\FlStr{P}$ is diffeomorphic to a circle and contains exactly two
zero dimensional strata.
\end{prop}

\dem
Assume $\dim \FlStr{P}=1$ and look at the dimensions
$\dim {\rm S}_{\pi_i(\underline{P}) }$ ($1\leq i \leq k$).
We see at once that there is exactly one index $i$ for which
$P_i$ is not a full node of $\underline{P}$.
We have $\dim {\rm S}_{\pi_j(\underline{P}) }=1$ for all $j\geq i$  and 
$\dim {\rm S}_{\pi_{j}(\underline{P}) }=0$ for all $j<i$.
Recall that full nodes are roots because we have assumed that
$\underline{P}\in \Kscr_{n,k}$.
For $j<i$ the node $P_j$ is root corresponding to a zero dimensional
connected component subtree, whence it has exactly one element.
The node $P_i$ must have exactly two distinct elements,
different from all elements in the previous nodes, because $\dim {\rm S}_{\pi_i(\underline{P}) }=1$.
We have two possibilities:
either for some $j>i$,  $P_{j}$ is a root with the same two elements in $P_i$, and we are in case 1., or else
every $P_{j}$ with $j>i$ is a root with exactly one element (distinct fro the previous ones), and we are in case 2..

The stratum $\FlStr{P}$ is diffeomorphic to ${\rm S}_{\pi_i(\underline{P}) }$
through the projection $\pi_i$, but this is clearly diffeomorphic to the $(2,1)$-stratum of the tree
$\underline{Q}=(\{1,2\})$, and a simple computation shows $\FlStr{Q}$ is the projective line,
therefore a circle. The fact that there is only room for two zero dimensional strata
inside $\FlStr{P}$ is obvious.
\cqd

\bigskip

We say that two distinct permutations $\pi,\pi'\in \Sperm_{n,k}$ {\em are linked}\, iff\,
there is a one-dimensional stratum $\underline{P}\in\Kscr_{n,k}$ such that
 $\underline{P}$ contains  $\underline{\pi}$ and $\underline{\pi'}$.
With this linking relation $\Sperm_{n,k}$ becomes a graph that we refer as the
{\em skeleton graph } of the stratification  $\Sscr_{n,k}$.
We shall use the notation  ${\rm R}(\pi)=\{\pi_1,\ldots, \pi_k\}$ and 
${\rm C}(\pi)=\{1,\ldots, n\}-{\rm R}(\pi)$.

\begin{prop}\label{Snk:linking}
Given two permutations $\pi,\pi'\in\Sperm_{n,k}$, 
$\pi$ and $\pi'$  are linked \, iff\,  either
\begin{enumerate}
\item[(1)] there is some $(i,j)\in \{1,\ldots, k\}\times {\rm C}(\pi)$
such that $\pi'$ is obtained from $\pi$ replacing $\pi_i$ by $j$, or else
\item[(2)] there is a pair $(i,j)$ with $1\leq i<j\leq k$
such that $\pi'$ is obtained from $\pi$ switching the entries $\pi_i$ and $\pi_j$ of $\pi$.
\end{enumerate}
\end{prop}

\dem
Assume $\pi,\pi\in\Sperm_{n,k}$ are linked, and let $\underline{P}\in\Kscr_{n,k}$
be the one-dimensional tree whose stratum contains both $\underline{\pi}$ and $\underline{\pi'}$.
Case 1. of proposition~\ref{one:dim:stratum:charac} implies case (1) here,
while case 2. of proposition~\ref{one:dim:stratum:charac} implies case (2) here.
The converse is also clear.
\cqd

Next we prove that the skeleton graph  of the stratification  $\Sscr_{n,k}$
is a regular graph whose degree is the dimension of the flag manifold $\Flag_{n,k}$.

\begin{prop}
$\Sperm_{n,k}$ is a regular graph of degree $d=k\,(2n-k-1)/2$.
\end{prop}

\dem
By proposition~\ref{Snk:linking}, the degree of graph $\Sperm_{n,k}$ at $\pi$ is the number of
pairs $(i,j)$ in case (2), which is $k\,(k-1)/2$, plus the number
 of pairs $(i,j)$ in case (1), which is $k\,(n-k)$. 
Therefore the degree of graph $\Sperm_{n,k}$ at $\pi$
is $d= k\,(k-1)/2 + k\,(n-k) = k\,(2n-k-1)/2$.
\cqd

\bigskip

With a straightforward adaptation of the argument used in the proof of
proposition~\ref{one:dim:stratum:charac} we can show that:

\begin{prop}\label{one:dim:sympl:stratum:charac}
A  symplectic stratum $\underline{P}\in \Kscr^{sp}_{n,k}$ has dimension one \, iff\,
either 
\begin{enumerate}
\item $\supp(\underline{P})$ has $k+1$ elements, all but one node $P_i$ has
exactly one element and is a full isotropic node in $\underline{P}$, 
while the exceptional node $P_i$ has two elements $P_i=\{a,b\}$, 
or else 
\item  $\supp(\underline{P})$ has $k$ elements, all but two nodes $P_i,\, P_j$ ($i<j$) have
exactly one element and are full isotropic nodes in $\underline{P}$, 
while the two exceptional nodes $P_i,\, P_j$ share the same two elements 
$P_i= P_j$ and $P_j$ is
a full isotropic node  in $\underline{P}$, or else
\item  $\supp(\underline{P})$ has $k+2$ elements, all but two nodes $P_i,\, P_j$ ($i<j$) have
exactly one element and are full isotropic nodes in $\underline{P}$, 
while the two exceptional nodes $P_i,\, P_j$ are such that $P_i\cap P_j=\emptyset$,
$\overline{P_i}= P_j$ and $P_j$ is
a full isotropic node  in $\underline{P}$.
\end{enumerate}
In any case $\FlStr{P}$ is diffeomorphic to a circle and contains exactly two
zero dimensional symplectic strata.
\end{prop}

\bigskip 

We say that two distinct isotropic permutations $\pi,\pi'\in \SJ_{n,k}$ {\em are linked}\, iff\,
there is a one-dimensional stratum $\underline{P}\in\Kscr^{sp}_{n,k}$ such that
 $\underline{P}$ contains  $\underline{\pi}$ and $\underline{\pi'}$.
With this linking relation $\SJ_{n,k}$ becomes a graph that we refer as the
{\em skeleton graph } of the stratification  $\Sscr^{sp}_{n,k}$.
As before we write  ${\rm R}(\pi)=\{\pi_1,\ldots, \pi_k\}$
and \,  $\overline{i}=i+n\,({\rm mod}\,2n)$.

\begin{prop}\label{SJnk:linking}
Given two permutations $\pi,\pi'\in \SJ_{n,k}$, $\pi$ and $\pi'$  are linked \, iff\,  
one of the following
cases occurs:
\begin{enumerate}
\item[(1)] there is $1\leq i\leq k$ such that
$\pi'$ is obtained from $\pi$ replacing
$\pi_i$ by $\overline{\pi_i}$;

\item[(2)]  there are $1\leq i\leq k$ and $j\notin {\rm R}(\pi)\cup \overline{ {\rm R}(\pi) }$
such that $\pi'$ is obtained from $\pi$ replacing
$\pi_i$ by $j$.

\item[(3)]  there are $1\leq i<j \leq k$
such that $\pi'$ is obtained from $\pi$ switching
$\pi_i$ with $\pi_j$;

\item[(4)]  there are $1\leq i<j \leq k$
such that $\pi'$ is obtained from $\pi$ replacing
$\pi_i$ by $\overline{\pi_j}$ and $\pi_j$ by  $\overline{\pi_i}$;

\end{enumerate}
\end{prop}
\dem
Assume $\pi,\pi\in\SJ_{n,k}$ are linked, and let $\underline{P}\in\Kscr^{sp}_{n,k}$
be the one-dimensional tree whose symplectic stratum contains both $\underline{\pi}$ and $\underline{\pi'}$.
Case 1. of proposition~\ref{one:dim:sympl:stratum:charac} 
with $P_i=\{a,b\}$ and $b=\overline{a}$  implies case (1) here;
case 1. of proposition~\ref{one:dim:sympl:stratum:charac} 
with $P_i=\{a,b\}$ and $b\neq \overline{a}$  implies case (2) here;
case 2. of proposition~\ref{one:dim:sympl:stratum:charac} implies case (3) here,
and finally case 3. of proposition~\ref{one:dim:sympl:stratum:charac} implies case (4) here.
\cqd

\bigskip

\begin{prop}
$\SJ_{n,k}$ is a regular graph of degree $d=k\,(2n-k)$.
\end{prop}

\dem
By proposition~\ref{SJnk:linking}, the degree of graph $\SJ_{n,k}$ at $\pi$ is the 
sum of the numbers of pairs $(i,j)$ in each of the four cases (1)-(4) above.
These numbers are: $k\,(k-1)/2$ in case (3),
$k$ in case (1),
$k\,(k-1)/2$ again in case (4)
and $k\,(2n-2k)$ in case (2).
Therefore the degree of graph $\SJ_{n,k}$ at $\pi$
is $d= k\,(k-1)/2 + k +  k\,(n-k) = k\,(2n-2k) = 2nk-k^2$.
\cqd

\bigskip
\bigskip

\section{Oriented Skeleton Graphs}

In the sequel we use the order  of the eigendirections $\Escr=\{E_i\}_i$ of the group $H_n$
fixed in the beginning of section~\ref{the:stratifications}.
Given a positive definite matrix $A\in H_n$, we say that $\Escr$ is $A$-ordered\, iff\,
denoting by $\lambda_i(A)$ the eigenvalue of $A$ associated with the eigenspace $E_i$
we have\, $\lambda_m(A)\geq \lambda_{m-1}(A)\geq \ldots\, \geq \lambda_2(A)\geq \lambda_1(A)$,
where $m=n$ or $m=2n$ (in the symplectic case).
We denote by $H_n(\Escr)$ the subset of all matrices $A\in H_n$ such that
 $\Escr$ is $A$-ordered. In all statements of this and the following section 
where a matrix $A\in H_n$ plays a role we shall assume that $A\in H_n(\Escr)$.

\bigskip

Next we are going to introduce orientations on $\Sperm_{n,k}$ and $\SJ_{n,k}$
which will make them oriented graphs with no cycles.

\begin{defi}\label{leads:to}
Given two permutations $\pi,\pi'\in\Sperm_{n,k}$, we say that 
$\pi$ {\em leads to} $\pi'$, and write $\pi \leadsto \pi'$, \, iff\,  either
\begin{enumerate}
\item[(1)] there is $(i,j)\in \{1,\ldots, k\}\times {\rm C}(\pi)$ with $\pi_i<j$
such that $\pi'$ is obtained from $\pi$ replacing $\pi_i$ by $j$, or else
\item[(2)] there is $(i,j)$ with $1\leq i<j\leq k$,
$\pi_i<\pi_j$  and $\pi'$ is obtained from $\pi$  switching the entries $\pi_i$ and $\pi_j$.
\end{enumerate}
\end{defi}
\begin{defi}\label{preceeds:to}
Given two permutations $\pi,\pi'\in\Sperm_{n,k}$, we say that 
$\pi$ {\em precedes } $\pi'$  \, iff\, either
\begin{enumerate}
\item[(1)] there is $(i,j)\in \{1,\ldots, k\}\times {\rm C}(\pi)$ with $j=\pi_i+1$,
and $\pi'$ obtained from $\pi$ replacing $\pi_i$ by $j$,\, or else
\item[(2)] there is a pair $(i,j)$ with $1\leq i<j\leq n$, $\pi_i<\pi_j$,
 $\pi'$  obtained from $\pi$  switching the entries $\pi_i$ and $\pi_j$, for which
 no $s\in \{1,\ldots, k\}$ exists with $\pi_i<\pi_s<\pi_j$.
\end{enumerate}
\end{defi}

Comparing the definition of the  "leads-to" relation with proposition~\ref{Snk:linking}
we see that given $\pi,\pi'\in\Sperm_{n,k}$,
$\pi$ and $\pi'$ are linked\, iff\, either $\pi$ leads to $\pi'$ or else
 $\pi'$ leads to $\pi$.
In other words,   the "leads-to" relation orients each and every edge of $\Sperm_{n,k}$,
making it an oriented graph, with an oriented edge from $\pi$ to $\pi'$ whenever $\pi \leadsto\pi'$.
The "precede" relation determines a subgraph of $\Sperm_{n,k}$ with an 
an oriented edge from $\pi$ to $\pi'$ whenever $\pi$ precedes $\pi'$.
\begin{lema}
Given $\pi, \pi'\in\Sperm_{n,k}$, if $\pi$ leads to $\pi'$ then 
there is a sequence of permutations
 $\pi^{(0)}, \pi^{(1)}, \ldots, \pi^{(m)}$ in $\Sperm_{n,k}$
such that
\begin{enumerate}
\item $\pi=\pi^{(0)}$ and $\pi'=\pi^{(m)}$,
\item $\pi^{(i)}$ precedes $\pi^{(i+1)}$, for $i=0,1,\ldots, m-1$.
\end{enumerate}
\end{lema}

This lemma proves that the transitive closures of  relations
"leads-to" and "precede" do coincide. We shall denote this
transitive closure by $\prec$,
and write $\pi \preceq \pi'$ \, iff\, $\pi=\pi'$ or $\pi\prec\pi'$.

\bigskip

Given a  permutation $\pi\in\Sperm_{n,k}$, 
consider the set $\Sigma_{n,k}(\pi)$ of all pairs $(i,j)$
satisfying either one of the following conditions:
\begin{enumerate}
\item[(1)] $1\leq i  \leq k$, $j\notin {\rm C}(\pi)$  and \, $\pi_i > j$;
\item[(2)] $1\leq i<j \leq k$ and \, $\pi_i > \pi_j$\;,
\end{enumerate}
and define the function  $H:\Sperm_{n,k}\to\N$ by\,
$ H(\pi):=\# \Sigma_{n,k}(\pi)$. 
Note that $H(\pi)$ counts the number of oriented edges of
$\Sperm_{n,k}$ which arrive at the vertex $\pi$.
Later we shall see (c.f. propositions~\ref{varphiA:Jacobian} and~\ref{QA:Hessian}) that
$H(\pi)$ is exactly the dimension of the stable manifold at $V_\pi$
of the diffeomorphism $\varphi_A:\Flag_{n,k}\to\Flag_{n,k}$ and of the
gradient flow of the function $Q_{A,b}:\Flag_{n,k}\to\R$, for any matrix
$A\in H_n(\Escr)$ with simple spectrum.

\bigskip

\begin{lema}
Given $\pi, \pi'\in\Sperm_{n,k}$, if $\pi$ precedes $\pi'$ then 
$H(\pi')= H(\pi)+1$.
\end{lema}

\bigskip

\begin{coro}
Given $\pi, \pi'\in\Sperm_{n,k}$, if $\pi\prec \pi'$ then $H(\pi) <H(\pi')$.
\end{coro}

\bigskip

\begin{coro}
The oriented graph  $\Sperm_{n,k}$ has no cycles.
Therefore,  $\preceq$ is a partial order on $\Sperm_{n,k}$.
\end{coro}

\begin{figure}[h]
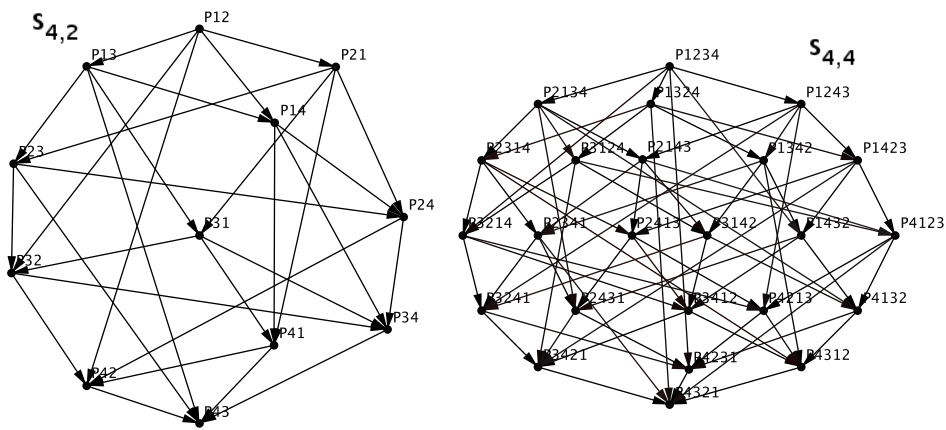

\figura{Sperm}{.5}
\caption{The graphs $\Sperm_{4,2}$ and $\Sperm_{4,4}$}
\end{figure}

We now turn our attention to the symplectic case.
Let ${\rm R}(\pi)=\{\pi_1,\ldots, \pi_k\}$, and 
$\widehat{{\rm C}}(\pi) := \{1,\ldots, 2n\}-\left( {\rm R}(\pi)\cup \overline{ {\rm R}(\pi)} \right)$.
Consider on $\Z_{2n}=\{1,2,\ldots, 2n\}$ the following order
\[ 1 \,\triangleleft\, 2 \,\triangleleft\,\cdots \,\triangleleft\, n \;\triangleleft\; 2n \,
\triangleleft\, 2n-1
\,\triangleleft\, \cdots \,\triangleleft\, n+1 \;.\]

\bigskip

\begin{defi}\label{SJ:leads:to}
Given two permutations $\pi,\pi'\in\SJ_{n,k}$, we say that 
$\pi$ {\em leads to} $\pi'$, and write $\pi \leadsto \pi'$, \, iff\, 
one of the following cases occurs:
\begin{enumerate}
\item[(1)] there is $1\leq i\leq k$ such that
$\pi_i \triangleleft \overline{\pi_i}$ and 
$\pi'$ is obtained from $\pi$ replacing
$\pi_i$ by $\overline{\pi_i}$;

\item[(2)]  there are $1\leq i\leq k$ and $j\in \widehat{{\rm C}}(\pi)$
such that $\pi_i \triangleleft j$ and  $\pi'$ is obtained from $\pi$ replacing
$\pi_i$ by $j$;

\item[(3)]  there are $1\leq i<j \leq k$
such that $\pi_i \triangleleft \pi_j$ and
$\pi'$ is obtained from $\pi$ switching
$\pi_i$ with $\pi_j$;

\item[(4)]  there are $1\leq i<j \leq k$
such that $\pi_i \triangleleft \overline{\pi_j}$   and
$\pi'$ is obtained from $\pi$ replacing
$\pi_i$ by $\overline{\pi_j}$ and 
$\pi_j$ by  $\overline{\pi_i}$.

\end{enumerate}
\end{defi}

\begin{defi}\label{SJ:preceeds:to}
Given two isotropic permutations $\pi,\pi'\in\SJ_{n,k}$, we say that 
$\pi$ {\em precedes } $\pi'$  \, iff\, any
one of the following situations  occurs:
\begin{enumerate}
\item[(1)] there is $1\leq i\leq k$ such that
$\pi_i$ and $\overline{\pi_i}$ are consecutive elements for $\triangleleft$
(which implies $\pi_i=n$ and $\overline{\pi_i}=2n$),  and 
$\pi'$ is obtained from $\pi$ replacing
$\pi_i$ by $\overline{\pi_i}$; 

\item[(2)]  there are $1\leq i\leq k$ and $j\in {\rm C}(\pi)$
such that $\pi_i$ and $j$ are consecutive elements for $\triangleleft$,
and  $\pi'$ is obtained from $\pi$ replacing $\pi_i$ by $j$;

\item[(3)]  there are $1\leq i<j \leq k$
such that $\pi_i \triangleleft \pi_j$,
no $s\in \{1,\ldots, k\}$ exists with $\pi_i\triangleleft \pi_s \triangleleft \pi_j$,
and $\pi'$ is obtained from $\pi$ exchanging
$\pi_i$ with $\pi_j$;

\item[(4)]  there are $1\leq i<j \leq k$
such that $\pi_i$ and $\overline{\pi_j}$ are consecutive elements for $\triangleleft$,
and $\pi'$ is obtained from $\pi$ replacing
$\pi_i$ by $\overline{\pi_j}$ and 
$\pi_j$ by  $\overline{\pi_i}$.
\end{enumerate}
\end{defi}

\bigskip

Comparing the definition of the  "leads-to" relation with proposition~\ref{SJnk:linking}
we see that given $\pi,\pi'\in\SJ_{n,k}$,
$\pi$ and $\pi'$ are linked\, iff\, either $\pi$ leads to $\pi'$ or else
 $\pi'$ leads to $\pi$.
In other words,   the "leads-to" relation orients each and every edge of $\SJ_{n,k}$,
making it an oriented graph, with an oriented edge from $\pi$ to $\pi'$ whenever $\pi \leadsto\pi'$.
The "precede" relation determines a subgraph of $\SJ_{n,k}$ with an 
an oriented edge from $\pi$ to $\pi'$ whenever $\pi$ precedes $\pi'$.

\begin{lema}
Given $\pi, \pi'\in\SJ_{n,k}$, if $\pi$ leads to $\pi'$ then 
there is a sequence of permutations
 $\pi^{(0)}, \pi^{(1)}, \ldots, \pi^{(m)}$ in $\Sperm_{n,k}$
such that
\begin{enumerate}
\item $\pi=\pi^{(0)}$ and $\pi'=\pi^{(m)}$,
\item $\pi^{(i)}$ precedes $\pi^{(i+1)}$, for $i=0,1,\ldots, m-1$.
\end{enumerate}
\end{lema}

This lemma proves that the transitive closures of  relations
"leads-to" and "precede" do coincide. We shall denote this
transitive closure by $\prec$,
and write $\pi \preceq \pi'$ \, iff\, $\pi=\pi'$ or $\pi\prec\pi'$.

\bigskip

Given a isotropic permutation $\pi\in\SJ_{n,k}$, 
consider the set $\Sigma_{n,k}(\pi)$ of all pairs $(i,j)$
satisfying either one of the following conditions:
\begin{enumerate}
\item[(1)] $1\leq i\leq k$,  $j=i+n$ and \, $\pi_i \triangleright \pi_{i+n}$;
\item[(2)] $1\leq i  \leq k$, $j\notin \widehat{{\rm C}}(\pi)$  and \, $\pi_i \triangleright j$;
\item[(3)] $1\leq i<j \leq k$ and \, $\pi_i \triangleright \pi_j$;
\item[(4)] $1\leq i < j \leq k$  and \, $\pi_i \triangleright \pi_{j+n}$\;,
\end{enumerate}
and define the function  $H^{sp}:\SJ_{n,k}\to\N$ by\,
$ H^{sp}(\pi):=\# \Sigma_{n,k}(\pi)$.
As in the general context, $H^{sp}(\pi)$ counts the number of oriented edges of
$\SJ_{n,k}$ which arrive at the vertex $\pi$.
We shall see (c.f. propositions~\ref{sympl:varphiA:Jacobian} and~\ref{sympl:QA:Hessian}) that
$H^{sp}(\pi)$ is exactly the dimension of the stable manifold at $V_\pi$
of the diffeomorphism $\varphi_A:\Fliso_{n,k}\to\Fliso_{n,k}$ and of the
gradient flow of the function $Q_{A,b}:\Fliso_{n,k}\to\R$, for any matrix
$A\in H_n(\Escr)$ with simple spectrum.

\bigskip

\begin{lema}
Given $\pi, \pi'\in\SJ_{n,k}$, if $\pi$ precedes $\pi'$ then 
$H^{sp}(\pi')= H^{sp}(\pi)+1$.
\end{lema}

\bigskip

\begin{coro}
Given $\pi, \pi'\in\SJ_{n,k}$, if $\pi\prec \pi'$ then $H^{sp}(\pi) <H^{sp}(\pi')$.
\end{coro}

\bigskip 

\begin{coro}
The oriented graph  $\SJ_{n,k}$ has no cycles.
Therefore  $\preceq$ is a partial order on  $\SJ_{n,k}$.
\end{coro}

\begin{figure}[h]
\figura{sp2}{.8}
\caption{The graph  $\SJ_{2,2}$ }
\end{figure}

The partial ordered sets $(\Sperm_{n,k},\preceq)$ and $(\SJ_{n,k},\preceq)$
are not lattices for $k>1$.
For instance, $\pi=(4,1)$ and $\pi'=(2,3)$ do not have a greatest lower bound (g.l.b.),
nor a least upper bound (l.u.b.), neither in $\Sperm_{4,2}$ nor in $\SJ_{2,2}$.
However, many subsets of $\Sperm_{n,k}$ and  $\SJ_{n,k}$ have g.l.b. and l.u.b.

\begin{prop}
Given a tree $\underline{P}\in\Kscr_{n,k}$, resp. an isotropic tree $\underline{P}\in\Kscr^{sp}_{n,k}$,
the set of all permutations $\pi\in\Sperm_{n,k}$, resp. $\pi\in\SJ_{n,k}$, such that
$\underline{\pi}\leq \underline{P}$ has a g.l.b. and l.u.b. in 
$(\Sperm_{n,k},\preceq)$, resp. $(\SJ_{n,k},\preceq)$.
\end{prop}

\dem
The proof goes by induction in $k$. 
Let $\Ascr_{\underline{P}}:= \{\,\pi\in\Sperm_{n,k}\,:\,  \underline{\pi}\leq \underline{P}\,\}$.
For $k=1$  $\underline{P}=(P_1)$, 
and the permutations  $\pi= {\rm g.l.b.}\,\Ascr_{\underline{P}}$ and
$\pi'= {\rm l.u.b.}\,\Ascr_{\underline{P}}$ are simply defined by
$\pi_1:= \min P_1$ and $\pi_1':= \max P_1$. In the symplectic case we must use the order
$\triangleleft$. Assume now that the statement of this proposition holds for $k-1$,
and take  $\underline{P}\in\Kscr_{n,k}$. By induction hypothesis there exist
$\pi= {\rm g.l.b.}\,\Ascr_{\pi_{k-1}(\underline{P})}$ and
$\pi'= {\rm l.u.b.}\,\Ascr_{\pi_{k-1}(\underline{P})}$.
We define 
$\pi_k:= \min P_k-{\rm R}(\pi)$ and $\pi_k':= \max P_k-{\rm R}(\pi)$, and set
$\pi:=(\pi_1,\ldots, \pi_{k-1},\pi_k)$ and $\pi':=(\pi_1',\ldots, \pi_{k-1}',\pi_k')$.
Then it is easily checked that
$\pi= {\rm g.l.b.}\,\Ascr_{\underline{P}}$ and
$\pi'= {\rm l.u.b.}\,\Ascr_{\underline{P}}$.
In the symplectic case we proceed as above but taking 
$\pi_k:= \min P_k-({\rm R}(\pi)\cup \overline{{\rm R}(\pi)})$ and
$\pi_k':= \max P_k-({\rm R}(\pi)\cup \overline{{\rm R}(\pi)})$.
\cqd

We shall use the following notation
\begin{align*}
\wedge \underline{P} &= {\rm g.l.b.} \{\,\pi\in\Sperm_{n,k}\,:\, \underline{\pi}\leq \underline{P}\,\} \;, \\
\vee \underline{P}&= {\rm l.u.b.} \{\,\pi\in\Sperm_{n,k}\,:\, \underline{\pi}\leq \underline{P}\,\} \;.
\end{align*}

\begin{prop} Given $A\in H_n(\Escr)$ with simple spectrum, $b\in\R^k_+$
with $b_k > \ldots > b_1 >0$, and a tree $\underline{P}\in\Kscr_{n,k}$, resp. an isotropic tree $\underline{P}\in\Kscr^{sp}_{n,k}$,
let $\pi = \wedge \underline{P}$ and $\pi'= \vee \underline{P}$. Then
denoting by $\mbox{leb}_{\underline{P}}$ the  Lebesgue
measure on $\FlStr{P}$,
\begin{enumerate}
\item[(a)] $V_\pi$ is the unique attractive fixed point for both $\varphi_A:\Flag_{n,k}\to\Flag_{n,k}$
and the gradient flow of $Q_{A,b}:\Flag_{n,k}\to\R$ over the invariant stratum $\FlStr{P}$.
Moreover, for both these dynamical systems,
 $V_\pi$ is the $\omega$-limit of $\mbox{leb}_{\underline{P}}$-almost every flag $V\in\FlStr{P}$;
\item[(b)] $V_{\pi'}$ is the unique repeller fixed point for both $\varphi_A:\Flag_{n,k}\to\Flag_{n,k}$
and the gradient flow of $Q_{A,b}:\Flag_{n,k}\to\R$ over the invariant stratum $\FlStr{P}$.
Moreover, for both these dynamical systems,
 $V_{\pi'}$ is the $\alpha$-limit of  $\mbox{leb}_{\underline{P}}$-almost every flag $V\in\FlStr{P}$.
\end{enumerate}
\end{prop}

\bigskip
\bigskip

\section{Morse Functions}

In this final section we prove that $Q_{A,b}:\Flag_{n,k}\to\R$ are
$\Z_2$-perfect Morse functions.

\bigskip

\begin{prop}[Fixed Points]
Given $A\in H_n$ with simple spectrum and $V\in\Flag_{n,k}$, resp. $V\in\Fliso_{n,k}$, \, $V=\varphi_A(V)$\, iff\,
$V=V_\pi$ 
for some permutation $\pi\in\Sperm_{n,k}$, resp.  isotropic permutation  $\pi\in\SJ_{n,k}$.
\end{prop}

\dem
Follows from lemma~\ref{fixed:point:char}.
\cqd

In particular $\varphi_A:\Flag_{n,k}\to\Flag_{n,k}$ has $n!/(n-k)!$ fixed points,
while $\varphi_A:\Fliso_{n,k}\to\Fliso_{n,k}$ has $(2n)(2n-2)\cdots (2n-2k+2)$ fixed points.

\bigskip 

\begin{prop}[Critical Points]
Given $A\in H_n$ with simple spectrum and $V\in\Flag_{n,k}$, resp. $V\in\Fliso_{n,k}$,
 \, $V$ is a critical point of $Q_{A,b}$ \, iff\,
$V=V_\pi$ 
for some permutation $\pi\in\Sperm_{n,k}$, resp. isotropic permutation $\pi\in\SJ_{n,k}$.
\end{prop}

\dem
Follows from lemma~\ref{crit:point:char}.
\cqd

In particular $Q_{A,b}:\Flag_{n,k}\to\R$ has $n!/(n-k)!$ critical points,
while $Q_{A,b}:\Fliso_{n,k}\to\R$ has $(2n)(2n-2)\cdots (2n-2k+2)$ critical points.

\bigskip

We use abstract lemmas to compute the eigenvalues at the fixed and critical points.
Assume $E\subset\R^n$ is an $A$-invariant $2$-dimensional plane,
and define $\Su^1(E)=\{\, x\in E\,:\, \nrm{x}=1\,\}$,
and  $\phi_A:\Su^1(E)\to\Su^1(E)$ by $\phi_A(x)= A\,x/\nrm{ A\,x }$.
Denote by $\lambda_i$ ($i=1,2$) the eigenvalues of $A$ on $E$, and by
$v_i$ ($i=1,2$) the corresponding eigenvectors of $A$, i.e.,
$A\,v_i=\lambda_i\, v_i$, with $v_i\in E$ and  $\nrm{v_i}=1$.

\begin{lema} \label{A:eigen:abstract}
If the eigenvalues are distinct, $\lambda_1\neq \lambda_2$,
then $\phi_A:\Su^1(E)\to\Su^1(E)$ has four fixed points:\, $\pm\,v_1$ and $\pm\,v_2$.
The eigenvalues of $\phi_A$  are:\,
$\lambda_2/\lambda_1$  at $\pm \,v_1$,
and  $\lambda_1/\lambda_2$ at $\pm\,v_2$.
\end{lema}

\bigskip

Consider the function $q:\Su^1(E)\to\R$,\, $q(x)=\frac{c}{2}\,\nrm{A\,x}^2$.

\begin{lema} \label{A:hess:Su:abs}
If the eigenvalues are distinct, $\lambda_1\neq \lambda_2$,
then $q:\Su^1(E)\to\R$ has four critical points:\, $\pm\,v_1$ and $\pm\,v_2$.
The eigenvalues of the Hessian of $q$  are:\,
$c\,((\lambda_2)^2-(\lambda_1)^2)$  at $\pm \,v_1$,
and  $c\,((\lambda_1)^2-(\lambda_2)^2)$  at $\pm\,v_2$.
\end{lema}

\bigskip

Consider now the function $q:\OO(2,E)\to\R$,\, $q(X)=\frac{1}{2}\,\nrm{A\,X\,D_c}^2$,
where $\OO(2,E)$ is the space of orthogonal linear maps $X:\R^2\to E$ onto the $2$-plane $E$,
and $D_c=\matriz{c_1}{0}{0}{c_2}$ with $c_1>c_2>0$.

\begin{lema} \label{A:hess:OO:abs}
If the eigenvalues are distinct, $\lambda_1\neq \lambda_2$,
then $q:\OO(2,E)\to\R$ has four critical points:\, $\pm\,v_1$ and $\pm\,v_2$.
The eigenvalues of the Hessian of $q$  are:\,
$c\,((\lambda_2)^2-(\lambda_1)^2)$  at $\pm \,v_1$,
and  $c\,((\lambda_1)^2-(\lambda_2)^2)$  at $\pm\,v_2$,
where $c= (c_1)^2-(c_2)^2>0$.
\end{lema}

\bigskip

Given two linked permutations $\pi,\pi'\in\Sperm_{n,k}$ there is a unique
one-dimensional stratum containing $\underline{\pi}$ and $\underline{\pi'}$.
We shall denote it by $\underline{E}(\pi,\pi')$.
By proposition~\ref{Snk:linking} in case (1) $E_i(\pi,\pi')= \{\pi_i,j\}$ and 
$E_s(\pi,\pi')=\{\pi_s\}$ for $s\neq i$, while in case (2) $E_i(\pi,\pi')=E_j(\pi,\pi')=\{\pi_i,\pi_j\}$ and 
$E_s(\pi,\pi')=\{\pi_s\}$ for $s\notin\{i,j\}$.
 $\underline{\pi}$ and $\underline{\pi'}$ are the unique zero dimensional strata contained in the circle
stratum $\underline{E}(\pi,\pi')$.
These correspond to attractive and repelling fixed points on $\underline{E}(\pi,\pi')$ for both
$\varphi_A$ and the gradient flow of $Q_{A,b}$.
In the following propositions $\lambda_i$ will always denote the eigenvalue
of the given matrix $A\in H_n$ associated with the eigendirection $E_i$.

\bigskip

Given a permutation $\pi\in\Sperm_{n,k}$ we call
 {\em eigenvector basis} at $V_\pi \in \Flag_{n,k}$ to any basis of $T_{V_\pi} \Flag_{n,k}$
whose vectors are tangent to the one dimensional strata ${\rm S}_{\underline{E}(\pi,\pi')}$
containing $\underline{\pi}=\{ V_\pi\}$.
A similar definition is given for every isotropic permutation $\pi\in\SJ_{n,k}$.

\begin{prop}[Jacobian at the Fixed Points]  \label{varphiA:Jacobian}
Given $\pi\in \Sperm_{n,k}$, the Jacobian matrix $D (\varphi_A)_{V_\pi}$ 
is a diagonal matrix w.r.t. any eigenvector basis  at $V_\pi$ with the
following eigenvalues:
\begin{enumerate}
\item \; $\lambda_{j}/\lambda_{\pi_i}$ for every pair $1\leq i \leq k$, $j\in\{1,\ldots, n\}-{\rm R}(\pi)$,
\item \; $\lambda_{\pi_j}/\lambda_{\pi_i}$ for every pair $1\leq i<j \leq k$.
\end{enumerate}
The Jacobian $D (\varphi_A)_{V_\pi}$ has \,$H(\pi)$\, eigenvalues $\lambda$ with $\lambda >1$  and\,
$d -H(\pi)$\, eigenvalues  $\lambda$ with $0\leq \lambda<1$,
where $d= k\,(2n-k-1)/2 = \dim \OO_{n,k}$.
\end{prop}

\dem
We have to compute the eigenvalue of  the Jacobian matrix $D (\varphi_A)_{V_\pi}$ 
along the each eigen-direction tangent to an invariant one-dimensional stratum
${\rm S}_{\underline{E}(\pi,\pi')}$ with $\pi'$ linked to $\pi$.
There are two cases according to proposition~\ref{Snk:linking}.
Consider the $A$-invariant $2$-plane:\, $E=E_{\{\pi_i, j\}}$ in case (1),
and \, $E=E_{\{\pi_i,\pi_j\}}$ in case (2),
and the projection $p:\widetilde{{\rm S}}_{\underline{E}(\pi,\pi')} \to \Su^1(E)$
defined by $p(X)=(X)_{i}$ in both cases. We can easily check that
$p\circ \varphi_A = \phi_A \circ p$. This shows that the eigenvalue of
$\varphi_A$ at the fixed point $X_\pi$ (along $\widetilde{{\rm S}}_{\underline{E}(\pi,\pi')}$) coincides
with the eigenvalue of $\phi_A$ at $p(X_\pi)$, and the claimed eigenvalues follow from
lemma~\ref{A:eigen:abstract}.
\cqd

\begin{prop} [Hessian at the Critical Points]  \label{QA:Hessian}
Given $\pi\in \Sperm_{n,k}$, the Hessian matrix ${\rm Hess}(Q_{A,b})_{V_\pi}$ 
is a diagonal matrix w.r.t. any eigenvector basis  at $V_\pi$ with the
following eigenvalues:
\begin{enumerate}
\item \; $2\,k^{-1}\, b_i^2\, (\,(\lambda_{j})^2 - (\lambda_{\pi_i})^2\,)$ 
for every pair $1\leq i \leq k$, $j\in\{1,\ldots, n\}-{\rm R}(\pi)$,
\item \; $2\,k^{-1}\,(b_i^2-b_j^2)\,(\,(\lambda_{\pi_j})^2 - (\lambda_{\pi_i}) ^2\,)$ 
 for every pair $1\leq i<j \leq k$.
\end{enumerate}
The Hessian ${\rm Hess}(Q_{A,b})_{V_\pi}$  has $H(\pi)$ positive eigenvalues and
$d -H(\pi)$ negative eigenvalues,
where $d= k\,(2n-k-1)/2 = \dim \OO_{n,k}$.
\end{prop}

\dem
We have to compute the eigenvalue of  the Hessian matrix ${\rm Hess}(Q_{A,b})_{V_\pi}$  
along the each eigen-direction tangent to an invariant one-dimensional stratum
${\rm S}_{\underline{E}(\pi,\pi')}$ with $\pi'$ linked to $\pi$.
There are two cases according to proposition~\ref{Snk:linking}.
Let $E$ be the $2$-plane defined in the proof of proposition~\ref{varphiA:Jacobian}.
In case (2) consider the projection $p:\widetilde{{\rm S}}_{\underline{E}(\pi,\pi')} \to \OO(2,E)$
defined by 
\[ p(X)(u_1,u_2) =  (X)_{i}\, u_1 +  (X)_{j}\, u_2 \,.\]
We can easily check that for $X$ over $\widetilde{{\rm S}}_{\underline{E}(\pi,\pi')}$
\begin{align*}
Q_{A,b}(X) &= \frac{1}{k}\,\sum_{i=1}^k b_i^2\,\nrm{ (A\,X)_i}^2 \\
 &= \text{ const} +  \frac{b_i^2}{k}\,\nrm{ (A\,X)_i}^2 + 
\frac{b_j^2}{k}\,\nrm{ (A\,X)_j}^2 \\
&= \text{ const} + \frac{1}{2}\,\nrm{ A\,p(X)\,D_c }^2 = 
\text{ const} + (q\circ p)(X)
\end{align*}
where $c_i:= 2\,b_i^2/k > c_j:= 2\,b_j^2/k >0$,
and $q(Y)=\frac{1}{2}\,\nrm{A\,Y\,D_c}^2$, for $Y\in\OO(2,E)$.
This shows that the eigenvalue of
${\rm Hess} (Q_{A,b})_X$ along $\widetilde{{\rm S}}_{\underline{E}(\pi,\pi')}$ coincides
with the corresponding eigenvalue of ${\rm Hess}(q)_{p(X)}$, 
and the claimed eigenvalues follow from
lemma~\ref{A:hess:OO:abs}.

\smallskip

In case (1) consider the projection $p:\widetilde{{\rm S}}_{\underline{E}(\pi,\pi')} \to \Su^1(E)$
defined by $p(X) =  (X)_{i}$.
We can check that for $X$ over $\widetilde{{\rm S}}_{\underline{E}(\pi,\pi')}$
\begin{align*}
Q_{A,b}(X) &= \frac{1}{k}\,\sum_{i=1}^k b_i^2\,\nrm{ (A\,X)_i}^2 \\
 &= \text{ const} +  \frac{b_i^2}{k}\,\nrm{ (A\,X)_i}^2  \\
&= \text{ const} + \frac{c}{2}\,\nrm{ A\,p(X) }^2 = 
\text{ const} + (q\circ p)(X)
\end{align*}
where $c:= 2\,b_i^2/k > 0$,
and $q(y)=\frac{c}{2}\,\nrm{A\,y}^2$, for $y\in\Su^1(E)$.
This shows that the eigenvalue of
${\rm Hess} (Q_{A,b})_X$ along $\widetilde{{\rm S}}_{\underline{E}(\pi,\pi')}$ coincides
with the corresponding eigenvalue of ${\rm Hess}(q)_{p(X)}$, 
and the claimed eigenvalues follow from
lemma~\ref{A:hess:Su:abs}.
\cqd

\bigskip

In the symplectic case, given linked isotropic permutations $\pi,\pi'\in\SJ_{n,k}$ 
we have four cases by proposition~\ref{SJnk:linking}.
In case (1) $E_i(\pi,\pi')= \{\pi_i,\overline{\pi_i}\}$ and 
$E_s(\pi,\pi')=\{\pi_s\}$ for $s\neq i$.
In case (2) $E_i(\pi,\pi')= \{\pi_i,j\}$ and 
$E_s(\pi,\pi')=\{\pi_s\}$ for $s\neq i$.
In case (3) $E_i(\pi,\pi')=E_j(\pi,\pi')=\{\pi_i,\pi_j\}$ and 
$E_s(\pi,\pi')=\{\pi_s\}$ for $s\notin\{i,j\}$.
Finally, in case (4) $E_i(\pi,\pi')=\{\pi_i,\overline{\pi_j} \}$,
$E_j(\pi,\pi')=\{\overline{\pi_i},\pi_j\}$  and 
$E_s(\pi,\pi')=\{\pi_s\}$ for $s\notin\{i,j\}$.

\bigskip

\begin{prop}[Symplectic Jacobian at the Fixed Points]  \label{sympl:varphiA:Jacobian}
Given $\pi\in \SJ_{n,k}$, the Jacobian matrix $D (\varphi_A)_{V_\pi}$ 
of $\varphi_A:\Fliso_{n,k}\to\Fliso_{n,k}$ is a diagonal matrix w.r.t. any eigenvector basis  at $V_\pi$ with the
following eigenvalues:
\begin{enumerate}
\item \; $\lambda_{\overline{\pi_i}}/\lambda_{\pi_i}$ for every $1\leq i \leq k$,
\item \; $\lambda_{j}/\lambda_{\pi_i}$ for every pair $1\leq i \leq k$, $j\in\{1,\ldots, 2n\}-\widehat{{\rm R}}(\pi)$,
\item \; $\lambda_{\pi_j}/\lambda_{\pi_i}$ for every pair $1\leq i<j \leq k$,
\item \; $\lambda_{\overline{\pi_j}}/\lambda_{\pi_i}$ for every pair $1\leq i<j \leq k$.
\end{enumerate}
The Jacobian $D (\varphi_A)_{V_\pi}$ has \,$H^{sp}(\pi)$\, eigenvalues $\lambda$ with $\lambda >1$  and\,
$d -H^{sp}(\pi)$\, eigenvalues  $\lambda$ with $0\leq \lambda<1$,
where $d= k\,(2n-k) = \dim \Osp_{n,k}$.
\end{prop}

\bigskip

\begin{prop} [Symplectic Hessian at the Critical Points]  \label{sympl:QA:Hessian}
Given $\pi\in \SJ_{n,k}$, the Hessian matrix ${\rm Hess}(Q_{A,b})_{V_\pi}$ 
of $Q_{A,b}:\Fliso_{n,k}\to\R$  is a diagonal matrix w.r.t. any eigenvector basis  at $V_\pi$ with the
following eigenvalues:
\begin{enumerate}
\item \; $2\,k^{-1}\, b_i^2\, (\,(\lambda_{\overline{\pi_i}})^2 - (\lambda_{\pi_i})^2\,)$  
for every  $1\leq i\leq k$,
\item \; $2\,k^{-1}\, b_i^2\, (\,(\lambda_{j})^2 - (\lambda_{\pi_i})^2\,)$  
for every pair $1\leq i \leq k$, $j\in\{1,\ldots, 2n\}-\widehat{{\rm R}}(\pi)$,
\item \; $2\,k^{-1}\,(b_i^2-b_j^2)\,(\,(\lambda_{\pi_j})^2 - (\lambda_{\pi_i}) ^2\,)$ for every pair
$1\leq i < j \leq k$,
\item \; $2\,k^{-1}\,(b_i^2-b_j^2)\,(\,(\lambda_{\overline{\pi_j}})^2 - (\lambda_{\pi_i}) ^2\,)$  for every pair 
$1\leq i < j \leq k$.
\end{enumerate}
The Hessian ${\rm Hess}(Q_{A,b})_{V_\pi}$  has $H^{sp}(\pi)$ positive eigenvalues and
$d -H^{sp}(\pi)$ negative eigenvalues,
where $d= k\,(2n-k) = \dim \Osp_{n,k}$.
\end{prop}

\bigskip

\begin{prop} \label{Morse:function}
Given $A\in H_n(\Escr)$ with simple spectrum and  $b=(b_1,\ldots, b_k)\in\R^k$ with $b_1>b_2>\ldots >b_k>0$, 
the functions $Q_{A,b}:\Flag_{n,k}\to\R$ and
$Q_{A,b}:\Fliso_{n,k}\to\R$
are Morse functions.
\end{prop}

\dem By proposition~\ref{crit:point:char} the critical points of $Q_{A,b}$
are exactly the points $V_\pi$ with $\pi\in\Sperm_{n,k}$, or $\pi\in\SJ_{n,k}$.
Then  propositions~\ref{QA:Hessian} and~\ref{sympl:QA:Hessian} show that every critical point
is non-degenerate. Therefore, $Q_{A,b}$ is a Morse function.

\cqd

\begin{prop} \label{Poincare:polys}
The Poincar\'{e} polynomial of $\Flag_{n,k}$ with coefficients in $\Z_2$ is
\begin{equation}\label{Poincare:Flag}
\Pscr_{n,k}(t)= \frac{(1-t^n)\,(1-t^{n-1})\,\cdots\, (1-t^{n-k+1}) }{(1-t)^k}\;.
\end{equation}
The Poincar\'{e} polynomial of $\Flag^{sp}_{n,k}$ with coefficients in $\Z_2$ is
\begin{equation}\label{Poincare:Fliso}
\Pscr^{sp}_{n,k}(t)= \frac{(1-t^{2n})\,(1-t^{2n-2})\,\cdots\, (1-t^{2n-2k+2}) }{(1-t)^{k}}\;.
\end{equation}
\end{prop}

\dem
See~\cite{B} or~\cite{MN}.
\cqd

\bigskip

We define the following coefficients:
Given $1\leq k\leq n$,\, $0\leq s \leq d= k\,(2n-k-1)/2$,
\begin{equation}
 \pol{n}{k}{s} := \#\{\,(s_1,\ldots, s_k)\,:\; s_1+\cdots +s_k=s \; \text{ and } \;
 0\leq s_i\leq n-i, \; \forall i \;\} \;.
\end{equation}
Analogously, given $1\leq k\leq n$,\, $0\leq s \leq d= k\,(2n-k)$,
\begin{equation}
 \polsp{n}{k}{s} := \#\{\,(s_1,\ldots, s_k)\,:\; s_1+\cdots +s_k=s \; \text{ and } \;
 0\leq s_i\leq 2n-2i, \; \forall i \;\} \;.
\end{equation}

\bigskip

\begin{lema} \label{pnki}
The Poincar\'{e} polynomials given in~(\ref{Poincare:Flag}) and~(\ref{Poincare:Fliso})
are respectively equal to
$$ \Pscr_{n,k}(t)= \sum_{i=0}^{k(2n-k-1)/2} \pol{n}{k}{i}\, t^i \quad \text{ and }\quad
\Pscr_{n,k}^{sp}(t)= \sum_{i=0}^{k(2n-k)} \polsp{n}{k}{i}\, t^i \;. $$
\end{lema}

\dem Follows easily from the relations
\[ \Pscr_{n,k}(t)=\prod_{i=1}^{k} \frac{(1-t^{n-i+1}) }{(1-t)}  =
\prod_{i=1}^{k} (1+t+t^2+\cdots + t^{n-i}) \]
\[ \Pscr^{sp}_{n,k}(t)=\prod_{i=1}^{k} \frac{(1-t^{2n-2i+2}) }{(1-t)}  =
\prod_{i=1}^{k} (1+t+t^2+\cdots + t^{2n-2i+1}) \]
\cqd

\bigskip

\begin{prop} The Morse polynomial of $Q_{A,b}:\Flag_{n,k}\to\R$
is equal to $\Pscr_{n,k}(t)$. In particular, this function is
$\Z_2$-perfect.
\end{prop}

\dem
We claim that for every $1\leq k\leq n$  and $0\leq i\leq d =k\,(2n-k-1)/2$,
\[  p^{n,k}_i = \#\,\{\, \pi\in\Sperm_{n,k}\,:\, H(\pi)=i \;\}  \;.\]
Thus, by proposition~\ref{QA:Hessian} and lemma~\ref{pnki},  
$\Pscr_{n,k}(t)$ is the Morse polynomial of 
$Q_{A,b}$. Let us now prove the claim.
For $1\leq k\leq n$,\, let
$$\Gamma_{n,k}= \{\,(s_1,\ldots, s_k)\,:\; 0\leq s_i\leq n-i, \; \forall\, 1\leq i\leq k \;\} \;.$$
We shall define a mapping $\phi:\Gamma_{n,k}\to\Sperm_{n,k}$  such that
$H\circ \phi(s_1,\ldots, s_k) = s_1+\cdots +s_k$, for every
$(s_1,\ldots, s_k)\in \Gamma_{n,k}$, and prove that $\phi$ is one-to-one.
Given $\pi\in\Sperm_{n,k}$
write  ${\rm C}_i(\pi):= \{1,\ldots, n\}-\{\pi_1,\ldots, \pi_{i-1}\}$.
Define $\psi:\Sperm_{n,k}\to \Gamma_{n,k}$ by
$\psi(\pi):=(s_1,\ldots, s_k)$, where 
$s_i:= \#\{ \, j\in {\rm C}_i(\pi) \,:\, \pi_i>j\,\}$.
By definition we have $0\leq s_i \leq n-i$, so that $\psi(\pi)= (s_1,\ldots, s_k)\in \Gamma_{n,k}$.
Conversely, 
$\pi = \phi(s_1,\ldots, s_k)$ is recursively defined as follows:
$\pi_i$ is taken to be the 
$(s_i+1)^{\text{th}}$-element of ${\rm C}_i(\pi)$.
In this way, there are exactly $s_i$ indices $j=i+1,\ldots, n$ such that
$j\in {\rm C}_i(\pi)$ and $\pi_i > j$, and this shows that $\phi$ is the inverse of $\psi$.
\cqd

\bigskip

\begin{prop} The Morse polynomial of $Q_{A,b}:\Fliso_{n,k}\to\R$
is equal to $\Pscr^{sp}_{n,k}(t)$. In particular, this function is
$\Z_2$-perfect.
\end{prop}

\dem
We claim that for every $1\leq k\leq n$  and $0\leq i\leq d =k\,(2n-k)$,
\[  \polsp{n}{k}{i} = \#\,\{\, \pi\in\Sperm_{n,k}\,:\, H^{sp}(\pi)=i \;\}  \;.\]
Thus, by proposition~\ref{sympl:QA:Hessian} and lemma~\ref{pnki},  
$\Pscr_{n,k}(t)$ is the Morse polynomial of 
$Q_{A,b}$. Let us now prove the claim.
For $1\leq k\leq n$,\, let
$$\Gamma_{n,k}^{sp}= \{\,(s_1,\ldots, s_k)\,:\; 0\leq s_i\leq 2n-2i, \; \forall\, 1\leq i\leq k \;\} \;.$$
We shall define a mapping $\phi:\Gamma_{n,k}^{sp}\to\SJ_{n,k}$  such that
$H^{sp}\circ \phi(s_1,\ldots, s_k) = s_1+\cdots +s_k$, for every
$(s_1,\ldots, s_k)\in \Gamma^{sp}_{n,k}$, and prove that $\phi$ is one-to-one.
Given $\pi\in\SJ_{n,k}$
write  $\widehat{{\rm C}}_i(\pi):= \{1,\ldots, 2n\}-\{\pi_1,\ldots, \pi_{i-1},
\overline{\pi_1},\ldots, \overline{\pi_{i-1}} \}$.
Define $\psi:\SJ_{n,k}\to \Gamma^{sp}_{n,k}$ by
$\psi(\pi):=(s_1,\ldots, s_k)$, where 
$s_i:= \#\{ \, j\in \widehat{{\rm C}}_i(\pi) \,:\, \pi_i>j\,\}$.
By definition we have $0\leq s_i \leq 2n-2i$, so that $\psi(\pi)= (s_1,\ldots, s_k)\in \Gamma^{sp}_{n,k}$.
Conversely, 
$\pi = \phi(s_1,\ldots, s_k)$ is recursively defined as follows:
$\pi_i$ is taken to be the 
$(s_i+1)^{\text{th}}$-element of $\widehat{{\rm C}}_i(\pi)$.
In this way, there are exactly $s_i$ indices $j=i+1,\ldots, n$ such that
$j\in \widehat{{\rm C}}_i(\pi)$ and $\pi_i > j$, and this shows that $\phi$ is the inverse of $\psi$.
\cqd

\bigskip

\thispagestyle{empty}

\end{document}